\documentclass[10pt, letterpaper, twoside]{article} 

\usepackage{amsmath, amsthm, amssymb,bbold,mathrsfs} 
\usepackage{graphicx,hyperref}
\usepackage[toc,page]{appendix}

\newtheorem{thm}{Theorem}[section]
\newtheorem{cor}{Corollary}[section]
\newtheorem{lem}{Lemma}[section]
\newtheorem{prop}{Proposition}[section]

\newtheorem{assumption}{Assumption}[section]

\theoremstyle{remark}
\newtheorem{rem}{Remark}[section]

\theoremstyle{definition}
\newtheorem{example}{Example}[section]

\numberwithin{equation}{section}
 
 \def\k{\theta}
 \def\K{\Theta}
 \def\0{\mathbf{0}}
 \def\b1{\mathbb{1}}
 
 \def\argmin{\mathop{\arg\min}}
 \def\E{\mathbb{E}}
 \def\M{\mathcal{M}}
 \def\U{\mathscr{U}}
 \def\B{\mathscr{B}}
 \def\A{\mathscr{A}}
 \def\P{\mathscr{P}}
 
 \def\1{\mathbf{1}}
 
\def\moveup{\vspace{-0.1cm}}
\def\mysmallskip{\vspace*{0.01in}}

\begin{document} 

 \bigskip
 \title{A Mixed Value and Policy Iteration Method for \\ Stochastic Control with Universally Measurable Policies}
 
\author{Huizhen Yu\thanks{Huizhen Yu was with the Lab.\ for Information and Decision Systems, M.I.T., and she is now with the Department of Computing Science, University of Alberta.
        {\tt\small janey.hzyu@gmail.com}}
\and
Dimitri P. Bertsekas\thanks{Lab.\ for Information and Decision Systems, M.I.T.
        {\tt\small dimitrib@mit.edu}}
}

\date{} 
 \maketitle
 
\begin{abstract}
We consider stochastic optimal control models with Borel spaces and universally measurable policies. For such models the standard policy iteration is known to have difficult measurability issues and cannot be carried out in general. We present a mixed value and policy iteration method that circumvents this difficulty. The method allows the use of stationary policies in computing the optimal cost function, in a manner that resembles policy iteration. It can also be used to address similar difficulties of policy iteration in the context of upper and lower semicontinuous models. We analyze the convergence of the method in infinite horizon total cost problems, for the discounted case where the one-stage costs are bounded, and for the undiscounted case where the one-stage costs are nonpositive or nonnegative.

For undiscounted total cost problems with nonnegative one-stage costs, we also give a new convergence theorem for value iteration, which shows that value iteration converges whenever it is initialized with a function that is above the optimal cost function and yet bounded by a multiple of the optimal cost function. This condition resembles Whittle's bridging condition and is partly motivated by it. The theorem is also partly motivated by a result of Maitra and Sudderth, which showed that value iteration, when initialized with the constant function zero, could require a transfinite number of iterations to converge. We use the new convergence theorem for value iteration to establish the convergence of our mixed value and policy iteration method for the nonnegative cost case.
\end{abstract}

\thispagestyle{myheadings}
\markboth{}{ \large  LIDS REPORT 2905}
\setlength{\unitlength}{1cm}
\begin{picture}(0,0)(0,0)
\put(-0.52,13.6){\fontsize{12}{14} \selectfont \it July 2013; Revised July 2014}
%\put(-0.52,13.1){\fontsize{12}{14} \selectfont \it Last updated December 2014}
\put(-0.52,13.1){\fontsize{12}{14} \selectfont \it Accepted by Math.\ Oper.\ Res.}
\end{picture}

\clearpage
   \tableofcontents
\clearpage

\section{Introduction} \label{sec1}
\markboth{\rm \S \thesection. Introduction}{\rm \S \thesection. Introduction}

 We consider discrete-time stochastic control problems with additive one-stage costs 
 in a general framework that involves Borel state and control spaces and universally measurable policies.
Historically, our framework traces back to the pioneering work on dynamic programming (DP) in Borel spaces by Blackwell \cite{Blk-finite,Blk-discount,Blk-positive,Blk-borel} and Strauch~\cite{Str-negative}, which was developed further, along several directions, through a sequence of subsequent works.
These include: the books by Hinderer~\cite{Hinderer70}, and Dynkin and Yushkevich~\cite{Dynkin79}, which considered a framework based on Borel measurable policies and the notion of almost-sure $\epsilon$-optimality; the work of Maitra~\cite{Maitra68}, Furukawa~\cite{Fur72}, Freedman~\cite{Fre74} and Sch\"{a}l~\cite{Schal75}, as well as Dynkin and Yushkevich~\cite{Dynkin79}, which studied Borel measurable policies and semicontinuous models; the work of Blackwell, Freedman and Orkin~\cite{BFO74}, which introduced a formulation involving analytic sets and analytically measurable policies;
and the work of Shreve and Bertsekas \cite{ShrB78,ShrB79}, and Bertsekas and Shreve~\cite[Part II]{bs}, which considered universally measurable policies.
Further research on alternative frameworks suitable for DP include: Shreve~\cite{Shr78} and Bertsekas and Shreve~\cite[Part II]{bs} on C-sets and limit-measurable policies, Blackwell \cite{Blk-bp} on Borel-programmable functions and Shreve~\cite{Shr81} on Borel-approachable functions.
 We refer to the monograph \cite{bs} and the papers \cite{ShrB78,Shr79} for a discussion of the differences between these frameworks, along with a review of the literature for the early period of the subject. We refer to the books \cite{Puterman94,HL96,HL99, Alt99, MDPhandbook}, \cite[Part III]{Mey08}, and the survey papers \cite{acmdp-survey93,Fein02} for more recent accounts and extensive references about the significant development of the field since then. 
In this paper, we will focus on the universally measurable policies framework of \cite{ShrB78,ShrB79,bs}, and three types of classical infinite horizon total cost problems: the discounted case where the one-stage costs are bounded, and the undiscounted case where the one-stage costs are all nonpositive or all nonnegative.

The early works of Blackwell and Strauch showed that taking Borel measurable policies as the only admissible policies does not lead to desirable results that are comparable with the ones available for problems where measurability is not a concern. In particular, a Borel measurable policy need not exist even when the control constraint set is Borel \cite{Blk-borel}. Moreover,
if we restrict attention to Borel measurable policies, there need not exist an everywhere $\epsilon$-optimal policy even in discounted problems~\cite{Blk-discount}. 
An important step toward a more satisfactory framework was taken by Blackwell, Freedman and Orkin~\cite{BFO74}. Studying finite horizon nonnegative reward problems, they introduced an approach based on analytic sets and semi-analytic functions (a family of functions whose level sets are analytic sets), and obtained optimality results for analytically measurable policies (a larger class of policies that includes Borel measurable ones). Their model still does not admit the existence of everywhere optimal policies or the existence of everywhere $\epsilon$-optimal nonrandomized policies among structured families of policies in general.  Building upon analytic sets and semi-analytic functions as in~\cite{BFO74}, a fuller framework was developed
in Shreve and Bertsekas \cite{ShrB78,ShrB79}, Bertsekas and Shreve~\cite[Part II]{bs}. 
In this framework, the class of admissible policies is enlarged to be the class of universally measurable policies,
structural properties of the optimal cost functions are derived, 
and selection theorems that stem from the Jankov-von Neumann theorem ensure the existence of everywhere $\epsilon$-optimal or optimal policies among structured families of policies (e.g., stationary, Markov or semi-Markov policies), both for finite horizon problems and for infinite horizon 
problems that we consider.

However, with analytically or universally measurable policies, standard policy iteration has measurability-related difficulties, as noted in \cite[p.\ 940]{BFO74} and \cite[p.\ 232]{bs}.
The selection of an admissible measurable policy can fail at the policy improvement step because the cost function of an analytically or universally measurable policy need not have the necessary structure for exact or $\epsilon$-exact selection of an improved policy. This causes the policy iteration procedure to break down.

A similar difficulty occurs in upper and lower semicontinuous models. There the selection of a Borel measurable policy at the policy improvement step may fail because the cost function of the current Borel measurable policy does not have adequate semicontinuity structure.

One of the major purposes of this paper is to provide an approach to circumvent the difficulty just discussed, and to allow stationary policies to be used in computing the optimal cost function, in a manner that resembles policy iteration (even when $\epsilon$-optimal stationary policies do not exist).
We refer to our approach as a mixed value and policy iteration method, as it combines characteristics of both value and policy iteration. Algorithmically, compared to standard policy iteration, the main difference of our method is in the policy evaluation phase: instead of computing the costs of a given policy, it solves exactly or approximately an optimal stopping problem defined by a stationary policy of interest and by a stopping cost that is an estimate of the optimal cost. The stopping costs are then adjusted and the procedure is repeated.
To avoid measurability issues, we exploit the fact that every universally measurable stationary policy has Borel measurable portions (see Prop.\ 3.1(b)), and we define the optimal stopping problems accordingly so that the iterative method just mentioned can operate within the family of functions with the desired semi-analytic structure.
Another critical feature of our approach results from the optimal stopping formulation and its inherent value iteration character: for convergence, it does not require the policies involved to improve successively on one another (this is generally impossible within our context).
This feature allows us to operate the method more flexibly and leads to a variety of algorithms.
As a result, we obtain policy iteration-like algorithms if we choose policies in a way analogous to policy improvement, using the Jankov-von Neumann type of selection theorems.

Similarly, for the lower semicontinuous model, we exploit the fact that Borel measurable policies have continuous portions (Lusin's Theorem; see e.g.,~\cite{Dud02}), and we use it to specialize our method to produce policy iteration-like algorithms that operate within the desired class of lower semicontinuous functions.
For the upper semicontinuous model, a similar, albeit less general, result is derived: with additional assumptions which may hold only for certain problems, we show how to construct policy iteration-like algorithms that operate within upper semicontinuous functions.

We establish the convergence of our method under certain initial conditions for the three types of infinite horizon total cost problems we consider. Our convergence results parallel those for standard value iteration for these problems.
 
For nonnegative cost models, the convergence analysis of our mixed value and policy iteration method relies on another main result of this paper, which is of independent interest. This is a new convergence theorem for value iteration.
It is well-known that for nonnegative cost models, value iteration need not converge to the optimal cost function. Conditions for convergence from below, which involve compactness-type assumptions on the control constraint set, have been given by Bertsekas \cite{b72} for a related special case of minimax reachability problems, by Sch{$\ddot{\mathrm{a}}$l \cite{Schal75} and Bertsekas \cite{b77} for cases where measurability issues are not a concern, and by Bertsekas and Shreve~\cite{bs} for the universally measurable policies framework of this paper. Sufficient conditions have also been studied by Whittle~\cite{Whit79,Whit80}.

Our theorem shows that value iteration converges whenever it is initialized with a function that lies above the optimal cost function and yet is bounded by a multiple of the optimal cost function.
 This condition resembles Whittle's bridging condition~\cite{Whit79,Har80} and is partly motivated by it. Whittle's condition, however, delineates a subset of nonnegative cost models in which value iteration converges when initialized with the constant function zero, whereas our theorem holds without model restrictions. In formulating the theorem, we were also partly motivated by a general convergence result of Maitra and Sudderth~\cite{MS92}, which showed that starting from the constant function zero, value iteration could require a transfinite number of iterations to converge.
We will give two alternative proofs of the new theorem for the convergence of value iteration (in the standard, non-transfinite form), one of which uses, among others, Maitra and Sudderth's result.

We also note that Meyn in his book~\cite[Chap.\ 9]{Mey08} uses systematically the space of functions bounded in absolute value by a multiple of the optimal cost function (or an optimal relative cost function) to analyze average cost and total cost problems. The other proof we give for the new convergence theorem is motivated by some arguments used there (this line of analysis can be extended to total cost problems without sign constraints on one-stage costs \cite{Yu-vi14}).

Using the new convergence theorem for value iteration, we are also able to show that for certain nonnegative cost models (which include countable-spaces problems with finite optimal costs), 
convergence of our mixed value and policy iteration method is maintained if the optimal stopping problems involved are solved approximately by solving associated linear programs.
This result can be contrasted with the fact that nonnegative cost models in general do not admit a linear programming formulation. It suggests that even when there are no measurability concerns, for the nonnegative cost models, the mixed value and policy iteration approach may provide computationally efficient algorithms that are based on linear programming.

The mixed value and policy iteration method of this paper evolved from the enhanced policy iteration algorithmic framework proposed and analyzed in our earlier works for finite-state and control problems \cite{BerY-discount,YuB-ssp} and for abstract DP problems \cite{BerY-adp} under discounted and undiscounted total cost criteria (see also the book accounts of these works in \cite{Ber-DP12,Ber-DP13}).
In the finite-spaces or abstract DP context of these papers, measurability is not an issue. Instead, our main motivation has been the development of an asynchronous distributed, possibly model-free, policy iteration-like method that involves a flexible form of policy evaluation. A stopping problem formulation of policy evaluation was first proposed in \cite{BerY-discount} for that purpose. It has led to computationally efficient policy iteration-like algorithms with guaranteed convergence for the finite-spaces and abstract DP problems considered in \cite{BerY-discount,YuB-ssp,BerY-adp}, even in the context of an asynchronous distributed computing environment.
The method in this paper is based on the same idea and shares many important features with its counterparts in our earlier works, although its form has been modified and extended, in order to overcome the measurability issues in the present general-spaces stochastic control context. By providing a Borel-space counterpart of the method, one of our purposes is also to demonstrate that the mixed value and policy iteration approach is useful for addressing issues of not only computational but also theoretical nature. 

Although the focus of this paper is not so much on computation as on analytical issues, we note that our results suggest several computational benefits of the proposed method. First, some of its policy iteration-like variants may in practice converge significantly faster than value iteration (which our method contains as a special case). As will be discussed in Section~\ref{sec3.2}, these variants resemble modified policy iteration, a well-known method that combines value and policy iteration, 
which many computational studies and some analyses have found to be more efficient than both value iteration and policy iteration (see \cite[Chap.\ 6]{Puterman94}). Our own computational experiments with our method for finite-spaces problems \cite{BerY-discount} are also consistent with this assessment. 
Second, compared to standard policy iteration, our method offers better convergence properties, thanks to its inherent value-iteration character. This is evidenced by the convergence results of this paper, especially in the cases of nonpositive and nonnegative models, for which it is well-known that standard policy iteration can have convergence difficulties even when it can be carried out. This advantage was also found earlier for certain total cost finite-spaces problems in our earlier work \cite{YuB-ssp}. Third, like its predecessors in \cite{BerY-discount,YuB-ssp,BerY-adp}, the proposed method is suitable for asynchronous distributed computation (see the discussion in Section~\ref{sec7}), although we do not discuss this possibility in the present paper. In a distributed computing environment, based on the analyses in \cite{BerY-discount,YuB-ssp,BerY-adp} and the convergence results of this paper, we expect some of our policy iteration-like algorithms to be more efficient than asynchronous value iteration and at the same time, to converge more reliably than asynchronous policy iteration. Finally, another distinctive feature of the proposed method is that it can solve problems that cannot be handled by standard policy iteration. These include, e.g., problems in which exact policy improvement is impossible or $\epsilon$-optimal stationary policies do not exist. Potentially one may extract near-optimal policies from the iterates produced by the method. For such problems, however, the overall efficiency of our method relative to value iteration still requires further study.

The paper is organized as follows. In Section~\ref{sec2}, we provide background. 
In Section~\ref{sec3}, we introduce the mixed value and policy iteration method, and derive various algorithmic versions. We give greater attention to policy iteration-like algorithms, and we discuss their relation with standard policy iteration.
In Section~\ref{sec4}, we prove convergence results for the proposed method, for discounted problems with bounded one-stage costs and for total cost problems with nonpositive one-stage costs.
In Section~\ref{sec5}, we consider total cost problems with nonnegative one-stage costs. We first prove the new convergence theorem for value iteration in Section~\ref{secP5.1}. We then derive convergence results for the proposed method in Section~\ref{secP5.2}. 
In Section~\ref{sec6}, we discuss the application of the mixed value and policy iteration approach in semicontinuous models,
and we also give a result on the structure of the optimal cost function and optimal policies for nonnegative upper semicontinuous models.
In Section~\ref{sec7}, we conclude the paper with remarks on extensions and future research directions.
Appendices \ref{appsec-optstop}-\ref{appsec-Pexample} collect some related formulations, proofs and illustrative examples; Appendix \ref{appsec-altprfP} gives an alternative proof of the new convergence theorem for value iteration.

 \section{Background} \label{sec2}
 \markboth{\rm \S \thesection. Background}{\rm \S \thesection. Background}
 
In this section we describe the stochastic control framework with universally measurable policies. We give a brief summary of basic optimality results for infinite horizon, discounted and undiscounted total cost problems. We then explain the measurability issues in standard policy iteration. 

\subsection{Preliminaries} \label{sec2.0}
In this subsection we introduce some concepts and terminology, including
universal $\sigma$-algebras, analytic sets and lower semi-analytic functions.
We also highlight some properties that are important and provide the basis for the stochastic control framework. 
Several definitions and notation that we will use throughout the rest of the paper are given at the end of this subsection.

We focus on separable metrizable spaces.  
For such a space $X$, we denote by $\B(X)$ the Borel $\sigma$-algebra and by $\P(X)$ the set of Borel probability measures on $\B(X)$.
The space $\P(X)$ equipped with the weak topology is separable and metrizable.  
Let $\mathcal{F}$ be a $\sigma$-algebra on $X$ and let $Y$ be a separable metrizable space.
A function (or mapping) $f: X \to Y$ is $\mathcal{F}$-measurable if $f^{-1}(B) \in \mathcal{F}$ for any $B \in \B(Y)$; it is Borel measurable if $\mathcal{F} = \B(X)$.
We define likewise $\mathcal{F}$-measurable functions from a subset $X' \subset X$ to $Y$, where the $\sigma$-algebra on $X'$ is taken to be the trace $\sigma$-algebra $\mathcal{F} \cap X' =\{ D \cap X' \mid D \in \mathcal{F} \}$.

In this paper we will need to consider $\sigma$-algebras finer than the Borel $\sigma$-algebra, one of which is the universal $\sigma$-algebra defined as follows.
Any Borel probability measure $p \in \P(X)$ can be extended to a probability measure on the $\sigma$-algebra $\B_p(X)$ generated by $\B(X)$ and all the subsets of $X$ that have $p$-outer measure zero, such that the extension agrees with the $p$-outer measure on $\B_p(X)$. This extension of $p$ is called the \emph{completion of $p$} \cite[Sec.\ 3.3]{Dud02} and will also be denoted by $p$.
The intersection $\U(X) = \cap_{p \in \P(X)} \B_p(X)$ is called the \emph{universal $\sigma$-algebra} \cite[Def.\ 7.18]{bs}.
Sets in $\U(X)$ and measurable functions on $(X,\U(X))$ are said to be \emph{universally measurable}, since they are measurable with respect to the completion of any $p \in \P(X)$.

We consider subsets of a Polish space -- a topological space that can be metrized by a metric under which it is separable and complete.
In this paper, a \emph{Borel space} refers to a Borel subset of a Polish space,%footnote starts
\footnote{The definition of Borel spaces given in~\cite{bs} is more general and yet essentially equivalent to the one we give here. 
Borel spaces, as defined here, are 
now commonly called \emph{standard Borel spaces} (see e.g.~\cite{Sriv-borel}).}
%footnote ends 
endowed with the relative topology and Borel $\sigma$-algebra.
Analytic sets in a Polish space can be defined in several equivalent ways (see e.g., \cite[Prop.\ 7.41]{bs}, \cite[Sec.\ 13.2]{Dud02}), 
and roughly speaking, they are the images of Borel subsets of some Polish space under continuous or Borel measurable functions.
More specifically, in a Polish space $Y$,
the empty set is analytic by definition, and
a nonempty set $D$ is \emph{analytic} if $D=f(B)$ for some Borel set $B$ in a Polish space and Borel measurable function $f:B \to Y$ \cite[Thm.\ 13.2.1(c')]{Dud02}.
In a Polish space every Borel set is analytic 
and every analytic set is universally measurable (\cite[Cor.\ 7.42.1]{bs}, \cite[Thm.\ 13.2.6]{Dud02}). 

If $X$ is a Borel space or an analytic set, 
the \emph{analytic $\sigma$-algebra} $\A(X)$ is the $\sigma$-algebra generated by the analytic subsets of $X$, and it lies between the Borel $\sigma$-algebra and the universal $\sigma$-algebra:
$\B(X) \subset \A(X) \subset \U(X)$ (the inclusions are strict if $X$ is an uncountable Borel space) \cite[p.\ 171]{bs}.%footnote starts
\footnote{By an analytic subset of $X$ we mean a subset of $X$ which, as a subset of the Polish space $Y$ that contains $X$, is analytic. The $\sigma$-algebras $\B(X), \A(X), \U(X)$ here are the trace $\sigma$-algebras $X \cap \B(Y), X \cap \A(Y), X \cap \U(Y)$, respectively.}
%footnote ends
Thus on $X$, a Borel measurable function is analytically measurable (i.e., $\A(X)$-measurable) and an analytically measurable function is universally measurable.

The class of analytic sets in a Polish space is closed under countable unions and countable intersections, and moreover, Borel preimages of analytic sets are also analytic (\cite[Cor.\ 7.35.2, Prop.\ 7.40]{bs}, \cite[Chap.\ 4]{Sriv-borel}). 
This gives rise to many nice properties of lower semi-analytic functions, that is, functions whose lower level sets are analytic. 
More specifically, a function $f: D \to [-\infty, \infty]$ is said to be \emph{lower semi-analytic} if $D$ is an analytic set and 
for every $c \in \Re$, the level set $\{ x \in D \!\mid f(x) < c\}$ of $f$ is analytic \cite[Def.\ 7.21]{bs}. (Equivalently, the epigraph of $f$, $\{(x, c) \!\mid x \in D, f(x) \leq c, c \in \Re\}$, is analytic; cf.\ \cite[p.\ 186]{bs}.)  
Every lower semi-analytic function is universally measurable, 
since analytic sets are universally measurable. 
Moreover, based on the properties of analytic sets, the following operations on lower semi-analytic functions result in a lower semi-analytic function (see \cite[Lemma 7.30]{bs}):\moveup\moveup
\begin{itemize}
\item[(i)]
If $f, g: D \to [-\infty, \infty]$ are lower semi-analytic functions, then $f+g$ is lower semi-analytic (here we define $\infty - \infty = - \infty + \infty = \infty$). In addition, if $f, g \geq 0$ or if $g$ is Borel measurable and $g \geq 0$, then $f g$ is lower semi-analytic (here we define $0\cdot \infty = 0\cdot (-\infty)=\infty \cdot 0 = (-\infty) \cdot 0 =0$).  
Note a particular implication of this: for any $B \in \B(D)$, $f \cdot \b1_B$ is lower semi-analytic, 
where $\b1_B$ denotes the indicator function for $B$:
$\b1_B(x) = 1$ if $x \in B$; $\b1_B(x) = 0$ if $x \not\in B$.\moveup
\item[(ii)] If $g:X \to Y$ is Borel measurable, where $X, Y$ are Borel spaces, and $f: g(X) \to [-\infty, \infty]$ is lower semi-analytic, then the composition $f \circ g$ is lower semi-analytic. \moveup 
\item[(iii)] For a sequence of lower semi-analytic functions $f_n: D \to [-\infty, \infty]$, $n \geq 1$, the functions $\inf_n f_n$, $\sup_n f_n$, $\liminf_n f_n$ and $\limsup_n f_n$ are all lower semi-analytic. (These are pointwise definitions.) \moveup\moveup
\end{itemize}

Several properties of analytic sets and lower semi-analytic functions play instrumental roles in the stochastic control framework we will introduce. 
They concern analytic sets in product spaces or functions involving two variables.
The first property is closely related to value iteration and the structure of the optimal cost function in the stochastic control context.
If $D$ is an analytic set in $X \times Y$, where $X, Y$ are Polish, 
the projection of $D$ on $X$, $\text{proj}_X(D) = \{ x \!\mid (x,y) \in D \ \text{for some} \ y \}$, is analytic \cite[Prop.\ 7.39]{bs}.
When applied to level sets of functions, an implication of this is that if $D \subset X \times Y$ is analytic and $f: D \to [-\infty, \infty]$ is lower semi-analytic, then after partial minimization of $f$ over the vertical sections $D_x$ of $D$ for each $x$, the resulting function $f^*: \text{proj}_X(D) \to [-\infty, \infty]$ given by
\begin{equation} \label{eq-minf}
   f^*(x) = \inf_{y \in D_x} f(x, y), \quad \text{where} \ D_x = \{ y \mid (x, y) \in D \}, 
\end{equation}   
is also lower semi-analytic \cite[Prop.\ 7.47]{bs}.
  
The Jankov-von Neumann selection theorem asserts that if $A$ is an analytic set in $X \times Y$, where $X, Y$ are Polish, then there exists an analytically measurable function $\phi: \text{proj}_X(D) \to Y$ such that the graph of $\phi$ lies in $D$, i.e., $(x, \phi(x)) \in D$ for all $x \in \text{proj}_X(D)$~\cite[Prop.\ 7.49]{bs}. For minimization problems of the form~(\ref{eq-minf}), the theorem is applied to the level sets or epigraphs of lower semi-analytic functions, and together with other properties,
it yields the existence of an analytically measurable $\epsilon$-minimizer and the existence of a universally measurable $\epsilon$-minimizer that attains the minimum $f^*(x)$ at every $x$ where $f^*(x)$ is attained by some $y \in D_x$. For details, see the selection theorems given in~\cite[Prop.\ 7.50(a)-(b)]{bs}. 
In the stochastic control context, this is closely related to the existence of optimal or nearly optimal policies and their structures.

Another important property of lower semi-analytic functions involves integration and stochastic kernels.
Let $X$ and $Y$ be Borel spaces.
In this paper, a \emph{Borel}, \emph{analytically} or \emph{universally measurable stochastic kernel} on $Y$ given $X$ is a mapping $\kappa (\cdot \!\mid \cdot) : \B(Y) \times X \to [0,1]$ such that:\moveup\moveup
\begin{itemize}
\item[(i)] For each $B \in \B(Y)$, the function $\kappa(B \!\mid \cdot): X \to [0,1]$ is $\mathcal{F}$-measurable with $\mathcal{F} = \B(X), \A(X)$ or $\U(X)$, respectively.\moveup
\item[(ii)] For each $x \in X$, 
$\kappa(\cdot \!\mid x)$ is a probability measure on $(Y, \B(Y))$.\moveup\moveup
\end{itemize}
Equivalently, the function $x \mapsto \kappa(\cdot \!\mid x)$ is $\mathcal{F}$-measurable from $X$ to $\P(Y)$~\cite[Prop.\ 7.26, Lemma 7.28, Prop.\ 11.6]{bs}.
If $f: X \times Y \to [0, \infty]$ is lower semi-analytic and $\kappa(dy \!\mid x)$ is a Borel measurable stochastic kernel on $Y$ given $X$, then the integral 
$$ \int_Y f(x, y) \, \kappa(dy \!\mid x) $$
as a function of $x$ is lower semi-analytic on $X$~\cite[Prop.\ 7.48]{bs}, where for each $x$, the integration is defined to be with respect to the completion of the Borel probability measure $\kappa(dy \!\mid x)$.
If instead $\kappa(dy \!\mid x)$ is analytically or universally measurable, then the above integral as a function of $x$ is universally measurable~\cite[Prop.\ 7.46, Sec.\ 11.2]{bs} but \emph{not} necessarily lower semi-analytic.
These facts are closely related to the structure of the cost functions and the selection of measurable policies in the stochastic control context.
  
For more properties of analytic sets and lower semi-analytic functions, see the paper \cite{BFO74} and the monograph \cite[Chap.\ 7]{bs}. (For general properties of analytic sets, see also the books \cite{Par67,Sriv-borel}.) 
Later in this paper, we will also encounter special classes of semicontinuous functions and continuous stochastic kernels when we consider semicontinuous control models, and we will give the definitions of these objects then in Section~\ref{sec6}.

Let us introduce now some notation that we will use throughout the paper. 
For a Borel space or analytic set $X$, we denote by $\M(X)$ (resp.\ $A(X)$) the set of functions $f: X \to [-\infty, \infty]$ that are universally measurable (resp.\ lower semi-analytic), and we denote by $A_b(X)$, $A_+(X)$ and $A_-(X)$ the subsets of functions in $A(X)$ that are bounded, nonnegative and nonpositive, respectively.
The set $\M_b(X)$ of bounded universally measurable functions on $X$ is a Banach space with the norm $\textstyle \| f \|_\infty = \sup_{x \in X} | f(x) |$ for $f \in \M_b(X)$.
It is worth noting that the set $A_b(X)$ is a closed subset of $\M_b(X)$ and hence for the metric $d_{sup}(f, f') = \| f - f' \|_\infty$, the space $(A_b(X), d_{sup})$ is a complete metric space. 
We will consider convergence of functions in such spaces for discounted control problems.
For undiscounted control problems, we will consider primarily pointwise convergence of functions on $X$. 
We write $f_n \to f$ for a sequence of functions $f_n$ converging pointwise to a function $f$, and if the convergence is monotonically from above or from below, we write $f_n \downarrow f$ or $f_n \uparrow f$, respectively. 
For $x \in X$, we denote by $\delta_x$ the Dirac measure that assigns probability $1$ to the point $x$. 
We use the symbol $\Re$ for the set of reals, and the symbol $\0$ for the constant function zero on any space given in the context of discussion.

 \subsection{Stochastic Control Model} \label{sec2.1}
 
 Our stochastic control model involves a state space $S$ and a control space $C$, which are assumed to be Borel spaces. 
We will write $x$ for a state in $S$ and $u$ for a control in $C$.
At each state $x \in S$, one can apply a control from a nonempty subset $U(x) \subset C$.
The set-valued function $U$ given by $x \mapsto U(x)$ specifies the control constraint for all states.
We assume that the graph of $U$,
$$ \Gamma = \{ (x, u) \mid x \in S, u \in U(x) \}, 
$$
is an analytic subset of $S \times C$.
Applying a control $u$ at a state $x$ incurs a possibly infinite one-stage cost and moves the system to another state $x'$. 
The one-stage cost is given by $g(x,u)$, where $g: \Gamma \to [ -\infty, \infty]$ is assumed to be a lower semi-analytic function.
The transition to state $x'$ is according to a Borel measurable stochastic kernel $q(dx' \mid x, u)$ on $S$ given $S \times C$.

We consider measurable policies defined as follows.
A \emph{universally measurable policy} is a sequence $\pi=(\mu_0, \mu_1, \ldots)$, where for each $k$,
$\mu_k\big(du_k \!\mid x_0, u_0, \ldots, u_{k-1}, x_k \big)$ is a universally measurable stochastic kernel on $C$ given $(S \times C)^{k} \times S$ such that 
\begin{equation}  \label{eq-control-constraint}
   \mu_k\big(U(x_k) \!\mid x_0, u_0, \ldots, u_{k-1}, x_k \big) = 1, \qquad \forall \, (x_0, u_0, \ldots, u_{k-1}, x_k) \in (S \times C)^k \times S.
\end{equation}   
Here the constraint (\ref{eq-control-constraint}) says that the set of non-admissible controls, $C \setminus U(x_k)$, has probability zero, and
this is meaningful since $\Gamma$ is analytic: 
each vertical section $U(x)$ of $\Gamma$ is universally measurable \cite[Lemma 7.29]{bs} and hence measurable with respect to the completion of the Borel probability measure $\mu_k(\cdot \mid x_0, u_0, \ldots, u_{k-1}, x_k )$.
In what follows, when no confusion arises, we will simply refer to universally measurable policies as policies. 

A policy $\pi$ is \emph{nonrandomized} if $\mu_k\big(\cdot \!\mid x_0, u_0, \ldots, u_{k-1}, x_k \big)$ is a Dirac measure that assigns probability one to some point in $U(x_k)$ for every $k$ and every $(x_0, u_0, \ldots, u_{k-1}, x_k)$.
A policy $\pi$ is \emph{semi-Markov} if for every $k$, the function $(x_0, u_0, \ldots, u_{k-1}, x_k) \mapsto  \mu_k(\cdot \!\mid x_0,  u_0, \ldots, u_{k-1}, x_k)$ depends only on $(x_0, x_k)$;
\emph{Markov} if for every $k$, the latter function depends only on $x_k$; \emph{stationary} if $\pi$ is Markov and $\mu_k = \mu$ for all $k$. For the stationary case, we simply write $\mu$ for $\pi = (\mu, \mu, \ldots)$. 
A nonrandomized stationary policy $\mu$ can be viewed as a mapping that maps $x \in S$ to a point in $U(x) \subset C$.  
We denote this mapping also by $\mu$, and we will use both notations $\mu(x)$, $\mu(du\!\mid x)$ in the paper, depending on the context. 

A policy $\pi$ is \emph{Borel measurable} or \emph{analytically measurable} if each component $\mu_k$ of $\pi$ is a Borel measurable or analytically measurable stochastic kernel; such a policy is by definition also universally measurable.
Because $\Gamma$ is analytic, by the Jankov-von Neumann selection theorem \cite[Prop.~7.49]{bs}, there exists at least one universally measurable (in fact, analytically measurable) nonrandomized stationary policy. A Borel measurable policy, however, may not exist \cite{Blk-borel}.

We denote by $\Pi$ the set of universally measurable policies. 
Given a policy $\pi \in \Pi$, the collection of stochastic kernels 
\begin{align*}
 & \mu_0(d u_0 \mid x_0), \  q(dx_1 \mid x_0, u_0), \  \mu_1(du_1, \mid x_0, u_0), \ q(dx_2 \mid x_1, u_1), \ \ldots, \qquad \\
 & \ldots, \ \mu_k\big(d u_k \mid x_0, u_0, \ldots, u_{k-1}, x_k \big), \ q(d x_{k+1} \mid x_k, u_k), \ \ldots,
\end{align*}  
uniquely determines, for each initial distribution $p_0$ of $x_0$, a probability measure $r(\pi, p_0)$ on the universal $\sigma$-algebra on $(S \times C)^\infty$ with the following property \cite[Prop.\ 7.45]{bs}:%footnote starts
\footnote{It is worth noting that the universal $\sigma$-algebra on $(S \times C)^\infty$ is not a product $\sigma$-algebra, so the existence of a unique probability measure $r(\pi, p_0)$ here does not follow immediately from the Ionescu Tulcea theorem.}
%footnote ends
with respect to  $r(\pi, p_0)$, 
the expectation $\E f$ for any nonnegative, universally measurable function $f: (S \times C)^{k+1} \to [0, \infty]$ equals the iterated integral 
$$ \int_S \int_C \cdots \int_{S} \int_{C}  f(x_0, u_0, \ldots, x_{k}, u_{k})  \mu_k(du_{k} \mid x_0, u_0, \ldots, x_k) \, q(dx_k \mid x_{k-1}, u_{k-1}) \, \cdots \,\mu_0(u_0 \mid x_0) \, p_0(dx_0).$$
Here and in what follows, an integral $\int f dp$ that involves a universally measurable function $f$ and a Borel probability measure $p$, is defined to be the integral of $f$ with respect to the completion of $p$. 

In general, for a measurable, extended real-valued function $f$, $\E f$ is defined as usual to be $\E f^+ - \E f^-$, where $f^+ = \max \{ 0, f \}$, $f^- = - \min \{ 0, f\}$. 
The convention $\infty - \infty =  - \infty + \infty =  \infty$ will be adopted, although in the control problems we consider, we will not encounter such summations.

\vspace*{-0.01in}
\subsubsection*{Infinite Horizon Total Cost Problems} 

In this paper we formulate various stochastic control problems as cost minimization problems.
We consider primarily three total expected cost criteria: 
discounted total cost problems with bounded one-stage costs (D), and undiscounted total cost problems with nonpositive one-stage costs (N) and with nonnegative one-stage costs (P). Specifically,
let $\alpha \in [0,1]$ be the discount factor.%\moveup\moveup
\begin{itemize}
\item[(D)] 
$\alpha < 1$ and $-b \leq g(x,u) \leq b$ for all $(x, u) \in \Gamma$, where $b \in \Re$.\moveup
\item[(N)] 
$\alpha = 1$ and $g \leq 0$.\moveup
\item[(P)] 
$\alpha = 1$ and $g \geq 0$.\moveup
\end{itemize}
We mention that for reward maximization (instead of cost minimization), the reverse terminologies are used in the literature \cite{Blk-positive,Str-negative,MS92,Puterman94}: case (N) here corresponds to the positive model and case (P) to the negative model considered there. 

In each of the (D)(N)(P) cases, we define the cost of $\pi \in \Pi$ for an initial state $x_0=x \in S$ to be
$$ J_\pi(x) = \E^\pi \left\{ \sum_{k=0}^\infty \alpha^k g(x_k,u_k) \right\},$$
the expectation of the universally measurable function $\sum_{k=0}^\infty \alpha^k g(x_k,u_k)$ 
with respect to the probability measure $r(\pi, \delta_{x})$, which is induced by $\pi$ and the initial distribution $\delta_x$ 
as described earlier. (Although $g$ is only defined on $\Gamma$, $\pi$ is a policy and satisfies the control constraint, so $(x_k,u_k) \in \Gamma$ for all $k$ with probability one and the expectation is thus well-defined.) 
By the bounded convergence theorem (for case (D)) and the monotone convergence theorem (for cases (N)(P)), we can also write $J_\pi(x)$ as
$$ J_\pi(x) = \sum_{k=0}^\infty \alpha^k \, \E^\pi \big\{ g(x_k,u_k) \big\},$$
where the expectation is with respect to the marginal of $r(\pi, \delta_{x})$ on the space $S \times C$ of $(x_k, u_k)$.
For all $\pi \in \Pi$, the cost functions $J_\pi$ are universally measurable \cite[p.\ 215]{bs}.

The optimal cost function is defined by the minimal cost of universally measurable policies for each state:
$$ J^*(x) = \inf_{\pi \in \Pi} J_\pi(x), \qquad \forall \, x \in S.$$
If $J_\pi(x) = J^*(x)$, $\pi$ is optimal for state $x$.
For $\epsilon > 0$, $\pi$ is said to be $\epsilon$-\emph{optimal} if for all $x \in S$,
$$ J_\pi(x) \leq \begin{cases}
     J^*(x) + \epsilon & \text{if} \ J^*(x) > - \infty;\\
     - 1 / \epsilon & \text{if} \ J^*(x) = - \infty.
     \end{cases}
$$

We mention that certain nonnegative or nonpositive discounted problems where the discount factor at each stage depends on the state transition (see e.g., \cite{Schal75,Whit79}), can be converted to the (N)(P) models given above.%footnote starts
\footnote{
Suppose the total cost of $\pi$ for each initial state $x_0=x$ is defined as
$$ J_\pi(x) =   \E^\pi \left\{ \hat g(x_0, u_0, x_1) + \sum_{k=1}^\infty  \Big(  \prod_{i=1}^{k} \beta(x_{i-1}, u_{i-1}, x_i) \Big) \cdot \hat g(x_k,u_k, x_{k+1}) \right\},$$
where $\beta: S \times C \times S \to [0,1]$ is the transition-dependent discount factor, and 
$\hat g: \Gamma \times S \to [-\infty, 0]$ or $[0,\infty]$ is the transition cost function. 
Assume that $\beta$ is Borel measurable and $\hat g$ is lower semi-analytic.
This problem can be converted to an equivalent problem of type (N) or (P) by introducing an absorbing cost-free state $\infty$, defining the new state transition kernel $\tilde q(\cdot \mid \cdot)$ by
\begin{align*}
 \tilde q( B \! \mid x, u)  = \int_B \beta(x, u, x') \, q(dx' \!\mid x, u), \ \ \ B \in \B(S), \qquad 
 \tilde q \big(\{\infty\} \! \mid x, u \big)  = 1 - \tilde q(S \! \mid x, u),  
\end{align*}
for $(x, u) \in \Gamma$, and letting the one-stage cost be $g(x, u) = \int_S \hat g(x, u, x') \, q(dx' \!\mid x, u)$ on $\Gamma.$
}
%footnote ends
Results for (N)(P) including those given in this paper are thus applicable to these problems with transition-dependent discounting as well.

\subsection{Optimality Properties} \label{sec2.2}

In each of the (D)(N)(P) cases, the optimal cost function $J^*$ is lower semi-analytic, and it satisfies the optimality equation
$J^* = T (J^*),$
where $T$ maps $A(S)$ into $A(S)$ and is given by
\begin{equation} \label{eq-T}
  T (J)(x) = \inf_{u \in U(x)} \left\{ g(x,u) +  \alpha \int_S  J(x') \, q(dx' \mid x, u) \right\}, \qquad  x \in S.
\end{equation}  
We will refer to $T$ as the optimal cost operator. We note that $T(J) \in A(S)$ for any function $J\in A(S)$, as just mentioned. This is a direct consequence of the properties of lower semi-analytic functions discussed in Section~\ref{sec2.0} and the stochastic control model given in Section~\ref{sec2.1}. More specifically, since the state transition kernel $q(dx' \!\mid x, u)$ is Borel measurable and the one-stage cost $g$ is lower semi-analytic, $g(x, u) + \alpha \int_S  J(x') \, q(dx' \mid x, u)$ as a function of $(x,u)$ is lower semi-analytic on $\Gamma$ by \cite[Prop.\ 7.48 and Lemma 7.30(4)]{bs}. Then, since partial minimization preserves lower semi-analyticity \cite[Prop.\ 7.47]{bs}, it follows that $T(J)$ is lower semi-analytic.

In each of the (D)(N)(P) cases, the cost function $J_\mu$ for a stationary policy $\mu$ is universally measurable. 
It satisfies a linear equation,
$J_\mu = T_\mu (J_\mu),$ 
where $T_\mu$ is a mapping from $\M(S)$ to $\M(S)$, given by
\begin{equation} \label{eq-Tmu}
  T_\mu (J)(x) = \int_C \left( g(x,u) +  \alpha \int_S  J(x') \, q(dx' \mid x, u) \right) \mu(du \mid x), \qquad  x \in S.
\end{equation}  

In terms of the convergence properties of the value iteration sequence $T^k(J)$ and the structures of the optimal policies, the (D)(N)(P) cases differ.
\moveup\moveup
\begin{itemize}
\item[(a)] For (D)(N), value iteration converges to $J^*$. In particular, in case (D),
$T^k (J) \to J^*$
for any bounded lower semi-analytic function $J$,  and in case (N), $T^k (\0) \downarrow J^*$.\moveup
\item[(b)] For (P), value iteration need not converge to $J^*$: 
$ T^k (\0) \uparrow J_\infty \leq J^*,$
where the pointwise limit $J_\infty$ of $\{T^k (\0)\}$ satisfies $J_\infty \leq T (J_\infty)$.\moveup\moveup
\end{itemize}
In all three cases, $\epsilon$-optimal nonrandomized policies exist for each $\epsilon > 0$; however, they can be taken to be stationary for (D), \emph{semi-Markov} for (N), and \emph{Markov} for (P). An $\epsilon$-optimal randomized Markov policy need not exist for (N) (a counterexample was given by van der Wal~\cite[Example 2.26]{vdW81}; see also \cite[p.\ 326]{Puterman94}).
If for each state $x$, an optimal policy exists, then:\moveup\moveup
\begin{itemize}
\item[(a)] For (D)(P), an optimal nonrandomized stationary policy exists.\moveup
\item[(b)] For (N), an optimal randomized semi-Markov policy exists.\moveup\moveup
\end{itemize}

We refer to \cite[Chap.\ 9]{bs} the optimality properties mentioned above, as well as finer characterizations of the optimal cost function and optimal policies. 

 \subsection{Measurability Issues in Standard Policy Iteration} \label{sec2.3}
 
In the policy iteration scheme, 
we repeat the following two steps:\moveup\moveup
\begin{itemize}
\item[(i)] Evaluate the cost function $J_\mu$ of a given stationary policy $\mu$.\moveup
\item[(ii)] Find a stationary policy $\mu'$ with
 $ T_{\mu'} (J_\mu) = T (J_\mu)$
and go to step (i) with $\mu = \mu'$.\moveup\moveup
 \end{itemize}
A variant of this scheme is the modified policy iteration \cite{Puterman94}:\moveup\moveup
\begin{itemize}
\item[(i')]  For a given stationary policy $\mu$ and a given function $J$, compute as an approximation of $J_\mu$,  
$J' = T_{\mu}^m(J)$ for some positive integer $m$.\moveup
\item[(ii')] Find a stationary policy $\mu'$ with
 $ T_{\mu'} (J') = T (J')$
and go to step (i') with $\mu = \mu'$ and $J = J'$.%\moveup\moveup
\end{itemize}
Both schemes break down, however, for the stochastic control model with universally measurable policies, due to measurability issues (cf.\ \cite[p.~940]{BFO74}, \cite[p.~232]{bs}). We explain the reasons below.

As defined in (\ref{eq-T}), $T$ maps a universally measurable function $J \in \M(S)$ to the function $T(J)$, possibly outside $\M(S)$.
For a stationary policy $\mu$, $J_\mu$ is universally measurable, so $T(J_\mu)$ is defined. 
But since $J_\mu$ need not be lower semi-analytic, even if $T(J_\mu)$ is universally measurable, a stationary, universally measurable policy $\mu'$ such that
\begin{equation} \label{eq-polite}
  T_{\mu'} (J_\mu) = T (J_\mu) \quad \text{or} \quad T_{\mu'} (J_\mu) \leq T (J_\mu) +\epsilon,  \quad \text{for some given} \ \epsilon > 0,
\end{equation}  
may not exist (\cite[p.~940]{BFO74}, \cite[p.~232]{bs}). When this happens, step (ii) of policy iteration cannot be carried out. The same issue also causes modified policy iteration to break down.

Blackwell et al.\ \cite[Example (48)]{BFO74} gave an example of an analytically measurable function $J$ on $[0,1]$ for which $T(J)$ is not Lebesgue measurable. If $J_\mu$ equals such $J$, then there is certainly no stationary (universally measurable) policy $\mu'$ that can satisfy $T_{\mu'} (J) = T (J)$, because $T_{\mu'} (J)$ is universally measurable, whereas $T(J)$ is not.  Moreover, since $T(J)$ is not universally measurable, for some $p \in \P(S)$, $T(J)$ is not integrable with respect to the completion of $p$. Hence, $T^2(J)$ as well as $(T_\mu \circ T)(J)$ for a stationary policy $\mu$ can be undefined for some states $x$ (cf.\ \cite[Example (48)]{BFO74}). This means that variants of policy iteration of the form 
$J_{k+1} = T_k (J_{k}),$ where some of the $T_k$'s equal $T$ and others equal $T_\mu$ for some stationary policy $\mu$, can also run into trouble. 

The measurability issues in policy iteration can be assumed away by imposing stronger model assumptions. Most of these assumptions involve continuity and semicontinuity conditions to ensure that policy iteration can operate within the class of Borel measurable policies (see e.g., \cite{HL96}). (For an interesting case where policy iteration can operate with universally measurable policies without any measurability difficulty, see the subsequent Example~\ref{ex-choiceB2}.)
Nevertheless, without strong model assumptions, even upper and lower semicontinuous models (see Section~\ref{sec6}) have the measurability issues discussed above. 
Such issues arise because the cost function $J_\mu$ of a Borel measurable policy $\mu$ is Borel measurable but not necessarily (upper or lower) semicontinuous in those models, so for policy improvement, special selection theorems for semicontinuous functions (see Section~\ref{sec6}) cannot be applied to ensure the existence of a Borel measurable $\mu'$ that satisfies (\ref{eq-polite}). Then the policy improvement step may generate an analytically or universally measurable policy, subjecting the subsequent iterations to the measurability difficulties described above.

 \section{A Mixed Value and Policy Iteration Method} \label{sec3}
 \markboth{\rm \S \thesection. Mixed Value and Policy Iteration}{\rm \S \thesection. Mixed Value and Policy Iteration}

In this section we introduce formally our mixed value and policy iteration approach, and discuss some specific algorithms.
For (D)(N)(P), recall that the relation $J^* = T(J^*)$ holds: 
 $$ J^*(x)  = \inf_{x \in U(x)} \left\{ g(x, u) + \alpha \int_S J^*(x') \, q(dx' \!\mid x, u) \right\}, \qquad \forall \, x \in S.$$
 We define $Q^* \in A(\Gamma)$ by
  \begin{equation} \label{eq-Q}
    Q^*(x, u) = g(x, u) + \alpha \int_S J^*(x') \, q(dx' \!\mid x, u), \qquad \quad (x, u) \in \Gamma.
 \end{equation}   
For each $(x, u) \in \Gamma$, we may view $Q^*(x,u)$ as the result of cost minimization over controllers that start at state $x$, apply control $u$, and then choose some policy. 
This interpretation of $Q^*(x,u)$ is better revealed in the following equation, which is equivalent to (\ref{eq-Q}) \cite[Cor.\ 9.5.2]{bs}:
\begin{equation} \label{eq-Qalt}
 Q^*(x, u) = g(x, u) + \alpha  \inf_{\pi \in \Pi} \int_S J_{\pi}(x') \, q(dx' \!\mid x, u), \qquad \quad (x, u) \in \Gamma.
\end{equation} 
(In the literature on learning and simulation-based DP, $Q^*(x,u)$ is known as the optimal Q-factor associated with $(x,u)$; see e.g., \cite{BET,SUB}.)
To simplify notation, for any function $Q$ on $\Gamma$, let
 $$ M(Q)(x) = \inf_{u \in U(x)} Q(x, u), \qquad x \in S.$$
 The mapping $M$ maps $A(\Gamma)$ into $A(S)$ \cite[Prop.\ 7.47]{bs}.
With this notation, we can write the optimality equation in two equivalent ways:
 \begin{equation}  \label{eq-JQ^*}
   J^* = T(J^*) \qquad \Longleftrightarrow \qquad  J^* = M(Q^*).
 \end{equation}

Our mixed value and policy iteration method operates on the product space $A(S) \times A(\Gamma)$ and aims to compute $(J^*, Q^*)$. The method combines characteristics of both value and policy iteration, and the combination has two crucial features. First, it uses portions of a universally measurable policy that are Borel, to preserve the lower semi-analytic properties of the functions involved, thereby overcoming the measurability issues in standard policy iteration. Second, thanks to its value iteration character, it does not rely strongly on the behavior of policies for convergence. In particular, the policies involved are not required to be successively improving -- a requirement that in general cannot be met in our context or in the case where the policies involved are restricted to be Borel measurable \cite{Blk-discount}. Our method gives rise to various policy iteration-like algorithms, whose convergence we will analyze in Sections~\ref{sec4} and~\ref{sec5}.

In what follows we first introduce a family of mappings underlying the method and discuss their basic properties and their relation to optimal stopping problems (Section~\ref{sec3.1}). 
We then give various forms of algorithms (Section~\ref{sec3.2}), with emphasis on policy iteration-like methods, their properties, and their differences from standard policy iteration.
 
 \subsection{Mappings Induced by Stationary Policies} \label{sec3.1}
 
We define a family of parametrized mappings $F_\k$, which will be used later in a step of the mixed value and policy iteration algorithms that is analogous to policy evaluation. The parameters $\k$ include a policy component and a set component paired with that policy, introduced for overcoming measurability issues. Specifically,
let $\Theta$ denote the set of all pairs $(\mu, B)$, where $\mu$ is a stationary policy and $B$ a Borel subset of $S$, such that the function $x \mapsto \mu(du \!\mid x)$ restricted to $B$ is Borel measurable (equivalently, $\mu(D\!\mid \cdot)$ is Borel measurable on $B$ for every $D \in \B(C)$). 
For each $\k =(\mu, B) \in \K$, we define $F_\k: \M(\Gamma) \times \M(S) \to \M(\Gamma)$ by
\begin{align}
  F_\k ( Q \, ; J)(x, u) & =  g(x, u) + \alpha \int_{S\setminus B}  J(x') \, q(dx' \!\mid x, u)  \notag \\
        &  \ \ \ +  \alpha \int_B  \int_C  \min \big\{ J(x') \, , \, Q(x', u') \big\} \, \mu(du' \!\mid x') \, q(dx' \!\mid x, u), \qquad (x, u) \in \Gamma,  \label{def-F}
 \end{align}
for all $Q \in \M(\Gamma)$ and $J \in \M(S)$.
Here the convention $\infty - \infty = - \infty + \infty = \infty$ is used. We also note that although $Q$ is defined only on $\Gamma$, the inner integral of the third term in (\ref{def-F}) is well-defined because $\mu$ satisfies the control constraint. (We could, for example, view this integral as an integral for the extension of $Q$ to $S \times C$ with $Q(x',u') = \infty$ outside $\Gamma$.)
 
For any stationary policy $\mu$, the trivial choice $B = \emptyset$ gives $\k=(\mu, \emptyset) \in \K$, 
but the corresponding mapping $F_\k$ does not depend on the policy $\mu$ at all.
To introduce greater dependence of $F_\k$ on $\mu$, we desire ``large'' sets $B$. By the nature of universally measurable policies, one can indeed find ``large'' $B$ with $(\mu, B) \in \K$ (see Prop.~\ref{prp-F1}(b) below and Examples~\ref{ex-choiceB2}, \ref{ex-choiceB1} in Section~\ref{sec3.2}).
If the policy $\mu$ is Borel measurable, then $(\mu, S) \in \K$, so one may let $B$ be the entire space.

An important property of $F_\k$ is that it preserves the lower semi-analyticity of functions. This will allow us to overcome the measurability difficulties that hamper standard policy iteration.

\begin{prop} \label{prp-F1} \hfill \moveup\moveup
\begin{itemize}
\item[(a)] For any $\k \in \K$ and $J \in A(S)$, $F_\k(\cdot \, ; J)$ maps $A(\Gamma)$ into $A(\Gamma)$.\moveup%\moveup
\item[(b)] For each stationary policy $\mu$, given any $p \in \P(S)$, 
there is a Borel set $B \subset S$ with $p(S\setminus B) = 0$ and $(\mu, B) \in \K$.\moveup
\end{itemize}
 \end{prop}

\begin{proof}
\noindent (a) Let $Q \in A(\Gamma)$. To show $F_\k(Q \, ; J) \in A(\Gamma)$, we show that each term in its definition~(\ref{def-F}) is a lower semi-analytic function on $\Gamma$.
The first term $g(x,u)$ is lower semi-analytic by definition. 
The second term equals $\alpha \int_S  \b1_{S \setminus B}(x') J(x') \, q(dx' \!\mid x, u)$. Here the function $\b1_{S \setminus B} \cdot J$ is lower semi-analytic given that $S \setminus B$ is a Borel set and $J$ is lower semi-analytic~\cite[Lemma 7.30(4)]{bs}, and the stochastic kernel $q(dx' \!\mid x, u)$ is by definition Borel measurable. So by \cite[Prop.\ 7.48]{bs} this integral as a function of $(x,u)$ is lower semi-analytic on $S \times C$ and hence lower semi-analytic on the analytic set $\Gamma$.
For the third term,
$    \alpha \int_B  \int_C  \min \big\{ J(x') \, , \, Q(x', u') \big\} \, \mu(du' \!\mid x') \, q(dx' \!\mid x, u),$ 
since $\mu$ satisfies the control constraint, we can write it equivalently as
\begin{equation}
   \alpha \int_S  \int_C  f(x',u')  \, \mu(du' \!\mid x') \, q(dx' \!\mid x, u), \label{eq-prfF1-2}
\end{equation}    
where $f: S \times C \to [-\infty, \infty]$ is given by
$f(x',u') = \b1_B(x') \cdot \min \{ J(x')  ,  Q^e(x', u') \}$ for $(x',u') \in S \times C$, with $Q^e$ being a lower semi-analytic extension of $Q$ to $S \times C$ defined as $Q^e(x,u) = \infty$ for $(x,u) \not\in \Gamma$. 
The function $f$ is lower semi-analytic, since $J$ and $Q^e$ are lower semi-analytic functions and $B$ is a Borel set \cite[Lemma 7.30(2),(4)]{bs}; thus $f$ is also lower semi-analytic on $B \times C$. The fact $(\mu, B) \in \K$ implies that $\mu(d u' \!\mid x')$ is a Borel measurable stochastic kernel on $C$ given $B$. Consequently, the inner integral $ \int_C  f(x',u')  \, \mu(du' \!\mid x')$ is lower semi-analytic on $B$ by~\cite[Prop.\ 7.48]{bs}. We also have $ \int_C  f(x',u')  \, \mu(du' \!\mid x') = 0$ for $x' \not\in B$. Therefore, $ \int_C  f(x',u')  \, \mu(du' \!\mid x')$  is lower semi-analytic on $S$. Then, since $q(dx' \!\mid x, u)$ is a Borel measurable stochastic kernel on $S$ given $S \times C$, the integral (\ref{eq-prfF1-2}) as a function of $(x,u)$ is lower semi-analytic on $S \times C$ by~\cite[Prop.\ 7.48]{bs} and hence lower semi-analytic on the analytic set $\Gamma$. 
This proves part (a).

\noindent (b) Since $\mu(du \!\mid x)$ is a universally measurable stochastic kernel on $C$ given $S$, 
by \cite[Lemma 7.28]{bs}, there is a Borel measurable stochastic kernel $\tilde \mu(du \!\mid x)$ 
with $\tilde \mu(du \!\mid x) = \mu(du \!\mid x)$ for $p$-almost every $x$. Part (b) then follows.
\end{proof}

In the discounted case (D), we work with $J \in A_b(S), Q \in A_b(\Gamma)$, the subsets of bounded lower semi-analytic functions. In the nonpositive case (N), we work with $J \in A_-(S), Q \in A_-(\Gamma)$, the subsets of nonpositive lower semi-analytic functions, whereas in the nonnegative case (P), we work with $J \in A_+(S), Q \in A_+(\Gamma)$, the subsets of nonnegative lower semi-analytic functions. By Prop.~\ref{prp-F1} and the definition of $F_\k( \cdot \, ; J)$, we see that in each of the (D)(N)(P) cases, $F_\k( \cdot \, ; J)$ maps the sets $A_b(\Gamma)$, $A_-(\Gamma)$ and $A_+(\Gamma)$ into themselves, for $J \in  A_b(S)$, $J \in A_-(S)$ and $J \in A_+(S)$, respectively.

For discrete spaces and abstract DP problems, where measurability is not a concern, we have considered in our earlier work \cite{BerY-discount,YuB-ssp,BerY-adp} mappings of the form $F_\k$, $\k = (\mu, S)$, without splitting the state space into two parts, $B$ and $S \setminus B$, according to the policy $\mu$. 
In the present context, however, in order for $F_\k$ to map lower semi-analytic functions to lower semi-analytic functions, it is important to introduce $B$ as a parameter component in defining $F_\k$.

\vspace*{-0.2cm}
\subsubsection*{Some Basic Properties of $F_\k$}

We now discuss a few basic properties of the mappings $F_\k$ and $F_\k(\cdot\,; J)$, relating to monotonicity and fixed point properties, and their relation 
with $(J^*, Q^*)$.
Let $F^n_\k(\cdot\,; J)$ denote the $n$-fold composition of $F_\k(\cdot\,; J)$, i.e.,
 $$  F^n_\k(Q\,; J) = \underset{n \text{\ times}}{\underbrace{ F_\k \big(  \cdots F_\k \big(F_\k }} (Q\,; J) \, ; J \big) \cdots \, ; J \big).$$
From its definition $F_\k$ is monotone: $F_\k(Q\,; J) \geq F_\k(Q'\,; J')$ if $J \geq J'$ and $Q \geq Q'$. 
Hence
\begin{equation}
    J \geq J', \ Q \geq Q'  \quad \Longrightarrow  \quad F^n_\k(Q\,; J) \geq F^n_\k(Q'\,; J'), \quad \forall \, n \geq 1. \label{eq-Fmono}
\end{equation}
We consider the pointwise limit 
$$Q_{\k, J} = \lim_{n \to \infty}  F^n_\k(\0\,; J),$$
the existence of which will be shown as part of the following proposition.

 \smallskip
\begin{prop} \label{prp-optstop-cost}
{\rm (D)(N)(P)} \ Let $J \in A_b(S)$ for (D), $J \in A_-(S)$ for (N), and $J \in A_+(S)$ for (P). Then $Q_{\k, J} = \lim_{n \to \infty}  F^n_\k(\0\,; J)$ is well-defined, is lower semi-analytic, and satisfies
\begin{equation} \label{eq-optstop-cost}
    Q_{\k, J} = F_\k(Q_{\k,J} \, ; J).
 \end{equation}  
For (D), it is the only solution of $Q = F(Q\,; J)$ in $A_b(\Gamma)$. 
\end{prop}

\begin{proof}
In case (D), $F_\k(\cdot\,; J)$ is a contraction mapping on the complete metric space $A_b(\Gamma)$: for any $Q, Q' \in A_b(\Gamma)$, we have
$$ \big|  \min \big\{ J(x) \, ,  Q(x, u) \big\} - \min \big\{ J(x) \, ,  Q'(x, u) \big\}  \big| \leq   \|Q - Q' \|_\infty, \qquad \forall \, (x,u) \in \Gamma,$$
and hence by direct calculation,
$\| F_\k(Q\,; J) - F_\k(Q'\,; J) \|_\infty \leq \alpha \|Q - Q' \|_\infty.$
It then follows from Banach's contraction principle \cite[p.\ 220]{Rudin-ana} that $F^n_\k(\0\,; J)$ converges to $Q_{\k, J}$, the unique solution of $Q = F(Q\,; J), Q \in A_b(\Gamma)$.

For case (N)(resp.\ (P)), $Q_{\k,J}$ is the pointwise limit of a sequence of nonincreasing nonpositive functions (resp.\ nondecreasing nonnegative functions). Equation~(\ref{eq-optstop-cost}) then follows by the definition of $F_{\k}(\cdot\,; J)$ and the monotone convergence theorem. That $Q_{\k,J} \in A(\Gamma)$ follows from \cite[Lemma 7.30(2)]{bs}.
\end{proof}
%\mysmallskip

We can relate $F_\k$ and $Q_{\k,J}$ to $Q^*$ as follows.

\mysmallskip
 \begin{prop} \label{prp-propertyF}
 {\rm (D)(N)(P)} \ Let $\k \in \K$, $J \in A(S)$, $Q \in A(\Gamma)$. 
 \moveup%\moveup
 \begin{enumerate}
 \item[(a)]  $F_\k(Q^*; J^*) = Q^*$. \moveup
 \item[(b)] If $J \geq J^*, Q \geq Q^*$, 
 then $F^n_\k(Q ; J) \geq Q^*$ for all $n \geq 1$.\moveup
 \item[(c)]  Let $J$ be as in Prop.~\ref{prp-optstop-cost} with $J \geq J^*$. Then $Q_{\k, J} \geq Q_{\k,J^*}$.\moveup
 \item[(d)] $Q_{\k,J^*} = Q^*$.
 \end{enumerate}
 \end{prop} 
 
\begin{proof}
Let $\k = (\mu, B)$. Since $J^*(x) = \inf_{v\in U(x)} Q^*(x,v) \leq Q^*(x, u)$ for all $(x, u) \in \Gamma$, we can rewrite the iterated integral in the sum (\ref{def-F}) defining $F_\k (Q^*; J^*)(x,u)$ as
$$  \int_B \int_C \min \big\{ J^*(x') \, , \, Q^*(x', u') \big\} \, \mu(du' \!\mid x') \, q(dx' \!\mid x, u) =  \int_B J^*(x') \, q(dx' \!\mid x, u).$$
Combining it with the second term in (\ref{def-F}), we obtain 
$$F_\k(Q^*\,;J^*)(x,u) = g(x,u) + \alpha \int_S J^*(x') \, q(dx' \!\mid x, u) = Q^*(x, u), \qquad \forall \, (x,u) \in \Gamma.$$
This proves part (a). Part (b) then follows from part (a) and the monotonicity of $F_\k$ (cf.\ Eq.~(\ref{eq-Fmono})).

For part (c), since $J \geq J^*$, we have $F^n_\k (\0\,; J) \geq F^n_\k(\0\,; J^*)$ for every $n$, by the monotonicity of $F_\k$ (cf.\ Eq.~(\ref{eq-Fmono})).
Then by Prop.~\ref{prp-optstop-cost},
$Q_{\k, J} = \lim_{n \to \infty} F^n_\k (\0\,; J) \geq \lim_{n \to \infty} F^n_\k (\0\,; J^*) = Q_{\k,J^*}.$

There remains to show $Q_{\k,J^*} = Q^*$ in part (d).
For case (D), this is true because by part (a), $Q^*$ is the solution of $F_\k(Q\,; J^*) = Q, Q \in A_b(\Gamma)$, 
whereas this equation has $Q_{\k,J^*}$ as its unique solution by Prop.~\ref{prp-optstop-cost}(D).
For case (N), we have $J^* \leq 0$ and consequently, $F_\k(\0\,; J^*) = Q^*$ by the definitions of $F_\k$ and $Q^*$. In view of part (a), this implies $F^n_\k (\0\,; J^*) = Q^*$ for every $n$, and hence 
$Q_{\k,J^*} = Q^*$ by Prop.~\ref{prp-optstop-cost}.
For case (P), we will show that $Q_{\k,J^*} = Q^*$ as Prop.~\ref{lma-P-optstop} in Appendix~\ref{appsec2} (the proof is not as simple as in (D)(N)).
 \end{proof}
%\mysmallskip 

\subsubsection*{Relation to Optimal Stopping Problems}
In the mixed value and policy iteration algorithms, we will use mappings $F_\k(\cdot\,; J)$, $\k \in \K, J \in A(S)$, and their fixed points $Q_{\k,J}$ in a step that is analogous to policy evaluation. Let us first discuss the relations of these objects with certain optimal stopping problems and with the original control problem, to provide more insights and intuition. 

\mysmallskip  
\begin{rem}[Optimal stopping problem corresponding to $F_\k(\cdot\,; J)$] \label{rmk-optstop1}
The form of the mapping $F_\k(\cdot \,; J)$ (cf.\ Eq.~(\ref{def-F})) suggests a connection to an optimal stopping problem defined by $\k =(\mu, B)$, $J$, and the parameters of the original control problem, with $J$ specifying the stopping costs.
We give a precise mathematical formulation in Appendix~\ref{appsec-optstop}, where we will also show that $F_\k(\cdot \,; J)$ can be viewed as a form of the optimal cost operator (Lemma~\ref{lma-stopmodel1}).
Here we give an intuitive description (see Figure~\ref{fig-stop} for an illustration).
In the optimal stopping problem associated with $\k=(\mu, B)$ and $J$, the states are the state-control pairs of the original control problem. 
Suppose we start from a state $(x,u)$ in $\Gamma$ at time $0$. At this time we must pay $g(x,u)$ and choose to continue (this corresponds to the first term in Eq.~(\ref{def-F})).
At time $1$, we first land at $x'$ according to $q(dx' \mid x, u)$. 
If $x' \in S \setminus B$, then we must pay $J(x')$ and immediately stop (this corresponds to the second term in Eq.~(\ref{def-F})).
If $x' \in B$, then $u'$ is generated according to $\mu(du' \mid x')$ and we land at $(x', u')$.
There, we can either stop and pay $J(x')$, or continue with the continuation cost $g(x', u')$ and then repeat the process just described for time $1$. (This corresponds to the third term in Eq.~(\ref{def-F}), and the minimization inside that term reflects the two choices, to stop or to continue.) 
Because of the correspondence between $F_\k(\cdot\,; J)$ and an optimal stopping problem, some of the theory for (D)(N)(P) with a finite number of controls can be applied to analyze the properties of $F_\k$ (see Appendices~\ref{appsec-optstop} and~\ref{appsec2}). \qed
\end{rem}

\begin{figure}%[thb] 
   \centering
   \includegraphics[width=3.8in]{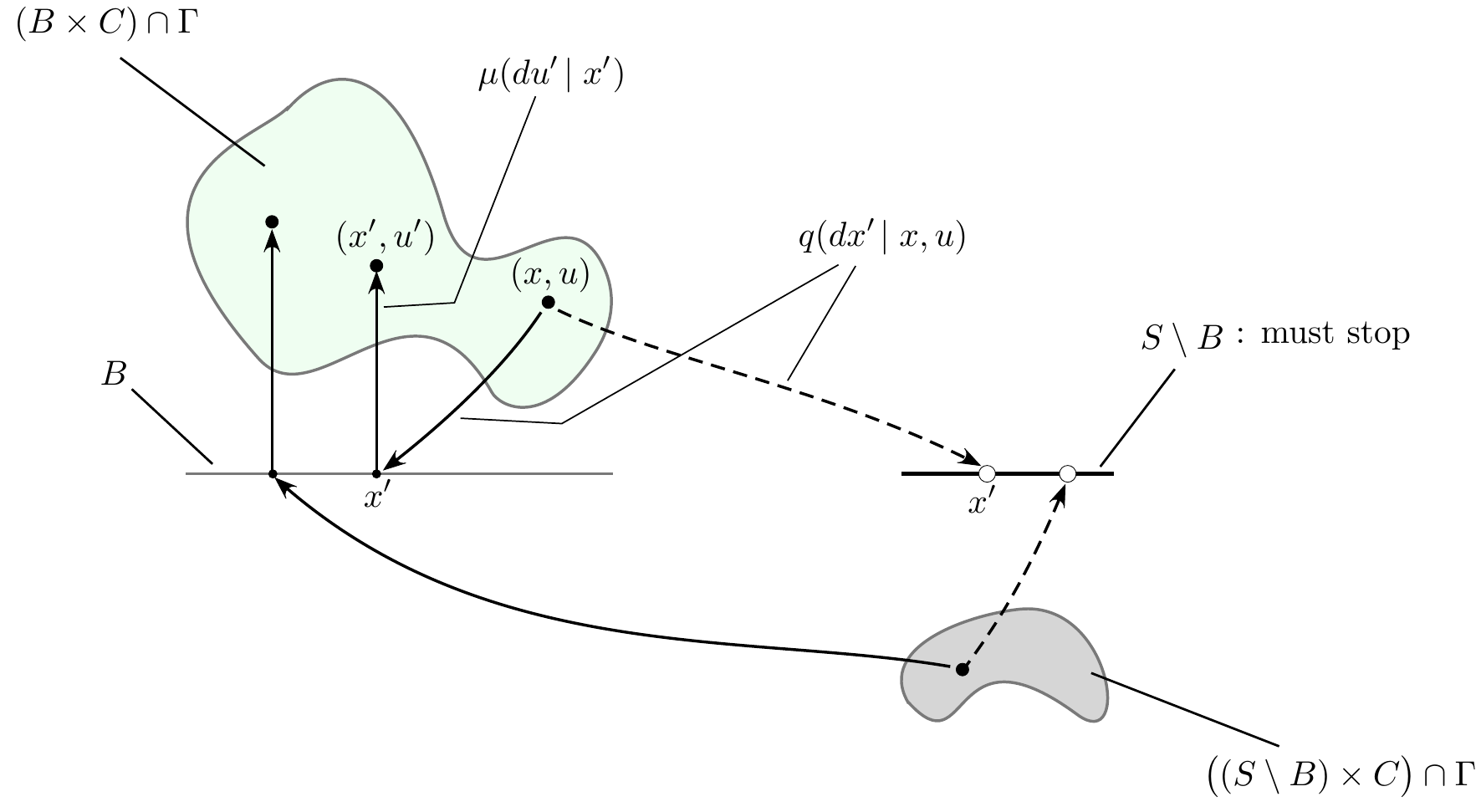} 
   \caption{Illustration of the system dynamics of an optimal stopping problem corresponding to $F_\k (\cdot\,; J)$ with $\k = (\mu, B) \in \K$.}
   \label{fig-stop}
\end{figure} 

\mysmallskip
\begin{rem}[Interpreting $Q_{\k,J}$ in the original control problem] \label{rmk-optstop2}
The function $Q_{\k, J}$ is the optimal cost function of the optimal stopping problem associated with $J$ and $\k=(\mu,B)$ mentioned in the preceding remark (see Cor.~\ref{cor-stopmodel3}, Appendix~\ref{appsec-optstop}). We discuss a special case of the stopping costs $J$. It allows us to relate $Q_{\k, J}$ to the total costs in the original problem. 

Suppose $J = J_\pi$ for some $\pi = (\pi_0, \pi_1, \ldots) \in \Pi$ in the original problem. 
Then for each $(x,u) \in \Gamma$, we may interpret $Q_{\k, J_\pi}(x,u)$ as the minimal cost over the subset of policies for the original problem which switch from
$\mu$ to $\pi$ at most once according to the following restrictions.
At time $0$ we must apply control $u$ at state $x$. From time $1$ on, we can either follow the stationary policy $\mu$ or use the policy $\pi$, which we must do if the state goes outside the set $B$. Once we start to use $\pi$ at time $\tau$, say, we must apply the randomized control rules $\pi_k(\cdot \mid x_{\tau}, u_{\tau}, x_{\tau+1}, u_{\tau+1}, \ldots, x_{\tau+k})$ for time $\tau+k, k \geq 0$, and continue in this way forever. 
In other words, if we can only apply policy $\mu$ or $\pi$ in the original problem in the manner just described, then $Q_{\k, J_\pi}(x,u)$ equals the minimal cost we can achieve by switching from following $\mu$ to following $\pi$ in an optimal way. 
Furthermore, an optimal or near-optimal time to switch to $\pi$ can be found from an optimal or near-optimal solution of the associated optimal stopping problem, if we interpret the action to stop as the decision to switch policies in the original problem.

We do not include a formal proof for the above interpretation of $Q_{\k, J_\pi}(x,u)$ in the paper; but we note that it is similar to the analysis we give in Appendix~\ref{appsec2}, and that the fact $Q_{\k,J^*} = Q^*$ (Prop.~\ref{prp-propertyF}(d)) is also a manifestation of this interpretation.
\qed
\end{rem}
%\smallskip

\subsection{Algorithms} \label{sec3.2}
We give first our mixed value and policy iteration algorithm in its basic form. The conditions needed for the convergence of this form of the algorithm are different for each of the (D)(N)(P) cases, and will be given in the subsequent Sections~\ref{sec4} and~\ref{sec5}.  
Our algorithm starts with a pair $(J_0, Q_0)$, which depending on whether case (D), (N), or (P) holds, must belong to $A_b(S) \times A_b(\Gamma)$, or $A_-(S) \times A_-(\Gamma)$, or $A_+(S) \times A_+(\Gamma)$, respectively. 
It then uses a sequence $\k_k = (\mu_k, B_k) \in \K$ to calculate iteratively functions $(J_k, Q_k) \in A(S) \times A(\Gamma)$, $k \geq 1$, as follows, with each iteration consisting of two steps, in analogy with the policy evaluation and improvement steps in standard policy iteration. 
The stationary policy $\mu_k$ and the Borel subset $B_k \subset S$ of each pair $\theta_k$ can be arbitrary here, but specific choices leading to policy iteration-like algorithms will be demonstrated in the subsequent two example algorithms.

\smallskip
{\samepage 
 \noindent {\bf Algorithm I ($k$th iteration; basic form):}
\moveup\moveup
\begin{itemize}
 \item Choose $\k_k = (\mu_k, B_k) \in \K$, 
 let
 \begin{equation}
      Q_{k+1} = F^{n_k}_{\k_k}(Q_k ; J_k)  \ \ \ \text{for some} \ n_k \geq 1, \quad  \text{or} \ \ Q_{k+1} = Q_{\k_k, J_k}, \qquad \quad \label{alg-1a}
 \end{equation}
 and let 
 \begin{equation}
 J_{k+1} = M (Q_{k+1}). \label{alg-1b}
 \end{equation}
 \end{itemize}
 }\smallskip

The fact that for all $k$, $(J_k,Q_k) \in A(S) \times A(\Gamma)$ can be seen from the inductive argument: by Props.~\ref{prp-F1} and \ref{prp-optstop-cost}, the function $Q_{k+1}$ is lower semi-analytic if $J_k, Q_k$ are lower semi-analytic, whereas the infimization (\ref{alg-1b}) results in a lower semi-analytic function $J_{k+1}$ by~\cite[Prop.\ 7.47]{bs}. 

The algorithm (\ref{alg-1a})-(\ref{alg-1b}) allows any choice of $\k_k = (\mu_k, B_k) \in \K$.
If we let $\k_k = (\mu_k, \emptyset)$ for each iteration $k$, the policy $\mu_k$ has no effect on the iterates, and the algorithm reduces to value iteration $J_{k+1} = T(J_k)$.
On the other hand, by Prop.~\ref{prp-F1}(b), we can choose sets $B_k$ that are not only nonempty but also large (cf.\ the subsequent Examples~\ref{ex-choiceB2}, \ref{ex-choiceB1}).
In what follows, we consider choices of $\mu_k$ based on $Q_k$ and a selection theorem of the Jankov-von Neumann type, and we derive policy iteration-like algorithms. 

Let $Q \in A(\Gamma)$ and $\epsilon > 0$. By a selection theorem for lower semi-analytic functions \cite[Prop.\ 7.50(b)]{bs},
we can select a universally measurable, nonrandomized stationary policy 
$\mu$ such that, with
$I = \big\{ x \in S  \, \big| \,  \argmin _{u \in U(x)} Q(x, u) \not=\emptyset \, \big\},$
\begin{align}
      \mu(x)  & \in   \argmin _{u \in U(x)} Q(x, u) \qquad \text{if} \ x \in I,   \label{eq-selmu1} \\ 
         Q(x, \mu(x)) & \leq \begin{cases}
           M(Q)(x) + \epsilon \ \ \  & \text{if}  \ x \not\in I, \ M(Q)(x)  > -\infty, \\
           - 1/\epsilon & \text{if} \ x \not\in I, \  M(Q)(x) = - \infty.
           \end{cases} \label{eq-selmu2}
\end{align}    
If we relax the condition (\ref{eq-selmu1}), then by \cite[Prop.\ 7.50(a)]{bs}, we can find instead an analytically measurable policy $\mu$ such that for all states $x$,
\begin{align}
  Q(x, \mu(x)) & \leq \begin{cases}
           M(Q)(x) + \epsilon \ \ \  & \text{if}   \ M(Q)(x)  > -\infty, \\
           - 1/\epsilon & \text{if}  \  M(Q)(x) = - \infty.
           \end{cases} \label{eq-selmu3}
\end{align}

By choosing the policies in the basic algorithm by ``approximate policy improvement'' based on the above selection theorem and for a desired value of $\epsilon$, we can obtain special forms of the basic algorithm that resemble to some degree standard forms of modified policy iteration.
Here is an example. The algorithm starts with a pair of functions $(J_0,Q_0)$ and a stationary policy $\mu_0$ as in Algorithm I. 
In the process of computing $J_{k+1}$ by partial minimization of $Q_{k+1}$, it also obtains the policy $\mu_{k+1}$ for the next iteration using the selection theorem above.  

 \smallskip
 {\samepage 
 \noindent {\bf Policy Iteration-Like Algorithm II ($k$th iteration):}
 \moveup\moveup
\begin{itemize}
\item Choose a Borel set $B_k \subset S$ so that $\k_k=(\mu_k,B_k) \in \K$, and compute $Q_{k+1}$ by Eq.~(\ref{alg-1a}).\moveup
\item Let $J_{k+1} = M(Q_{k+1})$ as in Eq.~(\ref{alg-1b}), and let $\mu_{k+1}$ be a nonrandomized stationary policy satisfying Eqs.~(\ref{eq-selmu1})-(\ref{eq-selmu2}), or Eq.~(\ref{eq-selmu3}), with $Q = Q_{k+1}$ and a desired value of $\epsilon$.\moveup
 \end{itemize}
 } \smallskip

If there exists at least one Borel measurable policy, we can further specialize the basic algorithm to use Borel measurable $\mu_k$ together with $B_k=S$ for every iteration or whenever this is desirable. 
As an example, we give below a policy iteration-like algorithm with Borel measurable policies. 
When the set $\Gamma$ is Borel, a nonrandomized Borel measurable policy is known to exist under fairly general conditions 
(see Appendix~\ref{appsec-borelpolicy}). Thus algorithms of this kind can be applied to a large class of problems.

\smallskip
{\samepage 
 \noindent {\bf Policy Iteration-Like Algorithm III with Borel Measurable Policies:}\\*[0.04in]
 \hspace*{0.2cm}  Let $\mu_0$ be a Borel measurable stationary policy (assumed to exist).\\*[0.03in]
 \hspace*{0.2cm} Iterate for each $k \geq 0$:\moveup\moveup
\begin{itemize}
\item For $\k_k=(\mu_k, S)$, compute $Q_{k+1}$ by Eq.~(\ref{alg-1a}). 
\moveup
\item Let $J_{k+1} = M(Q_{k+1})$ as in Eq.~(\ref{alg-1b}), and let $\mu_{k+1}'$ be a stationary policy satisfying Eqs.~(\ref{eq-selmu1})-(\ref{eq-selmu2}), or Eq.~(\ref{eq-selmu3}), with $Q = Q_{k+1}$ and a desired value of $\epsilon$.\moveup
\item
Select $p_{k+1} \in \P(S)$ and let $B \subset S$ be a Borel set such that $p_{k+1}(B) = 1$ and $(\mu_{k+1}', B) \in \K$ (cf.\ Prop.~\ref{prp-F1}(b)).
Define a Borel measurable policy $\mu_{k+1}$ by
\begin{equation} \label{eq-alg3}
\mu_{k+1}(du \!\mid x) = \begin{cases}
  \mu'_{k+1}(d u \!\mid x) & \text{on}  \ B, \\
  \bar \mu(du \!\mid x) & \text{on}  \ S \setminus B,
  \end{cases}
\end{equation}
where $\bar \mu$ is some Borel measurable stationary policy. (In particular, if $\bar \mu$ can be chosen to be nonrandomized,
then every $\mu_k, k \geq 1$, is a nonrandomized Borel measurable policy.)
%\moveup\moveup
\end{itemize}
} 
%\medskip

\subsubsection*{Comparison with Classical Value and Policy Iteration}

The convergence of the preceding algorithms will be the subject of the next two sections. Here, for a better understanding of the theoretical and computational properties of these algorithms, let us compare them to the classical value iteration and policy iteration.
First, let us try to relate informally the structure of the preceding Algorithms I-III to the classical algorithms, leaving aside for the moment 
the measurability and other mathematical issues that may prevent policy iteration from being carried out.

There are two basic choices at each iteration of the form (\ref{alg-1a})-(\ref{alg-1b}): the pair $(\mu_k,B_k)$ and the integer $n_k$. Both of these choices determine how closely our algorithms resemble value iteration and forms of policy iteration. 
In particular, as noted earlier, the choice $B_k=\emptyset$ makes Algorithms I-II identical to value iteration, regardless of the choice of $n_k$ and $\mu_k$. If we let $n_k=1$ and $Q_{k+1} = F_{\k_k}(Q_k\,; J_k)$ always, then regardless of the choice of $\mu_k$ and $B_k$, after time $0$, Algorithms I-III also become identical to value iteration.

To show the similarities and differences between Algorithms I-III and policy iteration, let us rewrite the latter in terms of mappings $F_\k$ and functions on $\Gamma$. 
Denote the constant function $+\infty$ simply by $+\infty$. The modified or exact policy iteration algorithm (cf.\ Section~\ref{sec2.3}) can be equivalently stated as follows.%footnote starts
\footnote{Here we look at the functions on the space $\Gamma$ produced by standard policy iteration just before the minimization in the policy improvement step takes place, and we write the policy evaluation and improvement steps equivalently in terms of these functions. Specifically, the modified policy iteration calculates the next cost function iterate $V'$ and policy $\mu'$ with $V' = T^n_\mu(V), T_{\mu'}(V') = T(V')$ for some integer $n \geq 1$. 
The improvement step, finding $\mu'$ with $T_{\mu'}(V') = T(V')$ (assuming such $\mu'$ exists), is to solve the partial minimization problem 
$ \mu'(x) \in \argmin_{u \in U(x)} Q'(x, u), x \in S$, 
for the function $Q'$ obtained by policy evaluation:
$Q'(x,u) = g(x,u) + \int_{S} T^n_\mu (V)(x') \, q(dx' \mid x, u)$, $(x, u) \in \Gamma.$
Expressed in terms of $F_{(\mu,S)}$, we have $Q' = F_{(\mu, S)}^{n+1}(Q \,; +\infty)$ for an initial function $Q$ with $Q(x, \mu(x)) = V(x), x \in S$. 
We also have from the definition of $\mu'$ that $Q'(x, \mu'(x)) = T_{\mu'}(V')(x), x \in S$. 
Putting these relations for each iteration of the algorithm together, we obtain the alternative description of modified policy iteration given in the discussion, where the integer $n_k$ corresponds to the integer $n$ here for iteration $k \geq 1$ and to $n+1$ for iteration $k=0$.
The equivalent form for the exact policy iteration is derived similarly.} 
%footnote ends 
Iterate for each $k \geq 0$ and a given nonrandomized stationary policy $\mu_k$:
\moveup\moveup
\begin{itemize}
\item[(i)] (Evaluation step) Let
$Q_{k+1} = F_{(\mu_k, S)}^{n_k}(Q_k \,; \mathbf{+\infty})$ for some $n_k \geq 1$ (in the case of modified policy iteration) or let $Q_{k+1} = Q_{(\mu_k, S), +\infty}$ (in the case of exact policy iteration).\moveup%\moveup 
\item[(ii)] (Improvement step) Let $\mu_{k+1}$ be such that
$\mu_{k+1}(x) \in \argmin_{u \in U(x)} Q_{k+1}(x, u)$ for all $x \in S.$\moveup
\end{itemize}
We now compare this description of standard policy iteration to Algorithms I-III, and make a few observations regarding its similarities, and its theoretical and computational differences from our methods.

\mysmallskip
\begin{rem}[Relation to standard policy iteration] \label{rmk-cmp-pi}
From the above description, we see that as $B_k$ is made larger and approaches $S$, {\it and} $\mu_k$ is chosen by ``approximate policy improvement" based on selection of the form (\ref{eq-selmu1})-(\ref{eq-selmu2}), our algorithms approach the exact form of policy iteration (if $Q_{k+1}=Q_{\k_k,J_k}$) or modified policy iteration (if $Q_{k+1}= F_{\k_k}^{n_k}(Q_k\,; J_k)$). Indeed, suppose that for all $k$, $\mu_k$ can be chosen by exact policy improvement (i.e., $I=S$ in Eq.~(\ref{eq-selmu1})) and $B_k$ can be chosen to be $S$. Then our algorithms become \emph{almost identical} to standard policy iteration, except for a difference in the evaluation step where the ``stopping costs'' are set to be $J_k$ instead of $+\infty$.   

Thus, if for all $k$, $J_k$ lies above the functions against which it is compared in our algorithms so that it acts like $+\infty$, our algorithms become identical to the standard policy iteration algorithms, under the preceding assumptions on $(\mu_k, B_k)$. As an example, under (D), the conditions on $J_k$ just mentioned are satisfied if $\mu_0$ is nonrandomized and the initial condition $J_0\geq J_{\mu_0}$, $J_0 \geq V_0 \geq T_{\mu_0}(V_0)$, where $V_0 = Q_0(\cdot, \mu_0(\cdot))$, holds; then, implemented with $Q_{k+1}=Q_{\k_k,J_k}$ or $Q_{k+1}= F_{\k_k}^{n_k}(Q_k\,; J_k)$, our method becomes identical to the classical policy iteration or modified policy iteration, respectively.

Being equivalent to standard policy iteration is, however, not always the goal of our algorithms. Indeed, in some cases of (N)(P), they can never be the same because of the sequence $\{J_k\}$ used by our algorithms. For example, under (N) exact policy iteration may cycle even for discrete-spaces problems \cite[Example 7.2.4, p.\ 296]{Puterman94}, and for convergence under (P), exact or modified policy iteration need conditions on the initial policy which may be unsatisfiable in some problems (cf.\ Cor.~\ref{cor-valite-P1}, Prop.~\ref{prp-valite-P3} in Section~\ref{secP5.1}). By contrast, our algorithms converge under general initial conditions that guarantee the convergence of value iteration (see Sections~\ref{sec4.2} and~\ref{secP5.2}). This more reliable convergence behavior can be attributed to the use of the ``stopping costs'' sequence $\{J_k\}$ in the evaluation step, which brings a value-iteration character to our algorithms. (Such advantages of our method have also been shown in an asynchronous distributed setting for discounted and total cost finite-spaces problems and abstract DP problems; see \cite{BerY-discount,YuB-ssp,BerY-adp}.)

Finally, since with nonrandomized policies $\mu_k$, our algorithms have similar forms as standard policy iteration, they have comparable per-iteration computation overhead. In particular, compared to exact policy iteration, computing $Q_{\k_k,J_k}$ involves solving an optimal stopping problem that has at most two controls at each state. Compared to modified policy iteration, computing $F_{\k_k}^{n_k}(Q_k\,; J_k)$ involves only an extra comparison of the functions $J_k$ and $Q_k$ before computing the expected value in the mapping $F_{\theta_k}$.  
\qed
\end{rem}
%\mysmallskip

As noted in the preceding discussion, when exact policy improvement can take place, the choice of $(\mu_k, B_k)$ and the use of $\{J_k\}$ in the evaluation step may cause differences between our algorithms and standard policy iteration. The use of $\{B_k\}$ lets our method overcome measurability issues, whereas the use of $\{J_k\}$ gives our method better convergence properties.
We now discuss further the high-level differences between our method and standard policy iteration.

\mysmallskip
\begin{rem}[Differences from standard policy iteration] \label{rmk-3.2}
We first note a generic difference that has to do with approximate policy improvement and can be related to a notion of almost-sure $\epsilon$-optimality. 
For simplicity, we consider Algorithm III and assume that it generates nonrandomized policies $\mu_k$ and sets $Q_{k+1} = Q_{\k_k, J_k}$ for all $k$, where $\k_k=(\mu_k,S)$. The mappings $F_{\k_k}$ take the simpler form
\begin{equation} \label{eq-FmuS-det}
      F_{(\mu,S)} \big(Q\,; J)(x,u) = g(x,u) + \alpha \int_{S} \min \big \{ J(x') \,, Q \big(x', \mu(x') \big) \big\}  \, q(dx' \!\mid x, u), \quad (x, u) \in \Gamma,
\end{equation} 
for the nonrandomized Borel measurable policies $\mu=\mu_k$.
Denote
$$ V_k(x) = \min \big \{ J_k(x) \,, Q_{\k_k, J_k} \big( x, \mu_k(x) \big) \big\}, \qquad x \in S.$$ 
From  Eq.~(\ref{eq-FmuS-det}) and the relation $Q_{\k_k, J_k} = F_{(\mu_k, S)}(Q_{\k_k, J_k}; J_k)$ (Prop.~\ref{prp-optstop-cost}), we see that for all $x \in S$,
\begin{align*}
   J_{k+1}(x) = M \big(Q_{\k_k, J_k} \big)(x) & = \inf_{u \in U(x)}  \left\{ g(x, u) + \alpha \int_{S} V_k(x')  \, q(dx' \!\mid x, u) \right\} = T (V_k)(x), \\
    Q_{\k_k, J_k}\big(x, \mu_{k+1}(x) \big) & = T_{\mu_{k+1}} (V_k) (x).
\end{align*} 
Assume $T (V_k)(x) > - \infty$ for all $x$, for simplicity. Since $\mu_{k+1}$ is chosen based on either Eqs.~(\ref{eq-selmu1})-(\ref{eq-selmu2}) or Eq.~(\ref{eq-selmu3}) (cf.\ the definition~(\ref{eq-alg3}) of $\mu_{k+1}$), it follows that for $k \geq 0$,
\begin{equation} \label{eq-rmk-pe}
   p_{k+1} \Big( \big\{ x \in S \mid T_{\mu_{k+1}} (V_k) (x) \leq T (V_k) (x) + \epsilon \big\} \Big) = 1,
\end{equation}   
where $p_{k+1}$ is the probability measure in Algorithm III. In other words, the generated policy $\mu_{k+1}$ is $\epsilon$-optimal for the optimization problem $T(V_k)(x)$ except on a $p$-null set of states with $p=p_{k+1}$---a property that bears similarity to the notion of ``$(p,\epsilon)$-optimal'' policies \cite{Blk-discount,Str-negative}.

Equation~(\ref{eq-rmk-pe}) shows that Algorithm III includes as a special case a generalized version of policy iteration, which involves $(p, \epsilon)$-optimal policy improvement as defined by Eq.~(\ref{eq-rmk-pe}), and is able to deal with situations where policy improvement cannot be exact either because in Eq.~(\ref{eq-selmu1}) $I\neq S$ or because of Borel measurability issues. By contrast, standard policy iteration does not allow such generalized forms of policy improvement.

Algorithm III also demonstrates important theoretical differences between our policy iteration-like algorithms and standard policy iteration, despite the many similarities the two share in their forms. Like Algorithms I-II, Algorithm III yields reliably $J^*$ in the limit in cases (D)(N)(P) under certain initial conditions (see Sections~\ref{sec4} and~\ref{sec5}), although it operates with Borel measurable policies only. By comparison, if $J^*$ is not Borel measurable, we cannot obtain $J^*$ by policy iteration or modified policy iteration operating with Borel measurable policies, since these algorithms keep the iterates $J_k$ in the set of Borel measurable functions. For an example, see \cite[Example 2]{Blk-discount} or \cite[Example 4.1]{Str-negative}.  
It is also impossible to obtain $J^*$ by exact policy iteration if for some state, there exists no \emph{stationary}, $\epsilon$-optimal policy (even if policy improvement can be carried out exactly). This can happen even in countable-state and control problems; see e.g., \cite[Example 6.1]{Str-negative}. Thus our algorithms overcome some fundamental theoretical limitations of standard policy iteration.
\qed
\end{rem}
%\mysmallskip

For many discrete-spaces problems, policy iteration has been preferred to value iteration because of its computational efficiency \cite{Puterman94}.
Let us discuss now the potential computational benefits of using our policy iteration-like algorithms in place of value iteration.

%\mysmallskip
\begin{rem}[Computational efficiency and the role of the parameter $n_k$]
We noted earlier that if $n_k = 1$ always, our algorithms are essentially identical to value iteration. 
The choice of $n_k$ in our algorithms plays a similar role as in modified policy iteration, and determines how frequently the high overhead minimization in the policy improvement step~(\ref{alg-1b}) is performed. In practice, values of $n_k$ that are larger than 1 are recommended in order to keep this overhead low, and experience has shown that modified policy iteration with well-designed choice of $n_k$ is computationally more efficient than both value iteration and exact policy iteration in finite spaces (see e.g., \cite[Chap.\ 6.5]{Puterman94}). An intuition behind this is that when the generated policies $\mu_k$ are close to an optimal policy, each relatively cheap iteration $T_{\mu_k}(V)$ in the evaluation step can act like the far more expensive value iteration $T(V)$. Based on their similarities to modified policy iteration, we expect our policy iteration-like algorithms, with a judicious choice of $n_k$ and other parameters, to have similar computational advantages over value iteration. This is also consistent with our observations from computational experimentation with discounted finite-spaces problems~\cite{BerY-discount}. For undiscounted problems under (P) or (N) and infinite-spaces problems in general, our algorithms are new and still need to be tested computationally, however.
\qed
\end{rem}

\vspace*{-0.15cm}
\subsubsection*{Choosing the Set $B$}

In Algorithms I-III, we repeatedly find, for a universally measurable policy $\mu$, a Borel set $B \subset S$ such that $\mu(du\!\mid \cdot)$ restricted to $B$ is Borel measurable.
As mentioned earlier, it is desirable to have a ``large'' set $B$ so that a large portion of the policy can be taken into account in the algorithms.
We may measure the ``largeness'' of $B$ with respect to a chosen probability measure $p$ on $S$ (cf.\ Prop.~\ref{prp-F1}(b)). 
The question is then how to choose the measure $p$ in a way coupled to the policy $\mu$. Let us discuss two natural possibilities. The first example describes a fortunate situation where there is a common $p$ suitable for all stationary policies. The second example constructs $p$ based on the policy $\mu$.

\mysmallskip
\begin{example}[Choice of $B$ when $q(dx'\!\mid x,u)$ has a density] \label{ex-choiceB2}
This example was suggested by an anonymous reviewer. 
Suppose that the state transition kernel $q(dx'\!\mid x,u)$ has a density with respect to some $p \in \P(S)$, i.e., it can be expressed as
$$  q(dx'\!\mid x,u) = f(x,u,x') \, p(dx'), \qquad \forall \, (x,u) \in S \times C, $$
for a nonnegative Borel measurable function $f(x,u,x')$.
Let us choose the set $B$ based on $p$, with $\k=(\mu, B) \in \K$ and $p(B) = 1$. Since $p(S \setminus B) = 0$, $q(S \setminus B \!\mid x,u) = 0$ and hence the second term in the definition (\ref{def-F}) of $F_\k( Q \, ; J)(x, u)$ vanishes, giving a simplified expression of $F_\k$ as
\begin{align*}
 F_\k ( Q \, ; J)(x, u)  & =  g(x, u) + \alpha \int_B  \int_C  \min \big\{ J(x') \, , \, Q(x', u') \big\} \, \mu(du' \!\mid x') \, q(dx' \!\mid x, u), \quad (x, u) \in \Gamma,
\end{align*}
where the integration over $B$ can also be replaced by the integration over $S$.

It is worth mentioning that in this example, standard policy iteration has no measurability difficulty.
Indeed, for any stationary policy $\mu$, one can find a Borel set $B \subset S$ with $p(B) = 1$ such that restricted to $B$, the universally measurable cost function $J_\mu$ is Borel measurable \cite[Lemma 7.27]{bs}. This implies that the function 
$g(x,u) + \alpha \int_S J_\mu(x') \, q(dx' \!\mid x, u) = g(x,u) + \alpha \int_B J_\mu(x') \, q(dx' \!\mid x, u)$
is lower semi-analytic, and hence at a policy improvement step a stationary policy satisfying Eq.~(\ref{eq-polite}) can be found.
\qed
\end{example}

%\mysmallskip
\begin{example}[Choice of $B$ based on the Markov chain induced by $\mu$] \label{ex-choiceB1}
Consider a Markov chain on $(S, \U(S))$ with state transition kernel $\kappa(dx' \!\mid x)$ defined by
$$\kappa(D \!\mid x) = \int_C q(D \!\mid x, u) \, \mu(du \!\mid x), \qquad D \in \U(S),$$
where $q(D \!\mid x, u)$ is the measure of $D$ with respect to the completion of $q(dx' \!\mid x, u)$.%footnote starts
\footnote{By definition $\kappa(dx' \!\mid x)$ has the following two properties, which make it a state transition kernel on $(S, \U(S))$. 
First, for any fixed state $x \in S$, $\kappa(\cdot \!\mid x)$ is a probability measure on $(S, \U(S))$.
Second, for any fixed set $D \in \U(S)$,  by~\cite[Prop.\ 7.46]{bs}, $q(D \!\mid x, u) = \int_S \b1_D(x') q(dx' \!\mid x, u)$ as a function of $(x,u)$ is universally measurable on $S \times C$, and therefore, $\kappa(D \!\mid x)$ as a function of $x$ is universally measurable on $S$ by~\cite[Prop.\ 7.46]{bs}.}
%footnote ends
The $n$-step transition kernels are: $\kappa^0(dx' \!\mid x) = \delta_x (dx')$,
$\kappa^{n}(d x' \!\mid x) = \int_S \kappa^{n-1}(d x' \!\mid y) \, \kappa (dy \!\mid x)$, $n \geq 1.$
Given a probability measure $\rho$ on $(S, \U(S))$ and $\beta \in (0,1)$,
define a probability measure $p$ on $(S, \U(S))$ by
$$ p(D) = ( 1 - \beta) \sum_{n=0}^\infty \beta^n \int_S \kappa^n(D \!\mid x) \, \rho(dx), \qquad D \in \U(S).$$
The measure $p$ reflects which sets of states are visited with positive probability under the policy $\mu$ if the initial distribution is $\rho$. This choice of $p$ was suggested by S.~Shreve.

Let us choose $B$ such that $\k=(\mu, B) \in \K$ and $p(B)=1$. Using the definition of $p$, it can be shown that $B$ contains a nonempty absorbing set of states, $\bar D = \{ x \in B \mid \kappa^{n}(S \setminus B \!\mid x) = 0, \, n \geq 1\}$; i.e., $\kappa(\bar D \!\mid x) = 1$ for $x \in \bar D$ (the proof is similar to the proof of \cite[Prop.\ 4.2.3(ii)]{MeT09}). It can also be shown that $p(\bar D) = 1$. To make $\bar D$ ``large,'' we can let $\rho$ have a large support. A nice case is when the state space $S$ is indecomposable for the above Markov chain (i.e., there do not exist two disjoint absorbing sets). Then $\bar D$ is indeed large regardless of our choice of the initial distribution $\rho$.%footnote starts
\footnote{In particular, $S$ is indecomposable if the Markov chain is $\psi$-irreducible (with $\psi$ being a maximal irreducibility measure). Then $\psi$ is absolutely continuous with respect to the probability measure $p$ given above \cite[Prop. 4.2.1(ii)]{MeT09}, and moreover, $p$ is itself a maximal irreducibility probability measure of the Markov chain (so $p$, $\psi$ are mutually absolutely continuous) if the initial distribution $\rho$ is an irreducibility measure~\cite[Prop.\ 4.2.2(iv)]{MeT09}.}
%footnote ends

To have a large absorbing set inside $B$ is an appealing property of the preceding way of choosing $p$ and $B$. To elaborate this point, suppose $\mu$ is nonrandomized. Then on the absorbing set $\bar D$, 
$$ F_\k ( Q \, ; J)(x, \mu(x))  =  g(x, \mu(x)) + \alpha \int_{\bar D}  \min \big\{ J(x') \, , \, Q(x', \mu(x')) \big\} \, q(dx' \!\mid x, \mu(x)), \quad \forall \, x \in \bar D.$$
(Because for $x \in \bar D$, $q(B \!\mid x, \mu(x)) = \kappa(B \!\mid x) = 1$ and consequently, the second term in the definition (\ref{def-F}) of $F_\k( Q \, ; J)(x, \mu(x))$ vanishes, similar to Example~\ref{ex-choiceB2}.) The above expression shows that for the state-control pairs $(x,\mu(x)), x \in \bar D$, the evaluation step (\ref{alg-1a}) with the mapping $F_\k(\cdot\,; J)$ acts like standard policy evaluation when $J$ is sufficiently large (cf.\ the alternative description of standard policy evaluation preceding Remark~\ref{rmk-cmp-pi}). This is desirable if we believe that the policy $\mu$ achieves better costs than the estimated costs $J$, for in that case our choice of $B$ allows the step~(\ref{alg-1a}) to evaluate $\mu$ faithfully on the absorbing set $\bar D$ of the Markov chain induced by $\mu$.
\qed
\end{example}

 \section{Convergence Analysis for Discounted Case (D) and Nonpositive Case (N)} \label{sec4}
\markboth{\rm \S \thesection. Convergence Analysis for Cases (D)(N)}{\rm \S \thesection. Convergence Analysis for Cases (D)(N)}

In this section, we analyze the convergence of the mixed value and policy iteration algorithms given in Section~\ref{sec3.2} for cases (D) and (N). We state convergence results for the basic algorithm (\ref{alg-1a})-(\ref{alg-1b}), since the two other policy iteration-like algorithms are its special cases. 

 \subsection{Discounted Case (D)}
 
In case (D) we consider the sets $A_b(S)$, $A_b(\Gamma)$ of bounded, lower semi-analytic functions.
As mentioned in Section~\ref{sec2.0}, they are closed subsets of the Banach spaces $\M_b(S)$ and $\M_b(\Gamma)$ respectively,
and endowed with the metric $d_{sup}(f, f') = \| f - f' \|_\infty$, the spaces
$(A_b(S), d_{sup})$ and $(A_b(\Gamma), d_{sup})$ are complete.
The algorithm (\ref{alg-1a})-(\ref{alg-1b}) works on the product space
$A_b(S) \times A_b(\Gamma)$ endowed with the metric
$$  d\big( (J,Q) \, ,  (J',Q') \big) =\|(J, Q) - (J', Q') \|_\infty := \max \big\{ \| J - J' \|_\infty \, ,  \| Q - Q' \|_\infty \big\},$$
which is also a complete metric space. 
The convergence results below use a contraction property of $F_\k$ and parallel those given in our earlier work \cite{BerY-discount} for discounted finite-state and control problems.

 \mysmallskip
 \begin{lem} \label{lma-D1} 
 {\rm (D)} \ Let $\k \in \K$, $J, J' \in A_b(S)$, and $Q, Q' \in A_b(\Gamma)$. Then
\begin{align*}
  \| F_\k(Q\, ; J) - F_\k(Q' \,; J')\|_\infty  & \leq \alpha \max \big\{ \| J - J' \|_\infty \, , \, \|Q - Q' \|_\infty \big\}, \\
 \| F_\k(Q \,; J) -  Q^*\|_\infty & \leq \alpha \max \big\{ \| J - J^* \|_\infty \, , \, \|Q - Q^* \|_\infty \big\}, \\
  \|Q_{\k, J} - Q^* \|_\infty & \leq \alpha \| J - J^* \|_\infty.
\end{align*} 
 \end{lem}
 
\begin{proof}
For every $(x, u) \in \Gamma$, 
$$ J(x) \leq J'(x) + \max \big\{ \| J - J' \|_\infty \, , \, \|Q - Q' \|_\infty \big\}, \quad Q(x, u) \leq Q'(x, u) + \max \big\{ \| J - J' \|_\infty \, , \, \|Q - Q' \|_\infty \big\},$$
so
$$ \min \big\{ J(x) \, ,  Q(x, u) \big\} - \min \big\{ J'(x) \, ,  Q'(x, u) \big\}  \leq  \max \big\{ \| J - J' \|_\infty \, , \|Q - Q' \|_\infty \big\}$$
and by symmetry,
$$ \big|  \min \big\{ J(x) \, ,  Q(x, u) \big\} - \min \big\{ J'(x) \, ,  Q'(x, u) \big\}  \big| \leq  \max \big\{ \| J - J' \|_\infty \, ,  \|Q - Q' \|_\infty \big\}.$$
Using the above inequality and the definition of $F_\k$ given in Eq.~(\ref{def-F}), a direct calculation then shows that for each $(x, u) \in \Gamma$,
\begin{align*}
    \Big|  F_\k(Q ; J)(x, u) - F_\k(Q' ; J')(x,u) \Big| & \leq \alpha \, \| J - J' \|_\infty \cdot q (S \setminus \!B \!\mid x, u )  \\
    & \ \ \ + \alpha  \max \big\{ \| J - J' \|_\infty \, , \, \|Q - Q' \|_\infty \big\} \cdot q(B \!\mid x, u)  \\
    & \leq \alpha \max \big\{ \| J - J' \|_\infty \, , \, \|Q - Q' \|_\infty \big\}.
\end{align*}
This proves the first inequality in the lemma, from which the second desired inequality follows by setting $(J', Q') = (J^*, Q^*)$ and using the fact $F_\k(Q^*; J^*) = Q^*$ (Prop.~\ref{prp-propertyF}(a)). To show the third inequality in the lemma, set $Q = Q_{\k,J}$ in the second inequality just proved. Since $Q_{\k,J} = F_\k(Q_{\k,J}\,; J)$ (Prop.~\ref{prp-optstop-cost}(D)),
we obtain
$\| Q_{\k, J} - Q^* \|_\infty \leq \alpha \max \big\{ \| J - J^* \|_\infty \, , \| Q_{\k, J} - Q^* \|_\infty \big\},$
and since $\alpha < 1$, this is equivalent to $\| Q_{\k, J} - Q^* \| \leq \alpha  \| J - J^* \|_\infty$.
\end{proof}

\mysmallskip
\begin{thm} \label{thm-D}
{\rm (D)} \ For any $J_0 \in A_b(S)$ and $Q_0 \in A_b(\Gamma)$, 
the sequence $\big\{(J_k, Q_k)\big\}$ generated by the iteration (\ref{alg-1a})-(\ref{alg-1b}) 
converges to $(J^*, Q^*)$, and
$$ \big\| (J_{k}, Q_{k}) - (J^*, Q^*) \big\|_\infty \leq \alpha^k \, \big\| (J_{0}, Q_{0}) - (J^*, Q^*) \big\|_\infty.$$
\end{thm} 

\begin{proof}
At iteration $k$, either $Q_{k+1} = F_\k^n(Q_k\,; J_k)$ or $Q_{k+1} = Q_{\k,J_k}$ for some $\k \in \K, n \geq 1$.
For the first case, applying the second inequality in Lemma~\ref{lma-D1} $n$ times, we have
$$ \big \| F_\k^n(Q_k\,; J_k) - Q^* \big\|_\infty \leq \alpha \max \big\{ \| J_k - J^* \|_\infty \, , \, \alpha^{n-1} \| Q_k - Q^* \|_\infty \big\}, $$ 
whereas for the second case, $\| Q_{\k,J_k} - Q^* \| \leq \alpha \| J_k - J^* \|_\infty$ by the third inequality in Lemma~\ref{lma-D1}. Thus in either case, 
$$ \| Q_{k+1} - Q^* \|_\infty \leq \alpha \max \big\{ \| J_k - J^* \|_\infty \, , \, \| Q_k - Q^* \|_\infty \big\}.$$
Since $J_{k+1} = M(Q_{k+1}), J^* = M(Q^*)$, and $M$ is nonexpansive, i.e., $\|M(Q) - M(Q')\|_\infty \leq \| Q - Q'\|_\infty$, we have
$$ \| J_{k+1} - J^* \|_\infty = \|M(Q_{k+1}) - M(Q^*) \|_\infty  \leq \alpha \max \big\{ \| J_k - J^* \|_\infty \, , \, \| Q_k - Q^* \|_\infty \big\}.$$
Combining the preceding two inequalities, we obtain
$$\big\| (J_{k+1}, Q_{k+1}) - (J^*, Q^*) \big\|_\infty \leq \alpha^{k+1} \big\| (J_{0}, Q_{0}) - (J^*, Q^*) \big\|_\infty,$$
which is the desired inequality and implies $(J_k,Q_k) \to (J^*, Q^*)$. 
\end{proof}
\mysmallskip

\begin{rem}[Finding near-optimal policies] \label{rem-policyD}
From the iterate sequence $\big\{(J_k,Q_k)\big\}$ generated by the algorithm, we may extract an asymptotically near-optimal sequence of nonrandomized stationary policies $\{\nu_k\}$ as follows. 
Specify $\epsilon > 0$ and let $\nu_k$ be such that $Q_k(x, \nu_k(x)) \leq M (Q_k)(x) + \epsilon$ for all $x \in S$ (cf.\ Eq.~(\ref{eq-selmu3})). 
The policies $\{\mu_k\}$ (resp.\ $\{\mu'_k\}$) generated in the policy iteration-like algorithm~II (resp.\ III), for example, already satisfy these conditions. 
To analyze $J_{\nu_k}$, note that by Theorem~\ref{thm-D}, $\| Q_k - Q^*\|_\infty \leq \alpha^k \Delta$ where $\Delta = \| (J_0, Q_0) - (J^*, Q^*) \|_\infty$. 
Hence for all $x \in S$,
$$ Q^*(x, \nu_k(x)) \leq Q_k(x, \nu_k(x)) + \alpha^k \Delta \leq M (Q^*)(x) + 2 \alpha^k \Delta + \epsilon = J^*(x) + 2 \alpha^k \Delta + \epsilon.$$
Since $Q^*(x, \nu_k(x)) = T_{\nu_k} (J^*)(x)$ (cf.\ Eqs.~(\ref{eq-Q}), (\ref{eq-Tmu})), the preceding inequality shows 
$$ T_{\nu_k} (J^*)(x) \leq J^*(x) +  2 \alpha^k \Delta + \epsilon, \qquad \forall \, x \in S.$$ By the monotonicity and contraction properties of $T_{\nu_k}$, this implies that $\|J_{\nu_k}-J^*\|_\infty \leq (2\alpha^k\Delta + \epsilon)/(1-\alpha)$ and hence
\begin{equation} \label{eq-policyD-asymopt}
    \limsup_{k \to \infty} \| J_{\nu_k} - J^* \|_\infty \leq \epsilon / (1-\alpha).
\end{equation}   
Alternatively, we may choose $\nu_k$ such that 
$\big\|T_{\nu_k}(J_k) - T(J_k) \big\|_\infty \leq \epsilon$ (this is possible by the selection theorem \cite[Prop.\ 7.50]{bs}). 
Then using the contraction property of $T$ and $T_{\nu_k}$, it can be shown (see e.g., \cite[p.\ 45]{Ber-DP13}) that
$\|J_{\nu_k} - J^*\|_\infty \leq \big( \epsilon + 2\alpha \,\|J_k-J^*\|_\infty \big) / (1-\alpha)$ for all $k$, so by Theorem~\ref{thm-D}, Eq.~(\ref{eq-policyD-asymopt}) holds and the sequence $\{\nu_k\}$ is asymptotically $\epsilon/(1-\alpha)$-optimal in that sense. 
 \qed
\end{rem} 
%\mysmallskip
 
 \subsection{Nonpositive Case (N)} \label{sec4.2}
 
In case (N) the one-stage cost function $g \leq 0$ and $J^* \leq 0$, $Q^* \leq 0$. 
The mixed value and policy iteration algorithm (\ref{alg-1a})-(\ref{alg-1b}) operates with nonpositive lower semi-analytic functions in $A_-(S)$ and $A_-(\Gamma)$. We will rely on the monotonicity and fixed point properties of $F_\k$ to ensure its convergence.

First, we derive some simple upper and lower bounds on the iterates generated by the algorithm. 
To simplify notation, let
 \begin{equation}  \label{eq-H}
    H(x, u, J) = g(x, u) + \int_S J(x') \, q(dx' \!\mid x, u), \qquad (x, u) \in \Gamma.
 \end{equation}
Expressed in these terms, $T(J)(x) = \inf_{u \in U(x)} H(x, u, J)$, the optimality equation $J^* = T(J^*)$ is  
$$J^*(x) = \inf_{u \in U(x)} H(x, u, J^*), \qquad x \in S,$$ and by the definition of $Q^*$ (cf.\ Eq.~(\ref{eq-Q})), 
\begin{equation}  \label{eq-QH}
   Q^*(x, u) = H(x, u, J^*), \qquad (x, u) \in \Gamma.
 \end{equation}
 
% \mysmallskip
 \begin{lem} \label{lma-simplebound-DNP} 
 {\rm (D)(N)(P)} \ Let $J \in A(S)$ be as in Prop.~\ref{prp-optstop-cost}, and let $Q \in A(\Gamma)$. Then 
 \begin{align}  \label{eq-upper1}
 F_\k(Q\,; J )(x,u) & \leq H(x, u, J),  & Q_{\k, J}(x, u) & \leq H(x, u, J), \quad \forall \, (x, u) \in \Gamma, \qquad \\
 M \big( F_\k(Q\,; J ) \big) & \leq T(J),  & M \big( Q_{\k, J} \big)  & \leq T(J). \label{eq-upper2} 
 \end{align}
\end{lem} 
 
\begin{proof} 
Since $\min\{ J(x') , Q(x',u') \} \leq J(x')$, 
we have that for all $(x, u) \in \Gamma$, 
\begin{align*}
    F_\k ( Q \, ; J)(x, u) & =  g(x, u) + \alpha \int_{S\setminus B}  J(x') \, q(dx' \!\mid x, u) \\
        &  \ \ \ 
        +  \alpha \int_B  \int_C  \min \big\{ J(x') \, , \, Q(x', u') \big\} \, \mu(du' \!\mid x') \, q(dx' \!\mid x, u) \\
        & \leq  g(x, u) + \alpha \int_{S\setminus B}  J(x') \, q(dx' \!\mid x, u)  +  \alpha \int_B  J(x')  q(dx' \!\mid x, u)  = H(x, u, J).
\end{align*}
From this inequality and Prop.~\ref{prp-optstop-cost},  it follows $Q_{\k,J} = F_\k(Q_{\k,J}\,; J ) \leq H(x, u, J)$. 
This establishes the inequalities in Eq.~(\ref{eq-upper1}). 
Minimizing over $U(x)$ for each $x$ in Eq.~(\ref{eq-upper1}), we obtain Eq.~(\ref{eq-upper2}).
\end{proof}
\mysmallskip
 
 We use the above bounds to upperbound the iterates of the algorithms. The next lemma applies also to (D)(P). For the algorithm that uses the second rule of (\ref{alg-1a}) to set $Q_{k+1} = Q_{\k_k,J_k}$ at some iterations, the second statement of the lemma will rely on Prop.~\ref{prp-propertyF}(d), which in the case (P) will be proved in Appendix~\ref{appsec2} as Prop.~\ref{lma-P-optstop}.
 
 \mysmallskip
 \begin{lem} \label{lma-bound-NP}
 {\rm (N)(P)} \ Let $\big\{(J_k, Q_k)\big\}$ be iterates generated by the iteration (\ref{alg-1a})-(\ref{alg-1b}) with $J_0 \in A_-(S)$, $Q_0 \in A_-(\Gamma)$ in case (N) and with $J_0 \in A_+(S)$, $Q_0 \in A_+(\Gamma)$ in case (P). Then for $k \geq 1$,
 \begin{equation} \label{eq-bound-NP}
     \qquad J_{k} \leq T^k(J_0), \qquad Q_k(x, u) \leq H(x, u, J_{k-1}), \ \ \  \forall \, (x, u) \in \Gamma.
 \end{equation}    
 If $J_0 \geq J^*, Q_0 \geq Q^*$, then we also have $J_k \geq J^*, Q_k \geq Q^*$.
 \end{lem}
 
 \begin{proof}
 For each $k \geq 0$, either $Q_{k+1} = F_\k^n(Q_k\,; J_k)$ or $Q_{k+1} = Q_{\k,J_k}$ for some $\k \in \K, n \geq 1$.
 By Eq.~(\ref{eq-upper1}), the right-hand side inequality for $Q_k$ in Eq.~(\ref{eq-bound-NP}) follows. 
 Since $J_{k+1} = M(Q_{k+1})$, we have, by Eq.~(\ref{eq-upper2}), $J_{k+1} \leq T(J_k)$ for all $k$.
 This implies $J_{k} \leq T^k(J_0)$ by the monotonicity of $T$.
 
 Let $J_0 \geq J^*$ and $Q_0 \geq Q^*$. We show by induction that $J_k \geq J^*, Q_k \geq Q^*$ for every $k$.
 Suppose this is true for some $k \geq 0$. Consider the $(k+1)$-th iterate. If $Q_{k+1} = F^n_\k(Q_k\,; J_k )$, then by the induction hypothesis, the monotonicity of $F_\k$ (cf.\ Eq.~(\ref{eq-Fmono})) and Prop.~\ref{prp-propertyF}(a), we have
$$Q_{k+1} = F^n_\k(Q_k\,; J_k ) \geq F^n_\k(Q^*\,;J^*) = Q^*.$$
If $Q_{k+1} = Q_{\k,J_k}$, then since $J_k \geq J^*$ by the induction hypothesis and $Q_{\k,J_k} \geq Q_{\k, J^*}$ by Prop.~\ref{prp-propertyF}(c), we have $Q_{k+1} \geq Q_{\k,J^*} = Q^*$ by Prop.~\ref{prp-propertyF}(d) (proved as Prop.~\ref{lma-P-optstop} for (P)). Thus in either case, $Q_{k+1} \geq Q^*$, and consequently, $J_{k+1} = M(Q_{k+1}) \geq M(Q^*) = J^*.$ This completes the induction.
  \end{proof}
% \mysmallskip
 
 The relation $J^*\leq J_k\le T^k(J_0)$ in Lemma \ref{lma-bound-NP}, which holds when $J_0\geq J^*$, is the key to our convergence analysis for cases (N) and (P). It implies that our method converges to $J^*$ from above whenever the ordinary value iteration method does. In case (N), we will exploit the generic convergence property of value iteration in the following theorem, whereas in case (P), we will derive sufficient conditions for convergence of value iteration from above in the next section.
 
 \mysmallskip
 \begin{thm} \label{thm-N}
{\rm (N)} \ For any $J_0 \in A_-(S)$ and $Q_0 \in A_-(\Gamma)$ such that $J_0 \geq J^*$ and $Q_0 \geq Q^*$,
the sequence $\big\{(J_k, Q_k)\big\}$ generated by the iteration (\ref{alg-1a})-(\ref{alg-1b}) 
converges to $(J^*, Q^*)$.
\end{thm}

\begin{proof}
We show first $J_k \to J^*$.
We have $J^* \leq J_k \leq T^k(J_0)$ by Lemma~\ref{lma-bound-NP}.  
Since $J^* \leq J_0 \leq 0$ by assumption and $T^k(\0) \downarrow J^*$ under (N),
we have $T^k(J_0) \to J^*$ and hence $J_{k} \to J^*$. 
Then, for each $(x, u) \in \Gamma$, by Fatou's lemma \cite[p.~131]{Dud02} (applied to nonpositive functions), 
$$ \limsup_{k \to \infty} H(x, u, J_k) \leq H \big(x, u, \limsup_{k \to \infty} J_k \big) = H(x, u, J^*) = Q^*(x,u)$$
(cf.\ Eqs.~(\ref{eq-H})-(\ref{eq-QH})).
Since $Q^*(x,u) \leq Q_{k+1}(x,u) \leq H(x,u, J_k)$ by Lemma~\ref{lma-bound-NP}, this implies the convergence 
$Q_{k} \to Q^*$. 
\end{proof} 
\mysmallskip
 
 \begin{rem}
Regarding near-optimal policies in case (N), recall that they are guaranteed to exist among semi-Markov policies, but not necessarily among stationary or Markov policies.
The construction of an $\epsilon$-optimal semi-Markov policy under (N) is much more involved than under (D)(P), and knowing the optimal cost function $J^*$ alone is insufficient (see the proof of~\cite[Prop.\ 9.20]{bs}), even if it were available. Moreover, even if an optimal stationary policy exists, it is possible that a policy $\mu$ satisfies $T_\mu(J^*)=T(J^*)$ without being optimal.%footnote starts
\footnote{As an example, let $S=\{0,1\}$ with state $0$ being cost-free and absorbing. 
At state $1$, there are two controls: control $1$ leads to state $1$ with cost $0$, and control $0$ leads to state $0$ with cost $-1$.
Then $J^*(0)=0$, $J^*(1)=-1$, and the suboptimal policy $\mu$ that makes self-transitions at state $1$ satisfies $T_\mu(J^*)=T(J^*)$.}
% footnote ends
Hence, we do not expect to have simple ways to obtain near-optimal policies from the iterate sequence $\big\{(J_k,Q_k)\big\}$ generated by our algorithm. Intuitively, if we start the algorithm with $J_0 =\0, Q_0=\0$ or cost functions of some policies, it seems possible to construct for each given state, history-dependent or semi-Markov policies that asymptotically become near-optimal for that state, by using the relations between the optimal stopping problems associated with the step (\ref{alg-1a}) and the original problem (cf.\ Remarks~\ref{rmk-optstop1}, \ref{rmk-optstop2}). Due to its complexity, however, we do not discuss this subject in this paper.
\qed
\end{rem}
%\mysmallskip

 \section{Convergence Analysis for Nonnegative Case (P)} \label{sec5}
 \markboth{\rm \S \thesection. Convergence Analysis for Case (P)}{\rm \S \thesection. Convergence Analysis for Case (P)}
 
In this section we consider the case (P) with nonnegative one-stage costs. We first prove a new convergence theorem for value iteration in Section~\ref{secP5.1}. Using this theorem, we then derive in Section~\ref{secP5.2} convergence results for the mixed value and policy iteration algorithms discussed in Section~\ref{sec3.2}, and for another variant algorithm which admits a linear programming implementation for a certain class of problems and thus has computational advantages.

 \subsection{A Convergence Theorem for Value Iteration} \label{secP5.1}
The nonpositive case (P) is more complex than (D)(N). Neither value iteration nor policy iteration are guaranteed to give us $J^*$, even if policy iteration encounters no measurability issues. 
For value iteration, as mentioned in Section~\ref{sec2.2}, for some $J_\infty \in A_+(S)$, we have
$ T^k (\0) \uparrow J_\infty \leq J^*,$
and it is possible that $J_\infty < J^*$. It is known that $J_\infty = J^*$ if $U(x)$ is a finite set for each $x \in S$, or more generally, if a compactness-type condition on the control constraint set holds~\cite[Prop.\ 9.17, Cor.\ 9.17.1]{bs}; but these conditions are restrictive.
For policy iteration, it can happen that for a suboptimal stationary policy $\mu$,
$ J_\mu = T_\mu(J_\mu) = T (J_\mu),$
even in finite-state and control problems,%footnote
\footnote{For a simple example, consider a problem with two states $\{0,1\}$. State $0$ is cost-free and absorbing. State $1$ has two controls $\{0,1\}$: the control $1$ leads to a zero-cost self-transition to state $1$, and the control $0$ leads to state $0$ with cost $1$. Then the nonrandomized stationary policy $\mu$ with $\mu(1)=0$ is suboptimal but satisfies $T_\mu(J_\mu) = T (J_\mu)$. See \cite[Example 7.3.4]{Puterman94} for a similar example. We also note that total cost finite-state and control problems can be solved by using the policy iteration algorithms of Veinott \cite{Vei66} and of Miller and Veinott \cite{MiV69} based on the concept of sensitive optimality (\cite{Vei69}; see also \cite[Sec.\ 10.3]{Puterman94}).} %footnote ends
and the method terminates with the suboptimal policy~$\mu$.

We thus look for ways to mitigate the difficulties.  
Any condition forcing $T^k (\0) \uparrow J^*$, however, seems restrictive, in view of Maitra and Sudderth's result \cite{MS92}.
They showed that $J^*$ can be obtained by applying $T$ a transfinite number of times, starting from the function $J = \0$, and in general, the number of times needed can be uncountably infinite \cite[p.\ 930]{MS92}.
This led us to consider ways to make value iteration converge from above instead of from below, which is also natural when using policy costs, since $J_\mu \geq J^*$. 
We will modify Whittle's bridging condition \cite{Whit79, Har80} to suit our purpose. 

Before proceeding, let us give a simple example to elucidate the behavior of value iteration just discussed. 
The example is from 
\cite[p.\ 215]{bs}. 
In this example $J_\infty < J^* \equiv \infty$. 
We illustrate how value iteration with transfinite recursion is able to obtain $J^*$ in the end, after countably many iterations.
This example falls into a special case analyzed in \cite[Sec.\ 5]{MS92}, which predicted, for a broad class of problems, that the number of iterations required for value iteration to converge from below is at most countably infinite.

\begin{example} \label{ex-P}
The state and control spaces are $S = \{0,1,2, \ldots\}, C = \{1, 2, \ldots\}$, and the control constraint is $U(x) = C$ for every $x \in S$. State transitions are deterministic and uncontrolled except at state $0$: applying control $u$ at state $x$, the successor state is $u$ if $x = 0$ and $x-1$ if $x \geq 1$. The one-stage cost is zero except at state $1$: $g(1, u) = 1$ for all $u$.
Write a function $J$ on $S$ in vector form as $J = \big(J(0), J(1), \ldots \big).$ 
The optimal cost function is 
$J^* = (\infty, \infty, \ldots)$
because under any policy, the system will visit state $1$ infinitely often and accumulate one more unit of cost at each visit.
      
The pointwise limit $J_\infty$ of $\{T^k(\0)\}$ is
$J_\infty = (0, 1, 1, \ldots),$ 
since $T^k(\0) = (0, 1, 1, \ldots, 1, 0, 0, \ldots)$ with $k$ $1$'s followed by all $0$'s.
As in \cite{KreP77}, set $J_{\infty0} = J_\infty$ and initiate value iteration with it.
This gives us $J_{\infty1} = \lim_{k \to \infty} T^k (J_{\infty0})$, which is 
$J_{\infty1} = (1, 2, 2, \ldots).$
Continuing in this way, we define recursively $J_{\infty (m+1)} = \lim_{k \to \infty} T^k (J_{\infty m})$ and we get
$J_{\infty (m+1)} = J_{\infty m} + 1.$
In the end, from the pointwise limit of the nondecreasing sequence $\{ J_{\infty m}\}$ we obtain $J^*$.      
\qed
\end{example} 
\mysmallskip
 
We now proceed to place a condition on the initial function $J_0$ for value iteration $T^k(J_0)$, to ensure the convergence of value iteration (from above, primarily) to $J^*$. This condition, given in the following theorem, is motivated by Whittle's bridging condition \cite{Whit79, Har80} (cf.\ Remark~\ref{rem-P3}) and its appealingly simple form. 
(The paper \cite{Whit79} called $J_0$ the ``terminal function'' instead of ``initial function,'' for the reason that $J_0$ can be viewed as setting the terminal costs for finite horizon problems.)
The implications of our theorem given below are, however, different from Whittle's \cite{Whit79, Har80}, as we will remark shortly. 

\mysmallskip
\begin{thm} \label{thm-valite-P}
{\rm (P)} The following hold:\\*[0.03in]
(a) For any $c > 1$, $T^k(c J^*) \downarrow J^*$.\\*[0.07in]
(b) $T^k (J) \to J^*$ for all $J \in A_+(S)$ such that 
\begin{equation} 
  \underline{J} \leq J \leq c J^*, \quad \text{for some} \ c > 1,  \notag
\end{equation}  
where $\underline{J} \in A_+(S)$ satisfies $\underline{J} \leq J^*, T^k (\underline{J}) \to J^*$. 
In particular, if $T^k(\0) \uparrow J^*$, then $T^k (J) \to J^*$ for all $J \leq c J^*, J \in A_+(S)$.\\*[0.07in]
(c) $J^*$ is the unique fixed point of $T$ within the set $\{ J \in A_+(S) \mid J \leq c J^* \ \text{for some} \ c > 1\}$. 
\end{thm}
\mysmallskip

We note that Theorem~\ref{thm-valite-P}(b)-(c) follows directly from Theorem~\ref{thm-valite-P}(a).
To see this, suppose part (a) is proved. Then under the assumptions of part (b), we have $T^k(\underline{J}) \leq T^k(J) \leq T^k(c J^*)$ by the monotonicity of $T$. Since
$T^k(\underline{J}) \to J^*$ by assumption and $T^k(c J^*) \downarrow J^*$ by part (a), part (b) follows.
For part (c), by \cite[Prop.\ 9.10(P)]{bs} we have the following implication,
$$ J \in A_+(S), \ J=T(J) \qquad \Longrightarrow \qquad J \geq J^*,$$
which together with part (a) implies the conclusion of part (c). 
Thus to prove Theorem~\ref{thm-valite-P}, it suffices to prove its part (a).

Before giving the proof, 
let us make several remarks about the implications of Theorem~\ref{thm-valite-P} and its relation with Whittle's bridging condition.

\mysmallskip
\begin{rem} 
In Theorem~\ref{thm-valite-P}(b), we can always let $\underline{J} = J^*$. Then Theorem~\ref{thm-valite-P}(b) reads as: 
\begin{equation}
   \qquad T^k (J) \to J^*, \qquad \forall \,  J \in A_+(S) \ s.t. \ J^* \leq J \leq c J^*, \ c > 1.
\end{equation} 
Indeed, in view of the result of Maitra and Sudderth~\cite{MS92} and the simple Example~\ref{ex-P}, $J^*$ may be the only function that can serve as an initial function $\underline{J}$ lying below $J^*$ from which VI converges as in Theorem~\ref{thm-valite-P}(b). \qed
\end{rem}

\mysmallskip
\begin{rem} \label{rem-P2}
Theorem~\ref{thm-valite-P}(a)-(b) roughly says that value iteration converges to $J^*$ if the initial function $J$ is ``commensurate'' with $J^*$.
In particular, if $J \geq J^*$, then on the set of states $x$ with finite $J^*(x)$, the shape of $J$ must be ``compatible'' with that of $J^*$, with $J(x) = 0$ whenever $J^*(x) = 0$. 
The theorem also implies that whenever the policy iteration algorithm gets stuck at a suboptimal policy $\mu$ with
$T_\mu (J_\mu) = T (J_\mu)$, $J_\mu$ must have a ``wrong shape'' relative to $J^*$. 

Of course it can be difficult to know even the ``shape'' of $J^*$. In Example~\ref{ex-P}, for instance, $J^* \equiv \infty$, so the only function between $J^*$ and $c J^*$, $c > 1$, is $J^*$ itself. In an example of Strauch~\cite[p.\ 881]{Str-negative} (see also \cite[p.\ 930]{MS92}), $J^*$ takes values in $\{0, 1\}$ and $T^k(\0) \not\to J^*$. The set $\{ x\in S \mid J^*(x) = 0\}$ is rather intricate. If we know this set of states, then with any initial function $J$ that takes the value $0$ on this set and the value $a \geq 1$ elsewhere, value iteration turns out to converge in one iteration in this example (see Appendix~\ref{appsec-Pexample}). 
\qed
\end{rem}

\mysmallskip
\begin{rem} \label{rem-P3}
Whittle's bridging condition is as follows: for some real $c$ and Markov policy $\pi$, either $J_\pi \leq c \, T^n(\0)$ for some $n$, or $J_\pi \leq c J_\infty$ and $J_\infty = T (J_\infty)$.
The condition implies that $J_\infty = J^*$ and $T^k( J) \to J^*$ for all $J \in A_+(S)$ with $J \leq a J^*$ for some real $a$ \cite{Whit79,Har80}. 
A similar, slightly weaker condition leading to the same conclusions is $J^* \leq c \, T^n(\0)$ for some $n$ (personal communication with E.~Feinberg).
The main difference between these results and Theorem~\ref{thm-valite-P} is that Theorem~\ref{thm-valite-P} does not place any condition on the model of the control problem. Instead, it restricts only the initial function for value iteration and it holds for all nonnegative control models. 
If the bridging condition or any other condition for $T^k(\0) \uparrow J^*$ holds, they can be used to set $\underline{J} = \0$ in the theorem, as stated in Theorem~\ref{thm-valite-P}(b). Then, the condition for $J$ becomes $0 \leq J \leq cJ^*$, the same as in \cite{Whit79,Har80}. 
\qed
\end{rem}
\mysmallskip

We now proceed to prove Theorem~\ref{thm-valite-P}(a). 
To this end, we start with two lemmas to characterize the pointwise limit of $\{T^k(c J^*)\}$.
The first lemma below is a basic fact;
the second one is important for our proof.

\mysmallskip
\begin{lem} \label{lma-P1}
If $J \in A_+(S)$ satisfies $T (J) \leq J$, then for some $J^\infty \in A_+(S)$, we have
$$ T^k(J) \downarrow J^\infty \quad \text{and} \quad T(J^\infty) \leq J^\infty.$$
\end{lem}
\begin{proof}
By the monotonicity of $T$, $T^k(J) \downarrow J^\infty$; by \cite[Lemma 7.30(2)]{bs}, $J^\infty \in A_+(S)$.
For every $k$, since $J^\infty \leq T^k(J)$, we have, by the monotonicity of $T$, $T(J^\infty) \leq T^{k+1}(J)$.
Hence $T(J^\infty) \leq J^\infty$.
\end{proof}

\mysmallskip
\begin{lem} \label{lma-P2}
Let $c > 1$. We have $T(c J^*) \leq c J^*$ and for some $J^\infty \in A_+(S)$,
$$ T^k(c J^*) \downarrow J^\infty, \qquad T ( J^\infty) = J^\infty, \qquad J^* \leq J^\infty   \leq c J^*. $$
\end{lem}

\begin{proof}
Since $c > 0$ and $J^* \in A_+(S)$, $c J^* \in A_+(S)$. 
Since $c > 1$ and the one-stage costs are nonnegative, it follows from the definition of $T$ and the fact $T(J^*) = J^*$ that $T (c J^*) \leq c J^*$.
Let $J^k = T^k(c J^*)$.  
By Lemma~\ref{lma-P1},
$$c J^* \geq J^k \downarrow J^\infty \geq J^* \quad \text{and} \quad T (J^\infty) \leq J^\infty,$$ 
where the inequality $J^\infty \geq J^*$ follows from the monotonicity of $T$ and the fact $T(J^*) = J^*$.
By rearranging the terms and using also the monotonicity of $T$, we have
$$ c J^* \geq J^\infty \geq T(J^\infty)  \geq J^*.$$

To prove $T(J^\infty) = J^\infty$, we now show $T (J^\infty) \geq J^\infty$, using the monotone convergence theorem. 
Consider an arbitrary $x \in S$. If $J^*(x) = \infty$, then $T(J^\infty)(x) = J^\infty(x) = \infty$ by the preceding relation. 
Suppose $J^*(x) < \infty$; we prove $T(J^\infty)(x) \geq J^\infty(x)$ below.

To simplify notation, consider the function 
$$ H(x, u, J) = g(x, u) + \int_S J(x') \, q(dx' \!\mid x, u), \qquad u \in U(x).$$
Recall that $T (J^*)(x) = \inf_{u \in U(x)} H(x, u, J^*)$ (cf.\ Eq.~(\ref{eq-T})).
Since $T (J^*)(x) = J^*(x) < \infty$, we have
$$ D(x) : = \big\{ u \in U(x) \mid H(x, u, J^*) < \infty \big\} \not= \emptyset.$$
For $u \in D(x)$, 
$$ H(x, u, c J^*) \leq c H(x, u, J^*) < \infty$$ 
(because $c > 1$ and $g \geq 0$), so in view of the relation $c J^* \geq J^k \downarrow J^\infty$, we have by the monotone convergence theorem \cite[p.\ 131]{Dud02}, 
\begin{equation} \label{prf-P0}
 H(x, u, J^\infty) = \lim_{k \to \infty} H(x, u, J^k).
\end{equation} 
Consequently,
\begin{equation} \label{prf-P1}
 H(x, u, J^\infty) \geq \limsup_{k \to \infty} \left\{ \inf_{u \in U(x)} H(x, u, J^k) \right\} = \lim_{k \to \infty} T (J^k)(x) = J^\infty(x), \qquad \forall \, u \in D(x).
\end{equation}
For $u \in U(x) \setminus D(x)$, 
$$ H(x, u, J^\infty) \geq H(x, u, J^*) = \infty.$$
Combining this with Eq.~(\ref{prf-P1}), we have
$$ T (J^\infty)(x) = \inf_{u \in U(x)} H(x, u, J^\infty) = \inf_{u \in D(x)} H(x, u, J^\infty) \geq J^\infty(x).$$
This completes the proof.
\end{proof} 
%\mysmallskip

We are now ready to prove Theorem~\ref{thm-valite-P}. We will give two different proofs, one here, and the other in Appendix~\ref{appsec-altprfP}. In the proof that follows, we will use a simple concavity property of $T$, which can be verified directly. (Hartley \cite{Har80} also used it in an alternative proof of Whittle's bridging condition.)
On the convex set $A_+(S)$, $T$ has the property that 
for any $\beta \in [0,1]$ and $J_1, J_2 \in A_+(S)$, 
\begin{equation} \label{eq-concT}
     T \big(\beta J_1 + (1 - \beta) J_2 \big) \geq \beta \, T (J_1) + ( 1 - \beta) \, T( J_2).
\end{equation}

We will also use Maitra and Sudderth's results \cite{MS92}.
Let $\omega_1$ be the first uncountable ordinal. 
For ordinals $\xi < \omega_1$,  define functions $J^\xi \in A_+(S)$ by transfinite recursion as follows. 
Let 
$$J^0 = T(\0), \qquad \quad J^\xi = T \Big( \sup_{\eta < \xi} J^\eta \Big), \quad \text{for} \ \xi > 0.$$ 
Also let 
$$ J^{\omega_1} = \sup_{\xi < \omega_1} J^\xi.$$
That all these functions are indeed in $A_+(S)$ is proved in \cite{MS92}. 
Moreover, Maitra and Sudderth \cite[Thm.\ 1.1]{MS92} proved  that 
\begin{equation} \label{eq-MS}
  T(J^{\omega_1}) = J^{\omega_1} = J^*.
\end{equation}  
(For ordinals, transfinite induction and transfinite recursion, see e.g., \cite[p.\ 27-28]{Kur-topology}, \cite[Secs.\ 1.3, A.3]{Dud02} or \cite[Chap.\ 1]{Sriv-borel}.)

%\mysmallskip
\begin{proof}[Proof of Theorem~\ref{thm-valite-P}]
Denote $J^\infty =  \lim_{k \to \infty} T^k(c J^*)$.  By Lemma~\ref{lma-P2},
$$T (J^\infty) = J^\infty, \qquad J^* \leq J^\infty \leq cJ^*.$$
We now prove $J^\infty = J^*$. This will prove part (a) and hence the entire theorem, as discussed earlier.

Let $\beta = 1/c < 1$. Since $J^\infty \leq c J^*$, $J^* \geq \beta J^\infty$, so by the monotonicity and concavity properties of $T$ (Eq.~(\ref{eq-concT})),
$$ T ( J^*) \geq T\big(\beta  J^\infty + (1-\beta) \, \0\big) \geq \beta \, T(J^\infty) + (1-\beta)\,T(\0).$$
Since $T(J^\infty) = J^\infty$ and $T(J^*) = J^*$, using the definition $J^0 = T(\0)$, we can write the above inequality equivalently as
\begin{equation} 
  J^* \geq \beta J^\infty + ( 1 - \beta)  J^0. \notag
\end{equation}  

We now apply transfinite induction. Suppose that for an ordinal $\xi \leq \omega_1$, 
$$ J^* \geq \beta J^\infty + ( 1 - \beta) J^\eta, \qquad \forall \, \eta < \xi.$$
Then
\begin{align*}
   J^*  \geq \sup_{\eta < \xi} \Big\{ \beta J^\infty + ( 1 - \beta)  J^\eta  \Big\} 
         =  \beta J^\infty + ( 1 - \beta) \,  \sup_{\eta < \xi} J^\eta.
\end{align*}
Consequently, by the monotonicity and concavity properties of $T$, 
\begin{align} 
 J^* & \geq \beta J^\infty + ( 1 - \beta) \,  T \Big(\sup_{\eta < \xi} J^\eta \Big) 
     = \beta J^\infty + ( 1 - \beta)   J^\xi,  \label{prf-P2}
\end{align} 
where in the first inequality we also used the fact $T(J^\infty) = J^\infty$ and $T(J^*) = J^*$, 
and in the equality we used the definition of $J^\xi$ for $\xi < \omega_1$, and the definition of $J^{\omega_1}$ for $\xi = \omega_1$, together with the fact $T(J^{\omega_1}) = J^{\omega_1}$ (\cite[Thm.\ 1.1]{MS92}; cf.\ Eq.~(\ref{eq-MS})).
This proves, by transfinite induction, that the inequality (\ref{prf-P2}) holds for all $\xi \leq \omega_1$, and in particular,
\begin{equation} \label{prf-P3}
   J^* \geq \beta J^\infty + ( 1 - \beta)  J^{\omega_1}.  
\end{equation}
Since $J^{\omega_1} = J^*$ by \cite[Thm.\ 1.1]{MS92} (cf.\ Eq.~(\ref{eq-MS})) and $J^* \geq 0$,
we have $ J^* \geq J^\infty$ by Eq.~(\ref{prf-P3}).
We also have $J^\infty \geq J^*$. Therefore, $J^\infty = J^*.$
\end{proof}
%\mysmallskip
 
We mention two immediate implications of Theorem~\ref{thm-valite-P}. 

%\mysmallskip
\begin{cor} \label{cor-valite-P2}
{\rm (P)} \ Suppose that the state space $S$ is finite and $J^*$ is real-valued. 
Let $\mathcal{V}$ be the set of nonnegative, real-valued functions $J$ such that $J(x) = 0$ for all $x\in S$ with $J^*(x)=0$.
Then $J^*$ is the unique fixed point of $T$ within $\mathcal{V}$. Moreover, $T^k(J) \to J^*$ for all $J \in \mathcal{V}$.
\end{cor}

\begin{proof} 
By Theorem~\ref{thm-valite-P}(c), $J^*$ is the unique fixed point of $T$ in $\{ J \geq 0 \mid \exists \, c \in \Re_+ \ s.t.\  J \leq c J^* \} = \mathcal{V}$.
By Theorem~\ref{thm-valite-P}(b), we have $T^k(J) \to J^*$ for all $J \in \mathcal{V}$ if $T^k(\0) \uparrow J^*$. 
The latter holds when $S$ is finite and $J^*$ is finite everywhere (cf.\ \cite[Thm.\ 7.3.10(a)]{Puterman94}). The reason is that $T^k(\0)$ converges to its limit $J_\infty$ uniformly, i.e., for any $\epsilon > 0$, $\| T^k(\0) - J_\infty \|_\infty \leq \epsilon$ for all $k$ sufficiently large. Thus with $\1$ denoting the constant function $1$, we have, by the monotonicity of $T$, that for all $k$ sufficiently large, 
$$ T(J_\infty) \leq T\big(T^k(\0) + \epsilon \1\big) \leq T^{k+1}(\0) + \epsilon \1. $$
Since $\epsilon$ is arbitrary, this implies $T(J_\infty) \leq J_\infty$. Since $J_\infty \leq T(J_\infty)$ \cite[Prop.\ 9.16]{bs}, we have 
$J_\infty = T(J_\infty)$ and hence by \cite[Prop.\ 9.16]{bs}, $J_\infty = J^*$, i.e., $T^k(\0) \uparrow J^*$. 
\end{proof}
%\mysmallskip

In connection with Cor.~\ref{cor-valite-P2}, we note that even when $S$ is finite, $J_\infty \not= J^*$ is possible if $J^*$ is not real-valued. As an example, let $S=\{0,1,2\}$, $C=(0,1) \cup \{t\}$ and $U(x) = C$ for all $x \in S$. States $0$, $1$ are absorbing; the self-transition costs are $g(0,u)=0$ and $g(1,u)=1$ for all $u$ (so $J^*(0)=0$, $J^*(1)=\infty$). At state $2$, for control $u \in (0,1)$, the one-stage cost is $g(2,u)=0$ and the next state is state $1$ with probability $u$ and state $0$ with probability $(1-u)$; and
for control $u=t$, $g(2,t)= 1$ and the next state is state $0$. Then $T^k(\0)(2) = 0$ for all $k$, so $J_\infty(2) = 0$; but $J^*(2) = 1$.

The next corollary relates to the convergence of the value functions generated by exact policy iteration (assuming it can be carried out). First, note that Theorem~\ref{thm-valite-P} implies also the convergence of value iteration for certain initial functions that cannot be bounded by a multiple of $J^*$. Specifically, by the theorem, if $J \in A_+(S)$ is such that $\underline{J} \leq T^n(J) \leq c J^*$ for some $n \geq 0$ and $c > 1$, where $\underline{J}$ is as in Theorem~\ref{thm-valite-P}(b), then $T^k (J) \to J^*$.

%\mysmallskip
\begin{cor}[Convergence result for policy iteration under (P)] \label{cor-valite-P1}
Suppose $\{\mu_k\}$  is a sequence of stationary policies that satisfy
$J_{\mu_k} \in A_+(S)$, $T_{\mu_{k+1}} (J_{\mu_k}) = T(J_{\mu_k})$ for $k \geq 0$. Then if $\mu_0$ is such that $T^n(J_{\mu_0}) \leq c J^*$ for some $n \geq 0$ and $c > 1$, we have $J_{\mu_k} \to J^*$.
\end{cor} 
\begin{proof} By Theorem~\ref{thm-valite-P}, the condition on $\mu_0$ implies that $T^k(J_{\mu_0}) \to J^*$. Since $T_{\mu_{k+1}} (J_{\mu_k}) = T(J_{\mu_k}) \leq T_{\mu_k}(J_{\mu_k})=J_{\mu_k}$,  by the monotonicity of $T_{\mu_{k+1}}$, we have that $T_{\mu_{k+1}}^n (J_{\mu_k})$ is nonincreasing as $n \to \infty$, so $J_{\mu_{k+1}} = \lim_{n \to \infty} T_{\mu_{k+1}}^n(\0) \leq  \lim_{n \to \infty} T_{\mu_{k+1}}^{n}\big(T(J_{\mu_k})\big) \leq T(J_{\mu_k})$. It then follows that $J^* \leq J_{\mu_{k+1}} \leq T^{k+1}(J_{\mu_0})$, and hence $J_{\mu_k} \to J^*$.
\end{proof}
%\mysmallskip

A similar conclusion holds for modified policy iteration (assuming it can be carried out).%footnote starts
\footnote{Specifically, suppose $\big\{(J_k, \mu_k)\big\}$ is a sequence of cost function and stationary policy pairs that satisfy
$J_{k} \in A_+(S)$, $J_{k+1} = T^{n_k}_{\mu_k}(J_k)$ for some $n_k \geq 1$, $T_{\mu_{k+1}} (J_{k+1}) = T(J_{k+1})$ for  all $k \geq 0$. Suppose in addition that $T_{\mu_0}(J_0) \leq J_0$ and $T^n(J_0) \leq c J^*$ for some $n \geq 0$ and $c > 1$. Then $J^* \leq J_{\mu_k} \leq J_k \to J^*$.}
%footnote ends
Note also that Cor.~\ref{cor-valite-P1} does not imply $J_\infty = J^*$. Neither does it imply the existence of an $\epsilon$-optimal stationary policy because the convergence of $J_{\mu_k}$ to $J^*$ is only pointwise, but it does imply the existence of an $\epsilon$-optimal stationary policy for each state.

Finally, we remark that 
when the problem does not admit a near-optimal stationary policy for some state, one cannot hope to find an initial function $J$ with the desired property $J \leq c J^*$ for some $c > 1$, among the cost functions of stationary policies. This is shown in the following proposition. 

%\mysmallskip
\begin{prop} \label{prp-valite-P3}
{\rm (P)} \ Suppose that for some $\bar x \in S$ and $\epsilon > 0$, an $\epsilon$-optimal stationary policy for $\bar x$ does not exist. Then there exists no stationary policy $\mu$ such that $J_\mu \in A_+(S)$, $T^k(J_\mu) \to J^*$ (hence there exists no $J_\mu \in A_+(S)$ with $T^n(J_\mu) \leq c J^*$ for some $n \geq 0$ and $c > 1$).
\end{prop}

\begin{proof}
To arrive at a contradiction, suppose $\mu$ is a stationary policy with $J_\mu \in A_+(S)$, $T^k(J_\mu) \to J^*$. 
Let $k$ be large enough so that $T^k(J_\mu)(\bar x) \leq J^*(\bar x) + \epsilon/2$.
By the selection theorem of \cite[Prop.\ 7.50]{bs}, there exist stationary policies $\mu_1, \mu_2, \ldots, \mu_k$ satisfying
$$ T_{\mu_i} \big(T^{i-1} (J_\mu) \big) \leq T^{i} (J_\mu) + \epsilon/(2k), \qquad i = 1, 2, \ldots, k.$$
Consider the Markov policy $\pi =(\mu_k, \mu_{k-1}, \ldots, \mu_1, \mu, \mu, \ldots)$. By a direct calculation we have
$$ J_\pi(\bar x)  = \big(T_{\mu_k} \circ T_{\mu_{k-1}} \circ \cdots T_{\mu_1} \big) (J_\mu) \leq T^k(J_\mu)(\bar x) + \epsilon/2 \leq   J^*(\bar x) + \epsilon.$$
Now consider an associated subproblem where at every state $x$, there are only $k+1$ ``controls'' corresponding to $\{\mu, \mu_1, \mu_2, \ldots, \mu_k\}$. Since the number of controls is finite, one can verify by a direct calculation that value iteration starting with the constant function zero does not have measurability issues and maintains the cost function iterates within the family of nonnegative universally measurable functions. Using a theorem for (P) \cite[Props.\ 5.10, 5.4]{bs}, we then obtain that for this subproblem, the optimal cost function is universally measurable and there exists 
an optimal nonrandomized stationary policy $\tilde \mu$. Clearly, $\tilde \mu$ corresponds to a universally measurable stationary policy in the original problem.
By the optimality of $\tilde \mu$ in the subproblem, we have $J_{\tilde \mu}(\bar x) \leq J_\pi(\bar x) \leq J^*(\bar x) + \epsilon,$ which contradicts the assumption that there exists no $\epsilon$-optimal stationary policy for state $\bar x$. This proves the main part of the proposition; the rest then follows by Theorem~\ref{thm-valite-P}(b).
\end{proof}
%\mysmallskip

We illustrate Prop.~\ref{prp-valite-P3} by an example based on \cite[Example 6.1]{Str-negative}, in which the cost function of every stationary policy, although not nearly optimal, is a fixed point of $T$. Let $S =\{0,1,2\}$, $C= (0,1)$ and $U(x) = C$ for all $x$. State $0$ is cost-free and absorbing. From state $2$, any control leads to state $1$ with cost $1$. For state $1$, under control $u$, we have probability $u$ to transition to state $0$ with transition cost $1$, and probability $(1-u)$ to transition to state $1$ with self-transition cost $0$. 
The optimal costs are $J^*(0)= J^*(1) = 0, J^*(2) = 1$. 
An $\epsilon$-optimal Markov policy, for example, is to apply at state $1$ control $u_k$ for the $k$th stage, with $\sum_k u_k \leq \epsilon$. 
No stationary policy is $\epsilon$-optimal for states $1$ and $2$: for any stationary policy $\mu$, transition from state $1$ to state $0$ occurs with probability one, so $J_\mu(0) = 0$, $J_\mu(1) = 1, J_\mu(2) = 2$, and moreover, $J_\mu$ is also a fixed point of $T$. 

\subsection{Convergence Properties of Mixed Value and Policy Iteration} \label{secP5.2}

We now consider the mixed value and policy iteration method in case (P).
Unlike case (N) where it is natural to apply the method with the initial iterate $J_0 = \0, Q_0 = \0$, 
here, as can be shown by a direct calculation, doing so reduces the method to value iteration $T^k(\0)$, and this is undesirable computationally even if $T^k(\0) \to J^*$, which is not guaranteed to hold. Our interest thus lies primarily in applying the method with an initial $(J_0, Q_0)$ above the optimal costs. 
We will apply Theorem~\ref{thm-valite-P} to analyze the convergence of the basic algorithm (\ref{alg-1a})-(\ref{alg-1b}) given in Section~\ref{sec3.2}.
We will also discuss another variant algorithm that connects to linear programming, and prove its convergence.

\smallskip
\begin{thm} \label{thm-P-mixed}
{\rm (P)} Let $J_0 \in A_+(S)$ and $Q_0 \in A_+(\Gamma)$. \\*[0.01in]
(a)  Suppose
\begin{equation} \label{eq-thm-inicond}
  J^* \leq J_0 \leq c J^* \ \ \text{for some} \ c > 1, \quad \text{and} \quad Q_0 \geq Q^*.
\end{equation}  
Then
the sequence $\big\{(J_k, Q_k)\big\}$ generated by the iteration (\ref{alg-1a})-(\ref{alg-1b}) 
converges to $(J^*, Q^*)$.\\*[0.07in]
(b) If $T^k(\0) \uparrow J^*$, then the initial condition (\ref{eq-thm-inicond}) in (a) on $(J_0, Q_0)$ can be relaxed to $J_0 \leq c J^*$.\\*[0.07in]
(c) Suppose $T^k\big(\underline{J}\big) \uparrow J^*$ for some $\underline{J} \in A_+(S)$. 
Then the conclusion of (a) holds for the iteration (\ref{alg-1a})-(\ref{alg-1b}) that always defines $Q_{k+1}$ using the first rule in (\ref{alg-1a}), under the initial condition that 
$$ \underline{J} \leq  J_0 \leq c J^* \ \ \text{for some} \ c > 1, \quad \text{and} \quad Q_0(x, u) \geq \underline{J}(x) \quad \forall \, (x, u) \in \Gamma.$$
\end{thm}
\mysmallskip 
 
 In either part of the theorem, it is assumed that $J_0 \leq c J^*$ for some $c > 1$.
 We prove first that under this condition on $J_0$, the limits of the iterates $(J_k,Q_k)$ can be upper bounded by $(J^*, Q^*)$.

\mysmallskip
\begin{lem} \label{lma-prfP-m1}
{\rm (P)} \ Let $J_0 \in A_+(S)$ and $Q_0 \in A_+(\Gamma)$. If  $J_0 \leq c J^*$ for some $c > 1$, then the sequence $\big\{(J_k, Q_k)\big\}$ generated by the iteration (\ref{alg-1a})-(\ref{alg-1b}) satisfies that
$$\limsup_{k \to \infty} J_k \leq J^*, \qquad 
 \limsup_{k \to \infty} Q_k \leq Q^*.$$
\end{lem}

 \begin{proof}
Let $J^k = T^k(c J^*)$. Since $J_0 \leq cJ^*$, we have $J_{k} \leq T^{k}(J_0) \leq J^k$ for every $k$,
by Lemma~\ref{lma-bound-NP} and the monotonicity of $T$. 
Since $J^k  \downarrow J^*$ by Theorem~\ref{thm-valite-P}(a), 
$\limsup_{k \to \infty} J_k \leq J^*$.

Consider now $Q_k(x,u)$ for each $(x, u) \in \Gamma$, and note that $Q^*(x,u) = H(x, u, J^*)$ by definition (cf.\ Eqs.~(\ref{eq-H}), (\ref{eq-QH})). 
By Lemma~\ref{lma-bound-NP}, for every $k \geq 0$, 
$$ Q_{k+1}(x, u) \leq H(x,u, J_k) \leq H(x, u, J^k).$$
If $Q^*(x, u) < \infty$, then we have
$$ \lim_{k \to \infty} H(x, u, J^k) = H \big(x, u, \lim_{k \to \infty} J^k \big) = H(x, u, J^*) = Q^*(x,u),$$
where the first equality follows from the monotone convergence theorem as we showed with Eq.~(\ref{prf-P0}) in the proof of Lemma~\ref{lma-P2}. 
By combining the preceding two relations, we obtain 
$$\limsup_{k \to \infty} Q_{k+1}(x, u) \leq Q^*(x,u).$$ 
This inequality also holds, trivially, if $Q^*(x,u) = \infty$. Therefore, $\limsup_{k \to \infty} Q_k \leq Q^*$. 
\end{proof}
\mysmallskip

We now proceed to prove the theorem by bounding the iterates from below.\footnote{For part (a), we will use the lower bounds given in Lemma~\ref{lma-bound-NP}, which rely on the relation $Q_{\k,J^*} = Q^*$ for all $\k \in \K$ (cf.\ Prop.~\ref{prp-propertyF}(d)). This relation will be proved as Prop.~\ref{lma-P-optstop} in Appendix~\ref{appsec2}, and it is needed in the analysis for the algorithm that can set $Q_{k+1}$ to be $Q_{\k_k,J_k}$ at some iterations.}

\begin{proof}[Proof of Theorem~\ref{thm-P-mixed}]
\noindent (a) Since $J_0 \geq J^*$ and $Q_0 \geq Q^*$,  we have $J_k \geq J^*, Q_k \geq Q^*$ by Lemma~\ref{lma-bound-NP}, and hence
$J_k \to J^*$, $Q_k \to Q^*$ by Lemma~\ref{lma-prfP-m1}.

\vspace*{0.05in}
\noindent (b) Starting with $J_0 \geq 0, Q_0 \geq 0$, let us prove by induction that for every $k \geq 0$,
\begin{equation} \label{eq-prfP-m1b}
   J_k \geq T^k(\0), \qquad  Q_k(x,u) \geq T^k(\0)(x), \quad \forall \, (x,u) \in \Gamma.
\end{equation}   
By Lemma~\ref{lma-prfP-m1} and the assumption $T^k(\0) \uparrow J^*$, the first inequality above will immediately imply that $J_k \to J^*$. 

To simplify notation, let
$\hat J_k=T^k(\0)$ and define ${\hat J}^e_k \in A_+(\Gamma)$ by
$$ {\hat J}^e_k(x,u) = {\hat J}_k(x), \quad \forall \, (x,u) \in \Gamma.$$ 
We will use the following facts.
For any $\k \in \K$, in view of $g \geq 0$ and the fact $\hat J_k \geq \hat J_{k-1} \geq \cdots \geq 0$,  
a direct calculation using the definition (\ref{def-F}) of $F_\k$ and its monotonicity shows that
\begin{equation} \label{eq-prfP-m2}
 F_\k\big(\0\,; {\hat J}_k\big) \geq {\hat J}^e_1, \quad F_\k\big({\hat J}^e_1 \,;  {\hat J}_k\big) \geq {\hat J}^e_2,  \  \ \cdots \ \  F_\k\big({\hat J}^e_{k-1} \,;  {\hat J}_k\big) \geq {\hat J}^e_k,
\end{equation} 
and that for every $n \geq 1$, 
\begin{equation} \label{eq-prfP-m3}
    F_\k^n\big({\hat J}^e_{k} \,;  {\hat J}_k\big) \geq    F_\k\big({\hat J}^e_k\,; {\hat J}_k \big). 
\end{equation}
In view of $g \geq 0$ and the definition of $H(x, u, J)$ (cf.\ Eq.~(\ref{eq-H})), a direct calculation shows that 
\begin{equation} \label{eq-prfP-m4}
 F_\k\big({\hat J}^e_k\,; {\hat J}_k \big)(x, u) = H\big(x, u, {\hat J}_k\big) \geq T \big(\hat J_k\big)(x) = \hat J_{k+1}(x), \qquad  \forall \, (x, u) \in \Gamma.
\end{equation} 

Now suppose Eq.~(\ref{eq-prfP-m1b}) holds for some $k \geq 0$. Consider the $k$th iteration of the algorithm.
We have either $Q_{k+1} = F^n_\k(Q_k\,; J_k )$ or $Q_{k+1} = Q_{\k, J_k}$ for some $\k \in \K$ and $n\geq 1$.
For the case $Q_{k+1} = F^n_\k(Q_k\,; J_k )$, we have 
$$  F^n_\k(Q_k\,; J_k ) \geq F^n_\k\big({\hat J}^e_k \,; {\hat J}_k\big) \geq  F_\k \big({\hat J}^e_k \,;  \hat J_k\big),$$
where the first inequality follows from the monotonicity of $F_\k$ (cf.\ Eq.~(\ref{eq-Fmono})) and the 
induction hypothesis that 
$J_k \geq \hat J_k, Q_k \geq {\hat J}_k^e$, and the second inequality follows from Eq.~(\ref{eq-prfP-m3}).
For the case $Q_{k+1} = Q_{\k, J_k}$, we have 
$$ Q_{\k, J_k} \geq F^{k+1}_\k \big(\0\,; J_k\big) \geq F^{k+1}_\k \big(\0\,; \hat J_k \big) \geq F_\k \big({\hat J}^e_k \,; \hat J_k\big),$$
where the first inequality holds because $F^n_\k (\0\,; J_k) \uparrow Q_{\k, J_k}$ as $n \to \infty$ (Prop.\ \ref{prp-optstop-cost}),
the second inequality follows from the induction hypothesis $J_k \geq \hat J_k$ and the monotonicity of $F_\k$ (cf.\ Eq.~(\ref{eq-Fmono})), and the third inequality follows from Eq.~(\ref{eq-prfP-m2}) and the monotonicity of $F_\k(\cdot\,; \hat J_k)$.
Thus in either case, we have 
$$ Q_{k+1} \geq F_\k \big({\hat J}^e_k \,; \hat J_k\big) \geq {\hat J}_{k+1}^e, \qquad J_{k+1} = M(Q_{k+1}) \geq {\hat J}_{k+1},$$
where Eq.~(\ref{eq-prfP-m4}) is used in the second inequality of the first relation above. 
This completes the induction and establishes Eq.~(\ref{eq-prfP-m1b}) for all $k$. 

We can now conclude that $J_k \to J^*$, as discussed earlier.
We prove $Q_k \to Q^*$ next.
As we just proved, $Q_{k+1} \geq F_\k \big({\hat J}^e_k \,; \hat J_k\big)$ for every $k$. By Eq.~(\ref{eq-prfP-m4}), this is equivalent to
\begin{equation} \label{eq-prfP-m5}
   Q_{k+1}(x, u) \geq H\big(x, u, {\hat J}_k\big), \qquad \forall \, (x, u) \in \Gamma.
\end{equation}   
Since $\hat J_k \uparrow J^*$ and $\hat J_k \geq 0$, by 
the monotone convergence theorem, 
$$H\big(x, u, \hat J_k\big) \uparrow H\big(x, u,  J^*\big) = Q^*(x, u)$$ (cf.\ Eqs.~(\ref{eq-H}), (\ref{eq-QH})).
Together with Lemma~\ref{lma-prfP-m1}, the preceding two relations imply that $Q_k \to Q^*$.

\vspace*{0.05in}
\noindent (c) 
By assumption $T^k(\underline{J}) \uparrow J^*$.
For the algorithm stated in (c), if we define $\hat J_0 = \underline{J}$, $\hat J_k = T^k(\underline{J})$ for $k \geq 1$, 
then the same arguments in the preceding proof for part (b) go through to establish that Eqs.~(\ref{eq-prfP-m3})-(\ref{eq-prfP-m4}) hold, that for every $k$,
$$   J_k \geq T^k(\underline{J}), \qquad  Q_k(x,u) \geq T^k(\underline{J})(x), \quad \forall \, (x,u) \in \Gamma,$$
and that $J_k \to J^*, Q_k \to Q^*$. The reason that the proof arguments for part (b) are applicable here is as follows. 
Among the crucial facts used in the proof of part (b), the only one that does not hold under the present initial condition on $J_0$ is the first inequality $F_\k\big(\0\,; {\hat J}_k\big) \geq {\hat J}^e_1$ in Eq.~(\ref{eq-prfP-m2}). This relation is needed in the convergence proof only when
$Q_{k+1}$ is generated by the second rule of (\ref{alg-1a}) as $Q_{k+1} = Q_{\k_k, J_k}$; but such cases are ruled out by the assumptions of part~(c).
 \end{proof}
 
\subsubsection*{A Variation of the Basic Algorithm (\ref{alg-1a})-(\ref{alg-1b})}

Let us consider a variation of the algorithm (\ref{alg-1a})-(\ref{alg-1b}), whereby instead of (\ref{alg-1a}), we use a different rule to update $Q_{k+1}$:\moveup\moveup 
\begin{itemize}
\item Choose $\k_k =(\mu_k,B_k) \in \K$, and find $Q_{k+1} \in A_+(\Gamma)$ such that
\begin{equation}
   Q_{k+1} \leq F_{\k_k}(Q_{k+1} ; J_k), \qquad Q_{k+1} \geq Q_{\k_k,J_k}. \label{alg-1av}
\end{equation}   
Then let 
 \begin{equation}
 J_{k+1} = M (Q_{k+1}). \label{alg-1bv}
 \end{equation}
 \end{itemize}

This algorithm is motivated by a computational issue in case (P).
Unlike (D)(N), control problems of type (P), even when the spaces $S, C$ are discrete, do not admit a linear programming formulation in general (cf.\ \cite[Prop.\ 9.10(P)]{bs}, \cite[Sec.\ 7.3.6]{Puterman94}). Thus to calculate $Q_{\k_k,J_k}$ in the algorithm (\ref{alg-1a})-(\ref{alg-1b}) without iterating $F_{\k_k}^n(\0\,; J_k)$ up to convergence, we cannot solve the optimal stopping problem associated with $(\k_k, J_k)$ by simply solving some linear program. 

On the other hand, an upper bound on $Q_{\k_k,J_k}$ will suffice if it also satisfies the first relation in (\ref{alg-1av}), as we show in the theorem below.
Unlike computing $Q_{\k_k,J_k}$, a solution to (\ref{alg-1av}) can be computed by solving a linear program associated with the optimal stopping problem defined by $(\k_k, J_k)$, under certain conditions that involve $(\k_k, J_k)$, as we will show in Lemma~\ref{lma-lp}, Appendix~\ref{appsec-optstop-lp}. 
 These conditions are satisfied, for example, if $S$ and $C$ are countable and $J_k$ is finite on $B_k$; see Remark~\ref{rmk-lp} in Appendix~\ref{appsec-optstop-lp}. Given that if $J_0 \leq c J^*$ for some $c > 1$, the algorithm (\ref{alg-1av})-(\ref{alg-1bv}) will generate $J_k$ with $J_k \leq c J^*$ throughout (see the theorem below), this means that the step (\ref{alg-1av}) can be carried out by linear programming for countable-spaces problems where $J^*$ is finite everywhere, in particular.

\mysmallskip
\begin{thm}  \label{thm-P-mixed2}
{\rm (P)} \ Under the same conditions as in Theorem~\ref{thm-P-mixed}(a) or (b), the sequence $\big\{(J_k,Q_k)\big\}$ generated by the iteration~(\ref{alg-1av})-(\ref{alg-1bv}) satisfies $J_k \leq c J^*$ for all $k$, and converges to $(J^*,Q^*)$.
\end{thm}

\begin{proof}
The proof is similar to that for Theorem~\ref{thm-P-mixed}(a)-(b). We will bound $(J_k,Q_k)$ from above and from below. 
By Lemma~\ref{lma-simplebound-DNP},
\begin{equation} 
 F_{\k_k}(Q_{k+1} ; J_k)(x,u) \leq H(x, u, J_k), \quad  \forall \, (x, u) \in \Gamma, \qquad   M \big( F_{\k_k}(Q_{k+1} ; J_k)\big) \leq T(J_k).  \notag 
 \end{equation}
From this and the upper bound on $Q_{k+1}$ given in Eq.~(\ref{alg-1av}), we have 
$$    Q_{k+1}(x, u) \leq H(x, u, J_{k}), \quad  \forall \, (x, u) \in \Gamma,\qquad J_{k+1}  \leq  T(J_k). $$
The second inequality above implies, by the monotonicity of $T$, that $J_k \leq T^k(J_0)$ for every $k$. Since $J_0 \leq c J^*$ by assumption, we obtain $J_k \leq T^k(c J^*) \leq c J^*$.

In the proof of Lemma~\ref{lma-prfP-m1}, use the preceding upper bounds on $J_k$ and $Q_{k+1}$ in place of Lemma~\ref{lma-bound-NP} (which states the same bounds for the basic algorithm). We obtain that the conclusion of Lemma~\ref{lma-prfP-m1} holds for the iteration~(\ref{alg-1av})-(\ref{alg-1bv}):
\begin{equation} \label{eq-prfP-m1-rep}
  \limsup_{k \to \infty} J_k \leq J^*, \qquad \limsup_{k \to \infty} Q_k \leq Q^*.
\end{equation}

Under the conditions of Theorem~\ref{thm-P-mixed}(a), we have 
$J_0 \geq J^*, Q_0 \geq Q^*$. Lemma~\ref{lma-bound-NP} shows that if $Q_{k+1} = Q_{\k_k,J_k}$ at every iteration of the algorithm, then $J_k \geq J^*, Q_k \geq Q^*$ for all $k$.
Since here we have $Q_{k+1} \geq Q_{\k_k,J_k}$ by Eq.~(\ref{alg-1av}), the same induction-based proof for Lemma~\ref{lma-bound-NP} shows that for the iteration~(\ref{alg-1av})-(\ref{alg-1bv}), we have $J_k \geq J^*, Q_k \geq Q^*$ for all $k$ as well. This together with Eq.~(\ref{eq-prfP-m1-rep}) implies that $J_k \to J^*, Q_k \to Q^*$.

Similarly, under the conditions of Theorem~\ref{thm-P-mixed}(b), the proof of Theorem~\ref{thm-P-mixed}(b) established the lower bounds (\ref{eq-prfP-m1b}), (\ref{eq-prfP-m5}) on $J_k, Q_k$ for the case where $Q_{k+1} = Q_{\k_k, J_k}$ at every iteration, and these lower bounds also hold for the iteration~(\ref{alg-1av})-(\ref{alg-1bv}) since $Q_{k+1} \geq Q_{\k_k, J_k}$. Together with Eq.~(\ref{eq-prfP-m1-rep}), they imply that $J_k \to J^*, Q_k \to Q^*$, as the proof of Theorem~\ref{thm-P-mixed}(b) showed.
\end{proof}
 \mysmallskip
 
\begin{rem}[Finding near-optimal policies]
Under (P), from the sequence $\{J_k\}$ generated by the algorithm (\ref{alg-1a})-(\ref{alg-1b}) or (\ref{alg-1av})-(\ref{alg-1bv}), in general one cannot extract easily near-optimal policies in the manner of Remark~\ref{rem-policyD}. Even if $J^*$ was available, an $\epsilon$-optimal stationary policy may not exist (see the discussion after Prop.~\ref{prp-valite-P3} or \cite[p.\ 145]{Ber-DP13} for an example). 
Nevertheless, let us describe first a favorable case where one can find stationary policies $\nu_k$ with $J_{\nu_k} \to J^*$.
In particular, let $\{J_k\} \subset A_+(S)$ be such that $J_k\to J^*$. 
Suppose $J_k$ and $\nu_k$ satisfy
\begin{equation} \label{eq-rem-policyP}
   T_{\nu_k}(J_k)=T(J_k)\leq J_k, \qquad \forall \, k \geq 1.
\end{equation}   
Then $J_{\nu_k} \to J^*.$ (To see this, note that by \cite[Prop.\ 9.11]{bs}, $J_{\nu_k}$ is the ``smallest" nonnegative function $J \in \M(S)$ satisfying $T_{\nu_k} (J) \leq J$, so the assumption implies that $J_{\nu_k}\le J_k$. Since $J_k\to J^*$, the result follows.) The assumption (\ref{eq-rem-policyP}), however, need not always hold for our algorithms. 

Consider now a more general case. 
Suppose $J_k \to J^*$ and for some $c > 1$, $J^* \leq J_k \leq c J^*$ for all $k \geq 1$. These conditions hold for our algorithms under the initial condition given in Theorem~\ref{thm-P-mixed}(a).  With such $\{J_k\}$, one can construct asymptotically near-optimal policies, even though these policies are not necessarily stationary or Markov and the constructions can be complicated. To simply discuss just the idea of the construction, let us assume in addition that $T^k(c J^*)$ converges to $J^*$ uniformly, i.e., $\epsilon_k = \sup_{x \in S} \big( T^k(c J^*)(x) - J^*(x) \big) \to 0$. Then given $J_k$, for any $\delta > 0$, we can find $n \geq 1$ such that $T^n(J_k) \leq J_k + \delta/2$ (because $T^n(J_k) \leq T^n(c J^*) \leq J^* + \epsilon_n \leq J_k + \epsilon_n$ and $\epsilon_n \to 0$), 
together with an $n$-stage policy $(\mu_1, \mu_2, \ldots, \mu_n)$ such that
\begin{equation} \label{eq-rem-policyP2}
    \big(T_{\mu_1} \circ T_{\mu_2} \circ \cdots T_{\mu_n}\big) (J_k)  \leq J_k + \delta.
\end{equation}
Using this inequality, for any $\epsilon' > 0$, an $(\epsilon_k + \epsilon')$-optimal Markov policy can be constructed as follows (similar to how a near-optimal Markov policy is formed when $J^*$ is given \cite[p.\ 238]{bs}). 
Let $\delta_m, m \geq 1,$ be a sequence of positive numbers with $\sum_m \delta_m = \epsilon'$. For each $m$, let $\pi^{(m)}=\big(\mu_1^{(m)}, \mu_2^{(m)}, \ldots, \mu_{n_m}^{(m)}\big)$ be an $n_m$-stage policy satisfying the inequality~(\ref{eq-rem-policyP2}) for $\delta = \delta_m$. Then the Markov policy $\pi = (\pi^{(1)}, \pi^{(2)}, \ldots)$ is $(\epsilon_k + \epsilon')$-optimal because $J_\pi \leq J_k + \epsilon' \leq J^* + (\epsilon_k+\epsilon')$.%footnote starts
\footnote{To see this,
write $V_\pi(n, J)$ for the cost function of $\pi$ in the $n$-stage problem, with the terminal cost at stage $n+1$ given by $J$; in particular, $V_\pi(0, J) = J$.
Consider the $n$-stage cost function of $\pi$, denoted $J_{\pi,n}$, for stages $n=N_1, N_2, \ldots,$ where $N_i = \sum_{m=1}^i n_m$.
Using the inequality (\ref{eq-rem-policyP2}) for each segment $\pi^{(m)}$ of the policy $\pi$, we have
$$ J_{\pi,N_i} \leq V_\pi(N_i, J_k) \leq V_\pi(N_{i-1}, J_k) + \delta_i \leq \cdots \leq V_\pi(0, J_k) + \sum_{m=1}^i \delta_m. $$ 
Since $J_{\pi,N_i} \to J_\pi$ as the sequence $\{N_i\}$ of stages increases to infinity, it follows that $J_{\pi} \leq J_k + \sum_m \delta_m \leq J^* + \epsilon_k + \epsilon'$.}

In the above, uniform convergence of $T^k(c J^*)$ and $J_k$ to $J^*$ was assumed. If the convergence is only pointwise, then instead of a Markov policy, one can similarly construct a history-dependent policy $\pi$ with $J_\pi(x) \leq J_k(x) + \epsilon', x \in S$, for a given $\epsilon' > 0$ and $k$. It also seems possible to extract such policies directly from some of our algorithms, without computing the value iteration $T^n(J_k)$, based on the optimal stopping interpretation of $F_\k$ for the step (\ref{alg-1a}) (cf.\ Remarks~\ref{rmk-optstop1}, \ref{rmk-optstop2}). We do not discuss this possibility in this paper, however, due to its complexity and space limit.
\qed
\end{rem} 
%\mysmallskip

\section{Applications in Semicontinuous Models} \label{sec6}
\markboth{\rm \S \thesection. Applications in Semicontinuous Models}{\rm \S \thesection. Applications in Semicontinuous Models}

We discuss in this section some direct applications of our results for two special cases of the stochastic control model given in Section~\ref{sec2.1}: the upper semicontinuous model and the lower semicontinuous model. We take both models from~\cite[Chap.\ 8]{bs}. 
To keep notation simple, we will consider a simplified version of the lower semicontinuous model. However, our results are applicable to the model as defined in \cite[Chap.\ 8]{bs}.

To apply the mixed value and policy iteration method in these semicontinuous models, it is desirable to work with semicontinuous functions instead of lower semi-analytic functions. Our strategy to keep the function iterates within the set of semicontinuous functions is to choose properly the parameters of the mappings $F_\k$ involved in the method. For the lower semicontinuous model, we will use Lusin's theorem for this purpose. For the upper semicontinuous model, our result is weaker in the sense that we will need additional conditions in order to obtain nontrivial policy iteration-like algorithms that operate within the set of upper semicontinuous functions.

In this section we will also give a result about the structure of $J^*$ and optimal policies for the upper semicontinuous model in case (P), as an application of Theorem~\ref{thm-valite-P}.

We need some definitions. 
Let $X$ be a metrizable topological space. A function $f: X \to [-\infty, \infty]$ is said to be \emph{upper semicontinuous} (u.s.c.) if for every $c \in \Re$, its upper level set $\{ x \in X \mid f(x) \geq c\}$ is closed in $X$. Equivalently, $f$ is u.s.c.\ if and only if for any sequence $\{x_n\}$ in $X$ converging to some $x \in X$, we have $\limsup_{n \to \infty} f(x_n) \leq  f(x).$
A function $f: X \to [-\infty, \infty]$ is said to be \emph{lower semicontinuous} (l.s.c.) if for every $c \in \Re$, its lower level set $\{ x \in X \mid f(x) \leq c\}$ is closed in $X$. Equivalently, $f$ is l.s.c.\ if and only if for any sequence $\{x_n\}$ in $X$ converging to some $x \in X$, we have $\liminf_{n \to \infty} f(x_n) \geq  f(x).$

Let $X$ and $Y$ be separable metrizable topological spaces. Let the topology on the space $\P(Y)$ of Borel probability measures on $Y$ be the weak topology. 
A stochastic kernel $\kappa(dy \!\mid x)$ on $Y$ given $X$ is \emph{continuous} if the function $\kappa(dy \!\mid \cdot): X \to \P(Y)$ is continuous. (In some literature such a kernel is said to be weakly continuous; see e.g., \cite{FKZ12}.) Similarly, if restricted to a subset $B \subset X$, the function $\kappa(dy \!\mid \cdot): B \to \P(Y)$ is continuous, we say $\kappa(dy \!\mid x)$ is \emph{continuous on $B$}.

\subsection{Upper Semicontinuous Models}

We consider the upper semicontinuous model as defined in \cite[Def.\ 8.8]{bs}.
Here, in addition to the model assumptions given in Section~\ref{sec2.1}, we assume that:\moveup\moveup
\begin{itemize}
\item[(a)]
The control constraint set $\Gamma$ is an open subset of $S \times C$.\moveup
\item[(b)] The state transition stochastic kernel $q(dx' \!\mid x, u)$ is continuous.\moveup
\item[(c)] The one-stage cost function $g$ is u.s.c.\ on $\Gamma$ and bounded above.\footnote{By $g$ being u.s.c.\ on $\Gamma$ we mean that $g$ is u.s.c.\ on the space $\Gamma$ endowed with the relative topology---relative to the topology on $S \times C$.}\moveup 
\moveup
\end{itemize}
It is known that under (D)(N), the optimal cost function $J^*$ is u.s.c. 
Starting with $J = \0$ for (N) and with any bounded u.s.c.\ function $J$ for (D), value iteration generates u.s.c.\ functions 
$T^k(J)$ converging to $J^*$. Moreover, there exists an $\epsilon$-optimal, nonrandomized Borel measurable policy which is stationary under (D) and semi-Markov under (N). (For these optimality results, see~\cite[Props.\ 8.7, 9.21]{bs}.)

Consider the mixed value and policy iteration algorithm~(\ref{alg-1a})-(\ref{alg-1b}).
By a selection theorem for u.s.c.\ functions~\cite[Prop.\ 7.34]{bs}, if $Q: \Gamma \to [-\infty, \infty]$ is u.s.c., then the function resulting from partial minimization,
$$M(Q)(x) = \inf_{u \in U(x)} Q(x, u), \qquad x \in S,$$ 
is u.s.c., and for any $\epsilon > 0$,
there exists a Borel measurable nonrandomized stationary policy $\mu$ such that for all $x \in S$,
\begin{align}
  Q \big(x, \mu(x) \big) & \leq \begin{cases}
           M(Q)(x) + \epsilon \ \ \  & \text{if}   \ M(Q)(x)  > -\infty, \\
           - 1/\epsilon & \text{if}  \  M(Q)(x) = - \infty.
           \end{cases} \label{eq-selmu-upp}
\end{align}
We are interested in applying the above selection theorem in place of the one given in Eqs.~(\ref{eq-selmu1})-(\ref{eq-selmu2}) to derive policy iteration-like algorithms%footnote starts 
\footnote{As discussed at the end of Section~\ref{sec2.3}, in the upper (or lower) semicontinuous model we consider, standard policy iteration also encounters measurability difficulties. Because under the model conditions here, a Borel measurable policy need not have an u.s.c.\ (or l.s.c.) cost function, so at the policy improvement step, we cannot apply the selection theorem just mentioned (or the one in the next subsection) to ensure that policy iteration can avoid the measurability difficulty discussed in Section~\ref{sec2.3} by generating only Borel measurable policies.\label{footnote-semicont}}
%footnote ends 
that are similar to Algorithm II in Section~\ref{sec3.2} but operate within the family of u.s.c.\ functions. 
One way to achieve this is to choose, at each iteration of~(\ref{alg-1a})-(\ref{alg-1b}), for a given stationary Borel measurable policy $\mu$, an appropriate set $B \subset S$ to form the parameter $\k=(\mu, B)$ in the mapping $F_\k$, as the following result shows.

\mysmallskip
\begin{prop}[Upper Semicontinuous Models] \label{prp-Uppsc}
Let $\k=(\mu, B)$ for an open subset $B$ of $S$ and a Borel measurable stationary policy $\mu$ such that $\mu(du \!\mid \cdot)$ is continuous on $B$. Then for any functions $J, Q$ that are u.s.c.\ and bounded above,
$F_\k(Q\,; J)$ 
is u.s.c.\ and bounded above.
\end{prop}

\begin{proof}
Since $g$ is u.s.c.\ and bounded above by our model assumption, to show that $F_\k(Q\,; J)$ is u.s.c.\ and bounded above, it suffices to show that the sum of the two integral terms in the definition (\ref{def-F}) of $F_\k(Q\,; J)$ is u.s.c. and bounded above. 
To this end, let us rewrite this sum as 
\begin{equation} \label{eq-prfUpp1}
  \alpha \int_S \big( \phi(x') \cdot \b1_B(x') + J(x') \cdot \b1_{S \setminus B}(x') \big) \, q(dx' \!\mid x, u),
\end{equation}  
where the function $\phi(x')$ is given by
\begin{equation} \label{eq-prfUpp2}
 \phi(x') = \int_C \min \{J(x'), Q(x',u')\} \, \mu(du'\!\mid x'), \quad x' \in S. 
\end{equation} 
We prove first that $\phi(x')$ is u.s.c.\ on $B$. 
Since $J, Q$ are u.s.c.\ and bounded above,
the function $\min \{J(x), Q(x,u)\}$ is u.s.c.\ and bounded above on $\Gamma$. Note that since $\Gamma$ is an open subset of $S \times C$, we may extend the function $\min \{J(x), Q(x,u)\}$ to an u.s.c.\ function on $S \times C$ that is bounded above, and view the integral defining $\phi(x')$ as the integral of this extension. This will not change the value $\phi(x')$, since $\mu$ satisfies the control constraint. 
Then since the function $x \mapsto \mu(du\!\mid x)$ is continuous on $B$, we can apply \cite[Prop.\ 7.31(b)]{bs} to conclude that $\phi(x')$ is u.s.c.\ on $B$.

Denote $\psi(x') = \phi(x') \cdot \b1_B(x') + J(x') \cdot \b1_{S \setminus B}(x') $ for $x' \in S$. We prove that $\psi(x')$
is u.s.c.\ on $S$. Consider any sequence $\{x_n\}$ in $S$ converging to some $\bar x \in S$.
By the definition of $\phi$, we have $\phi(x') \leq J(x')$ for all $x' \in S$.  
Therefore, in the case $\bar x \not\in B$, we have 
$$ \limsup_{n \to \infty} \psi(x_n) \leq \limsup_{n \to \infty}  J(x_n) \leq J(\bar x) = \psi(\bar x),$$
where the second inequality follows from the u.s.c.\ property of $J$.
In the other case $\bar x \in B$, since $B$ is open, $x_n \in B$ for $n$ sufficiently large, so we have
$$ \limsup_{n \to \infty} \psi(x_n) =  \limsup_{n \to \infty} \phi(x_n) \leq \phi(\bar x) = \psi(\bar x),$$
where the inequality holds since $\phi$ restricted to $B$ is u.s.c., as we proved earlier.
This proves that the function $\psi$ is u.s.c.  Clearly $\psi$ is bounded above. 
Then, using also the model assumption that the stochastic kernel $q(dx' \!\mid x, u)$ is continuous, we have, by \cite[Prop.\ 7.31(b)]{bs}, that the integral~(\ref{eq-prfUpp1}) as a function of $(x,u)$ is u.s.c.\ and clearly bounded above. 
This proves the proposition.
\end{proof}
%\mysmallskip

\begin{rem}
In the preceding proposition, we need an \emph{open} set $B \subset S$ such that restricted to $B$, the function $x \mapsto \mu(du\!\mid x)$ is continuous. The empty set $B = \emptyset$ trivially fulfills this requirement, but then $F_\k(Q\,; J)$ does not depend on the policy $\mu$ at all, and this does not serve our purpose of using policy iteration-like algorithms as an alternative to value iteration. As emphasized earlier, it is desirable to have a ``large'' set $B$. Lusin's theorem~\cite[Thm.\ 7.5.2]{Dud02} tells us that there exists a \emph{closed} set $\bar B \subset S$ restricted to which, the function $\mu(du\!\mid \cdot)$ is continuous, where $\bar B$ can be chosen to be very ``large,'' with its measure arbitrarily close to $1$ for any given Borel probability measure on $S$. We can let $B = \text{int}(\bar B)$, the interior of $\bar B$; however, despite the ``large'' size of $\bar B$, $\text{int}(\bar B)$ may still be empty. Because of this, the preceding proposition is not as powerful as its counterpart, the subsequent Prop.~\ref{prp-Lowsc}, for the lower semicontinuous model to be discussed shortly, where Lusin's theorem can be directly applied to yield a large set $B$ with the desired continuity property.\qed
\end{rem}
\mysmallskip

Based on Prop.~\ref{prp-Uppsc}, we see that to keep the iterates $J_k, Q_k$ of the basic algorithm~(\ref{alg-1a})-(\ref{alg-1b}) within the set of functions that are u.s.c.\ and bounded above, we can start with $J_0, Q_0$ that are u.s.c.\ and bounded above, choose the parameters $\k_{k}=(\mu_{k}, B_{k})$ to satisfy the condition in Prop.~\ref{prp-Uppsc}, and update $Q_{k+1}$ using always the first rule in (\ref{alg-1a}), thereby resulting in u.s.c.\ functions $Q_{k+1}$ and $J_{k+1}$. For cases (D)(N), it is not hard to show that the second rule in (\ref{alg-1a}), $Q_{k+1} = Q_{\k_k, J_k}$, also makes $Q_{k+1}$ u.s.c.\ and therefore can be used. (For case (P), however, we do not know if $Q_{\k_k, J_k}$ is u.s.c.\ in general.) Then, combined with the selection of $\mu_k$ based on Eq.~(\ref{eq-selmu-upp}) for $Q=Q_k$, we obtain, similar to Algorithm II in Section~\ref{sec3.2}, policy iteration-like algorithms that operate with Borel measurable policies $\{\mu_k\}$ and u.s.c.\ functions $\{(J_k,Q_k)\}$.

We conclude this subsection with an optimality result for the upper semicontinuous model in case (P). To our knowledge, here there is no guarantee that $J^*$ is u.s.c.; however, an application of Theorem~\ref{thm-valite-P} shows the following.

\mysmallskip
\begin{prop}[Case (P) in Upper Semicontinuous Models]
Suppose that $J^*$ is bounded above and for some open set $B \subset S$ and $\delta > 0$, $B \supset \{ x \in S \mid J^*(x) < \delta \}$ and $J^*$ is u.s.c.\ on $B$. Then $J^*$ is u.s.c.\ and for any $\epsilon > 0$, there exists an $\epsilon$-optimal, Borel measurable nonrandomized Markov policy.
\end{prop}

\begin{proof}
Let $a = \sup_{x \in S} J^*(x) < \infty$.
Let $J(x) = J^*(x)$ if $x \in B$ and $J(x) = a$ otherwise. Since $J^*$ is u.s.c.\ on the open set $B$, $J$ is by definition u.s.c.\ and bounded above.
Consequently, for all $k$, $T^k(J)$ is u.s.c.\ and bounded above by \cite[Props.\ 7.31, 7.34]{bs}. 
We also have $J^* \leq J \leq c J^*$ for $c \geq \max\{ 1, a/\delta\}$, so by Theorem~\ref{thm-valite-P}(b), $T^k(J) \to J^*$. 
Then, using the fact that $T^k(J)$ is u.s.c.\ and $T^k(J) \geq J^*$ for all $k$, it follows that $J^*$ is u.s.c.%footnote starts
\footnote{Here we used the fact that if $\{f_n\}$ is a sequence of u.s.c.\ functions on a metrizable space $X$ converging pointwise to $f$ with $f_n \geq f$ for all $n$, then $f$ is u.s.c. To see this, let $\{x_k\}$ be a sequence in $X$ converging to $x \in X$. We have for every $n$,
$ \limsup_{k \to \infty} f(x_k) \leq \limsup_{k \to \infty} f_n(x_k) \leq f_n(x)$, and hence $\limsup_{k \to \infty} f(x_k) \leq \lim_{n \to \infty} f_n(x) = f(x)$. This shows that $f$ is u.s.c.
} 
%footnote ends
The assertion of the existence of an $\epsilon$-optimal, Borel measurable nonrandomized Markov policy then follows from a selection theorem for u.s.c.\ functions~(\cite[Prop.\ 7.34]{bs}; cf.\ Eq.~(\ref{eq-selmu-upp})) and the same proof argument as that for \cite[Prop.\ 9.19(P)]{bs}.
\end{proof}
%\smallskip

\subsection{Lower Semicontinuous Models}

We now consider the lower semicontinuous model from \cite[Def.\ 8.7]{bs}.
For simplicity, in addition to the model assumptions given in Section~\ref{sec2.1}, let us assume that:\moveup\moveup
\begin{itemize}
\item[(a)]
The control space $C$ is compact, and the control constraint set $\Gamma$ is a closed subset of $S \times C$.\moveup
\item[(b)] The state transition stochastic kernel $q(dx' \!\mid x, u)$ is continuous.\moveup
\item[(c)] The one-stage cost function $g$ is l.s.c.\ on $\Gamma$ and bounded below.\moveup\moveup
\end{itemize}
This is a special case of the model defined in \cite[Def.\ 8.7]{bs}, but our discussion below applies to that more general model. 
Let us also mention that there have been substantial efforts in the literature to weaken the assumptions (a) and (c) above.
For these more general lower semicontinuous models and the most recent results, we refer to the paper by Feinberg, Kasyanov and Zadoianchuk \cite{FKZ12}. In principle, the approach we describe here is applicable to these models as well to address the measurability issues in standard policy iteration (cf.\ the discussion at the end of Section~\ref{sec2.3} and Footnote~\ref{footnote-semicont}), although the subject is beyond the scope of the present paper.

It is known that under the assumptions (a)-(c) above, the optimal cost function $J^*$ is l.s.c.\ for the models (D)(P). 
Starting with $J = \0$ for (P) and with any bounded l.s.c.\ function $J$ for (D), value iteration generates l.s.c.\ functions 
$T^k(J)$ converging to $J^*$. There exists an optimal, Borel measurable nonrandomized stationary policy under (D)(P). (For these optimality results, see~\cite[Prop.\ 8.6 and Cor.\ 9.17.2]{bs}.)

Consider the mixed value and policy iteration algorithm~(\ref{alg-1a})-(\ref{alg-1b}). In what follows, we apply arguments similar to those for the upper semicontinuous model, and we obtain policy iteration-like algorithms that keep iterates $(J_k, Q_k)$ within the set of l.s.c.\ functions.
Specifically, by a selection theorem for l.s.c.\ functions~\cite[Prop.\ 7.33]{bs}, we have that if $Q: \Gamma \to [-\infty, \infty]$ is l.s.c., then the function $M(Q)(x) = \inf_{u \in U(x)} Q(x,u)$ is l.s.c.\ on $S$ and for any $\epsilon > 0$,
there exists a Borel measurable nonrandomized stationary policy $\mu$ such that
\begin{equation} \label{eq-lscsel}
   Q(x, \mu(x))  = M(Q)(x), \qquad x \in S.
\end{equation}   
Thus, if at the $k$th iteration of the algorithm~(\ref{alg-1a})-(\ref{alg-1b}), the function $Q_k$ is l.s.c., then we can choose a Borel measurable policy $\mu_k$ according to the above selection theorem with $Q=Q_k$, similar to the policy iteration-like Algorithm II in Section~\ref{sec3.2}. In order for the subsequent iterates $Q_{k+1}, J_{k+1}$ to be l.s.c., we can choose an appropriate set $B_k$ based on the following result, when forming the parameter $\k_k=(\mu_k, B_k)$ for the mapping $F_{\k_k}$ in the algorithm.

Let $\mu$ be a Borel measurable stationary policy. 
We know from Lusin's theorem~\cite[Thm.\ 7.5.2]{Dud02} that there exists a closed subset $B \subset S$ such that restricted to $B$, the function $x \mapsto \mu(du\!\mid x)$ is continuous, and moreover, for any given $p \in \P(S)$, the set $B$ can be chosen to have $p(B)$ arbitrarily close to $1$. 

%\mysmallskip
\begin{prop}[Lower Semicontinuous Models] \label{prp-Lowsc}
Let $\k=(\mu, B)$ for a closed subset $B$ of $S$ and a Borel measurable stationary policy $\mu$ such that $\mu(du \!\mid \cdot)$ is continuous on $B$. Then for any functions $J, Q$ that are l.s.c.\ and bounded below, 
$F_\k(Q\,; J)$ is l.s.c.\ and bounded below.
\end{prop}

\begin{proof}
Similar to the proof of Prop.~\ref{prp-Uppsc}, it suffices to show that the integral (\ref{eq-prfUpp1}) as a function of $(x,u)$ is l.s.c.\ and bounded below on $\Gamma$.
We prove first that the function $\phi(x')$ given by Eq.~(\ref{eq-prfUpp2}) is l.s.c.\ on $B$. 
Since $J, Q$ are l.s.c.\ and bounded below,
the function $\min \{J(x), Q(x,u)\}$ is l.s.c.\ and bounded below on $\Gamma$. We may extend the function $\min \{J(x), Q(x,u)\}$ to an l.s.c.\ function on $S \times C$ that is bounded below, by defining its values outside $\Gamma$ to be $\infty$, and we can view the integral defining $\phi(x')$ as the integral of this extension. This will not change the value $\phi(x')$, since $\mu$ satisfies the control constraint. 
Then, since the function $x \mapsto \mu(du\!\mid x)$ is continuous on $B$, we can apply \cite[Prop.\ 7.31(a)]{bs} to conclude that $\phi(x')$ is l.s.c.\ and bounded below on $B$.

Denote $\psi(x') = \phi(x') \cdot \b1_B(x') + J(x') \cdot \b1_{S \setminus B}(x') $ for $x' \in S$. We prove that $\psi(x')$
is l.s.c.\ on $S$. Consider any sequence $\{x_n\}$ in $S$ converging to some $\bar x \in S$.
If $\bar x \not\in B$, then since $S \setminus B$ is open, we have $x_n \in S \setminus B$ for all $n$ sufficiently large, so
$$ \liminf_{n \to \infty} \psi(x_n) = \liminf_{n \to \infty}  J(x_n) \geq J(\bar x) = \psi(\bar x),$$
where the inequality follows from the l.s.c.\ property of $J$.
Suppose now $\bar x \in B$. There exists a subsequence $\{x_{n_i}\}$ of $\{x_n\}$ such that $\liminf_{n \to \infty} \psi(x_n) = \lim_{i \to \infty} \psi(x_{n_i})$ and either (i) $x_{n_i} \in B$ for all $i$ or (ii) $x_{n_i} \not\in B$ for all $i$.
Then in case (i), we have
$$\liminf_{n \to \infty} \psi(x_n) =  \lim_{i \to \infty} \psi(x_{n_i}) =  \lim_{i \to \infty} \phi(x_{n_i}) \geq \phi(\bar x) = \psi(\bar x),$$
where the inequality holds since $\phi$ restricted to $B$ is l.s.c., as we proved earlier.
In case (ii), we have
$$ \liminf_{n \to \infty} \psi(x_n) = \lim_{i \to \infty} \psi(x_{n_i}) =  \lim_{i \to \infty} J(x_{n_i}) \geq J(\bar x) \geq \phi(\bar x) = \psi(\bar x),$$ 
where the first inequality holds since $J$ is l.s.c., and the second inequality holds since by the definition of $\phi$, we have $\phi(x') \leq J(x')$ for all $x' \in S$.  Thus we have proved that the function $\psi$ is l.s.c.  Clearly $\psi$ is bounded below. 
Then, using also the fact that the stochastic kernel $q(dx' \!\mid x, u)$ is continuous, we have, by \cite[Prop.\ 7.31(a)]{bs}, that the integral~(\ref{eq-prfUpp1}) as a function of $(x,u)$ is l.s.c.\ and bounded below. This proves the proposition.
\end{proof}
%\mysmallskip

Based on Prop.~\ref{prp-Lowsc}, to obtain policy iteration-like algorithms that keep the iterates $J_k, Q_k$ within the set of functions that are l.s.c.\ and bounded below, we can specialize the basic algorithm~(\ref{alg-1a})-(\ref{alg-1b}) as follows:
Start with $J_0, Q_0$ that are l.s.c.\ and bounded below. For $k \geq 0$, form the parameter $\k_{k}=(\mu_{k}, B_{k})$ by first choosing a Borel measurable $\mu_k$ based on the selection theorem (\ref{eq-lscsel}) with $Q=Q_k$, and next choosing a closed set $B_k$ based on Lusin's theorem as described earlier. Then update $Q_{k+1}$ using the first rule in (\ref{alg-1a}). This results in l.s.c.\ functions $Q_{k+1}$ and $J_{k+1}$. For cases (D)(P), the second rule in (\ref{alg-1a}), $Q_{k+1} = Q_{\k_k, J_k}$, can also be shown to make $Q_{k+1}$ l.s.c.\ and therefore can be used. For case (N), however, we do not know if $Q_{\k_k, J_k}$ is l.s.c.\ in general. 

\section{Concluding Remarks} \label{sec7}
\markboth{\rm \S \thesection. Concluding Remarks}{\rm \S \thesection. Concluding Remarks}

In this paper we have addressed the long-standing open issue of constructing a valid policy iteration algorithm for total cost Borel-space stochastic DP with universally measurable policies. Our approach is based on a mixed value and policy iteration idea. It makes critical use of the fact that any universally measurable policy has Borel measurable portions, to maintain cost function iterates within the set of lower semi-analytic functions. 
It employs an algorithmic framework that combines the characteristics of both value and policy iteration, to allow stationary policies to be used in computing the optimal cost function. 
Our approach can also address similar policy iteration issues that arise in upper and lower semicontinuous models. By choosing algorithmic parameters accordingly, we have shown how to obtain policy iteration-like algorithms that can keep the cost function iterates within the desired family of semicontinuous functions.

Besides its ability to handle measurability or non-measurability related structural restrictions in stochastic DP problems, our mixed value and policy iteration approach has attractive computational features. We have discussed them in Sections~\ref{sec1} and \ref{sec3.2}; here we emphasize again two specific advantages. First, it yields policy iteration-like algorithms that have comparable computation complexity to standard policy iteration (which is usually much more efficient than value iteration) and yet have much more reliable convergence properties. Second, it works naturally with approximate policy improvement, and is capable of solving difficult problems that cannot be handled by standard policy iteration, such as problems whose optimal or near-optimal policies are non-stationary.

A number of issues are still not fully understood and deserve further research, however. One important question is how to choose the parameters $(\mu, B)$, especially the set $B$, for the policy iteration-like algorithms. We have discussed in the paper two possible choices of $B$ for a given policy $\mu$. It is still unclear how to choose $B$ to facilitate fast convergence of the algorithms. It is also desirable to be able to find such sets $B$ easily
so that the computation involved does not tamper with the overall efficiency of the method.

Another question relates to finding near-optimal policies, especially in problems with non-stationary near-optimal policies.
Our approach yields function sequences that converge pointwise to the optimal cost function for discounted, nonpositive and nonnegative cost models. It also yields asymptotically optimal policy sequences for discounted cost models. For nonpositive and nonnegative cost models, it seems plausible that we can extract nearly optimal policies from the data produced by the algorithms, but the details still need to be worked out. 

For nonnegative cost models, we have derived a new sufficient condition for convergence of value iteration, which involves the initial function only. It applies to all nonnegative models (countable space or uncountable Borel space models), and it provides, in addition, a new characterization of the set of functions within which the optimal cost function is the unique solution of the optimality equation. 
Using this condition, our method is also shown to produce in the limit the optimal cost function when initialized properly. Obtaining useful initial functions satisfying this condition is generally an open question at present, however. 
Moreover, in the literature on average cost problems, it has been shown that one can speed up the convergence of (relative) value iteration significantly, by choosing an initial function that qualitatively matches the shape of the optimal (relative) value function (see Chen and Meyn~\cite{Chm-vi} and Meyn \cite[Example 9.5.7]{Mey08}). Whether the same can be said for total cost problems is another interesting open question, which we aim to address in the future.

We conclude the paper with a discussion of some other applications of our approach and future research directions.

\smallskip
\noindent {\bf Asynchronous computation}
\mysmallskip

One may consider asynchronous distributed computation in the framework of universally measurable policies, by combining the approach and analysis given in this paper with arguments used in our earlier works \cite{BerY-discount,YuB-ssp,BerY-adp}.
We  discuss the subject briefly here, 
focusing on issues related to universal measurability in a simplified setting.

Suppose that instead of the basic algorithm~(\ref{alg-1a})-(\ref{alg-1b}), at each iteration $k$, we only compute $Q_{k+1}(x,u)$ for a subset $\Gamma_k$ of state-control pairs in $\Gamma$ and compute $J_{k+1}(x)$ for a subset $S_k$ of states in $S$. (For the rest of the states $x$ or state-control pairs $(x,u)$, we let $J_{k+1}(x)=J_k(x)$, $Q_{k+1}(x,u)=Q_k(x,u)$.)
This is the type of operations that would be performed in a distributed computation environment, where a single processor handles only part of a computation task and processors share results with each other. 

As before, with universally measurable policies, we need to keep the iterates within the set of lower semi-analytic functions. 
To meet this requirement, we can let $S_k$ be a Borel subset of $S$ and let $\Gamma_k = R_k \cap \Gamma$, where $R_k$ is a \emph{Borel subset of $S \times C$}.  
This will keep $Q_{k+1} \in A(\Gamma)$. The reason is that if $Q, Q' \in A(\Gamma)$ and $R$ is a Borel set in $S \times C$, then the function
$$ Q \cdot \b1_{R \,\cap\, \Gamma} + Q' \cdot \b1_{ (S \times C \setminus R) \, \cap \, \Gamma}$$
is lower semi-analytic by \cite[Lemma 7.30(4)]{bs},
because $\b1_{R \, \cap \, \Gamma}(x,u)$ and $\b1_{ (S \times C \setminus R) \, \cap \,\Gamma}(x,u)$ are nonnegative Borel measurable functions on $\Gamma$.
Similarly, the reason for $J_{k+1} \in A(S)$ is that if $J, J' \in A(S)$ and $D \subset S$ is Borel, then the function $J \cdot \b1_{D} + J' \cdot \b1_{S \setminus D}$ is lower semi-analytic. 

Further studies of such asynchronous algorithms, in particular, stochastic asynchronous Q-learning versions (along the lines of our earlier works \cite{BerY-discount,YuB-ssp}), are important subjects for future investigation.

\smallskip
\noindent {\bf More elaborate variants}
\mysmallskip

In this paper we have focused on the mappings $F_\k, \k \in \K$, defined by~(\ref{def-F}), where we partition the state space into two subsets. The same idea leads to more elaborate mappings, which can also be used in the mixed value and policy iteration approach. We give one such example here, in which we partition the state-control space $S \times C$.

For a stationary universally measurable policy $\mu$, let $R\subset S \times C$ be a Borel set such that $B = \text{proj}_S(R)$ is Borel and the function $x \mapsto \mu(du \mid x)$ is Borel measurable on $B$.
For any such pair $\hat \k = (\mu, R)$, 
we may consider a mapping $F_{\hat \k}$ defined as: for all $J \in A(S), Q \in A(\Gamma)$,
\begin{align}
  F_{\hat \k} ( Q \, ; J)(x, u) & =  g(x, u) + \alpha \int_S \int_C \b1_{S \times C \setminus R}(x', u') \cdot J(x')  \, \mu(du' \!\mid x') \, q(dx' \!\mid x, u) \notag \\
    & \ \ \ +  \alpha \int_S  \int_C  \b1_R(x',u') \cdot \min \big\{ J(x') \, , \, Q(x', u') \big\} \, \mu(du' \!\mid x') \, q(dx' \!\mid x, u), \quad (x, u) \in \Gamma. \notag 
\end{align}
It can be equivalently written as
\begin{align}    
F_{\hat \k} ( Q \, ; J)(x, u)  & =  g(x, u) + \alpha \int_{S\setminus B}  J(x') \, q(dx' \!\mid x, u)   + \alpha \int_B  J(x') \cdot \mu \big(C \setminus R_{x'} \!\mid x' \big) \,  q(dx' \!\mid x, u)    \notag \\
        &  \ \ \ +  \alpha \int_B  \int_{R_{x'}}  \min \big\{ J(x') \, , \, Q(x', u') \big\} \, \mu(du' \!\mid x') \, q(dx' \!\mid x, u), \qquad (x, u) \in \Gamma,  \label{def-Fvar}
 \end{align}
where $B = \text{proj}_S(R)$ and $R_x = \{ u \in C \!\mid (x, u) \in R \}$ is the vertical section of $R$ at $x$.
That the function $F_{\hat \k}(Q\,; J)$ is lower semi-analytic can be established similarly to Prop.~\ref{prp-F1}(a), using the arguments in its proof, together with the fact that restricted to $B$, $\mu(C \setminus R_{x'} \!\mid x')$ is a nonnegative Borel measurable function \cite[Cor.\ 7.26.1]{bs} and hence the term $\int_B  J(x') \cdot \mu\big(C \setminus R_{x'} \!\mid x' \big) \,  q(dx' \!\mid x, u)$ in (\ref{def-Fvar}) as a function of $(x,u)$ is lower semi-analytic.

Let us also mention that all the algorithms we gave earlier use the formula $J_{k+1} = M(Q_{k+1})$ to compute the cost function iterates. Suppose in a given problem one knows that $\underline{J} \leq J^* \leq \bar {J}$ for certain $\underline{J}, \bar {J} \in A(S)$. Then it is natural and more efficient to constrain $J_{k+1}$ to be within these upper and lower bounds of the optimal cost function by setting $J_{k+1} = \max \big\{ \underline{J}, \min \{ \bar J, M(Q_{k+1}) \} \big\}$. 

\smallskip
\noindent {\bf Extensions to other models}
\mysmallskip

Finally, we note that while we have focused on the three classical total cost problems in this paper, the technique we used to handle the measurability issues in policy iteration can be applied to other types of stochastic control problems. 
These include partially observable problems, which can be reduced to equivalent completely observable problems \cite[Chap.\ 10]{bs}, and also include discounted problems with unbounded one-stage costs, and undiscounted total cost problems without sign constraints on the one-stage costs.
Convergence properties of the mixed value and policy iteration method for such models are worthy of further study. 
Extensions to average cost problems and constrained problems are also important subjects for future research.

\section*{Acknowledgments}
We thank Prof.\ Steven Shreve for a helpful discussion and a suggestion about how to choose the probability measures for the algorithms in Section~\ref{sec3.2}, which we described in Example~\ref{ex-choiceB1}.
We also thank Prof.\ Eugene Feinberg, with whom our recent correspondence about Borel models stimulated this research, and Prof.\ Sean Meyn, who pointed us to related work in his book~\cite{Mey08}, which motivated the alternative proof for Theorem~\ref{thm-valite-P} we gave in Appendix~\ref{appsec-altprfP}.
We appreciate Prof.\ Sanjoy Mitter's helpful feedback on our early draft, and we owe many thanks to the anonymous reviewers for their critical and suggestive comments about our manuscript, including the suggestion of Example~\ref{ex-choiceB2}. 
This work was done when the first author was with the Laboratory for Information and Decision Systems, MIT. The research of both authors was supported by the Air Force Grant FA9550-10-1-0412.
 
\addcontentsline{toc}{section}{References} 
\markboth{\rm References}{\rm References}
\bibliographystyle{plain} 
\let\oldbibliography\thebibliography
\renewcommand{\thebibliography}[1]{%
  \oldbibliography{#1}%
  \setlength{\itemsep}{0pt}%
}
{\fontsize{9}{11} \selectfont
\bibliography{polite_borelDP}}
 
\clearpage 
\begin{appendices}

% add "Appendix" to the section heading
\newcommand{\appsection}[1]{\let\oldthesection\thesection
  \renewcommand{\thesection}{%Appendix 
  \oldthesection}
  \section{#1}\let\thesection\oldthesection}

\appsection{Optimal Stopping Problems Associated with $F_\k$} \label{appsec-optstop}
\markboth{\rm Appendix \thesection. Optimal Stopping Problems Associated with $F_\k$}{\rm Appendix \thesection. Optimal Stopping Problems Associated with $F_\k$}

In this appendix, for a given $\k=(\mu, B) \in \K$, $J \in A(S)$, and a control problem of type (D), (N) or (P),
we formulate an associated optimal stopping problem of the same type.
We establish the relation between its optimal cost function and 
the pointwise limit $Q_{\k, J} = \lim_{k \to \infty} F_\k^k(\0 \,; J)$, and we show that the mapping $F_\k(\cdot \, ; J)$ can be viewed as a form of the optimal cost operator for this problem. 
(Other formulations of the optimal stopping problem are also possible and equivalent for our purpose. We will focus only on one here.)
In addition we describe a linear program in case (P) and show that under certain conditions, it yields an upper bound on $Q_{\k,J}$ that can be used in a mixed value and policy iteration algorithm discussed in Section~\ref{secP5.2}. 

\subsection{Formulation} \label{secA.1}

As before we assume that the given function $J$ is such that $J \in A_b(S)$ in case (D),  $J \in A_-(S)$ in case (N), and $J \in A_+(S)$ in case (P). 
The function $J$ will define the stopping costs, whereas the policy $\mu$ will be used to define the dynamics of the unstopped process.

\medskip
{%\samepage
\noindent {\bf Optimal Stopping Problem Associated with $J$ and $(\mu, B) \in \K$}\moveup\moveup
\begin{itemize}
\item State space $S^o = \big(S\times C) \cup \{ \infty \}$, with $\infty$ representing an absorbing, cost-free state. (The topology on $S^o$ consists of the open sets in $S \times C$, the set $\{\infty\}$ and their unions.)\moveup
\item Control space $C^o= \{0, 1\}$, with $0$ representing ``to stop'' and $1$ ``to continue.''\moveup
\item Control constraint: $U^o(\infty) = \{0,1\}$ and
$$ U^o\big((x, u)\big) = \{0,1\} \  \text{on} \  B \times C, \qquad U^o\big((x, u)\big) = \{0\} \  \text{on} \  (S \setminus B) \times C. \vspace*{-0.03in} $$
\item One-stage costs: $g^o(\infty,0) = g^o(\infty,1) = 0$ and
\begin{align*} 
g^o\big((x, u),0 \big) & = J(x) \quad \ \ \ \forall \, (x, u) \in S \times C, \\
   g^o \big( (x, u), 1\big) & = g(x, u)\quad \ \! \forall \, (x,u) \in (B \times C) \cap \Gamma, \\
  g^o\big((x, u), 1) & = K \ \  \ \quad \ \ \ \, \forall \, (x,u) \in (B \times C) \setminus \Gamma,
\end{align*}
where $K=0$ for (N), $K = +\infty$ for (P), and $K \geq \max \{ \| g\|_\infty, \|J \|_\infty\}$ for (D).\moveup
\item State transition stochastic kernel $q^o(\cdot \mid \cdot)$ on $S^o$ given $S^o \times C^o$:
for any Borel set $D \subset S^o$ and any $(x, u) \in S \times C$,
$q^o \big(D \!\mid \infty, 0)  = q^o \big(D \!\mid \infty, 1) = q^o \big(D \!\mid (x, u), 0 \big) = \delta_\infty(D),$ and
$$   q^o \big(D \!\mid (x, u), 1 \big)  = \int_S \int_C \b1_{D\setminus \{\infty\}}\big((x', u')\big) \,\tilde \mu(du' \!\mid x') \, q(dx' \!\mid x, u), $$
where $\tilde \mu(du' \!\mid x')$ is a Borel measurable stochastic kernel on $C$ given $S$ 
such that $\tilde \mu(du' \!\mid x') = \mu(du' \!\mid x')$ for all $x' \in B.$\moveup%\moveup
\end{itemize}
}\mysmallskip

The state transition kernel required in the above 
can be chosen, for instance, by letting
$\tilde \mu(du' \!\mid x') = \mu(du' \!\mid x')$ for $x' \in B$ and $\tilde \mu(du' \!\mid x') = p(du')$ for $x' \not\in B$, where $p$ is any Borel probability measure on $C$.
Note also that with the control~$1$, for any $(x, u) \in S \times C$ and Borel $D \subset S^o$, we have
\begin{align}
    q^o \big(D \!\mid (x, u), 1 \big) & =  \int_B \int_C \b1_{D\setminus \{\infty\}}\big((x', u')\big)  \mu(du' \!\mid x') \, q(dx' \!\mid x, u) \notag \\
    &  \ \ \ \ + \int_{S \setminus B} \int_C \b1_{D\setminus \{\infty\}}\big((x', u')\big) \,\tilde \mu(du' \!\mid x') \, q(dx' \!\mid x, u). \label{eq-stopmodel-q}
\end{align}     

The above formulation fits the general stochastic control model described in Section~\ref{sec2.1}.
In particular, the graph of the control constraint $U^o$ is an analytic set, the one-stage cost function $g^o$ is lower semi-analytic, and the state transition kernel $q^o$ is Borel measurable.
For a state $z \in S^o$, we denote the cost of a universally measurable policy $\pi^o$ by $V_{\pi^o}( z)$. 
It is as defined in Section~\ref{sec2.1} and can be expressed as follows.
For $k \geq 0$, let $(z_k, u^o_k)$ denote the state and control at time $k$, and let $\tau$ be the time when the process is stopped: $\tau = \min \big\{ k \geq 0 \mid u_k^o = 0 \big\}$ where we define $\tau = \infty$ if $\big\{ k \geq 0 \mid u_k^o = 0 \big\} = \emptyset$.
For each $k$, let $(x_k,u_k) = z_k$ if $z_k \in S \times C$, and let $(x_k,u_k)$ equal some fixed state in $S \times C$ if $z_k = \infty$. Then for $z_0=(x,u) \in S \times C$, $V_{\pi^o}\big( (x, u) \big)$ can be expressed as
\begin{equation}
  V_{\pi^o}\big( (x, u) \big) = \E^{\pi^o} \left\{ \sum_{k=0}^\infty \alpha^k g^o\big( z_k, u^o_k\big) \right \} = \E^{\pi^o} \left\{ \sum_{k=0}^{\tau -1} \alpha^k g^o\big( (x_k,u_k), 1\big)  + \alpha^{\tau} J(x_\tau) \right \}. \label{eq-expV0}
\end{equation}  
Note that in the above, $(x_k,u_k)$ is meaningfully defined on $\{ \tau \geq k \}$.

Denote the optimal cost function by $V^*$ and the optimal cost operator by $T_o$. 
By the theory for (D)(N)(P) in the case where the number of controls at each state is finite \cite[Props.\ 9.8, 9.14, Cor.\ 9.17.1]{bs},
we have 
$$T_o^k(\0) \to V^*, \qquad V^* = T_o \big( V^* \big),$$
where the convergence is monotonic for (N)(P),
and furthermore, for (D)(P), there exists an optimal nonrandomized stationary policy. 
(In case (N), an optimal policy need not exist even when the control space is finite; see \cite[Ex.\ 4.1, p.\ 181]{Ber-DP13} for such an example.)

Let $V_k = T_o^k(\0)$, $k \geq 0$, be the optimal $k$-stage cost functions.
To simplify notation we will write $V(x,u)$ for $V\big((x,u)\big)$.
Clearly, for the absorbing state $\infty$ and for the states in $(S \setminus B) \times C$, where the only control is to stop, 
we have for all $k \geq 1$,
\begin{align} 
  V^*(\infty)  & =  V_k(\infty) = 0, 
    \qquad V^* (x,u)  = V_k (x,u) = J(x), \quad  \forall \, (x, u) \in (S \setminus B) \times C. \label{eq-stop5}
\end{align}   

Next we will calculate the optimal costs for states in the set $(B \times C) \cap \Gamma$ and relate the results to $Q_{\k,J}$ and $F_\k(\cdot \,; J)$. For our purposes, the set $(B \times C) \setminus \Gamma$ of states can be ignored, not only because they are outside the control constraint set of the original problem, but also because in the optimal stopping problem, they are formulated to be unreachable 
from the rest of the states.
In particular, if the starting state $(x, u)$ is in $(B \times C) \cap \Gamma$, then since the policy $\mu$ satisfies the control constraint of the original problem, we see from the expression (\ref{eq-stopmodel-q}) for the state transition probability $q^o(\cdot \!\mid (x,u),1)$ that the probability of the next state being in $(B \times C) \setminus \Gamma$ is zero.
If the starting state $(x, u)$ is in $(S \setminus B) \times C$, then the control $1$ (to continue) is not allowed according to the control constraint $U^o$, so the next state is $\infty$. Therefore, the set $(B \times C) \setminus \Gamma$ is not reachable from the rest of the states.

Since at time $k$, the continuation cost is $g^o((x_k,u_k), 1) = g(x_k,u_k)$ if $(x_k, u_k) \in (B \times C) \cap \Gamma$, 
the preceding discussion also shows that 
for each $(x,u) \in \Gamma$, 
the expected total cost of ${\pi^o}$ for the initial distribution $p^o(\cdot) = q^o(\cdot \!\mid  (x, u), 1)$, with respect to the probability measure induced by $\pi^o$ and $p^o$, is 
 \begin{equation}
   V_{\pi^o\!,p^o} =  \E^{\pi^o\!,p^o} \left\{ \sum_{k=0}^{\tau -1} \alpha^k g (x_k,u_k)  + \alpha^{\tau} J(x_\tau) \right \} \label{eq-expV}
\end{equation}
(cf.\ Eq.~(\ref{eq-expV0})).
We will use the expression (\ref{eq-expV}) later to derive an expression for $Q_{\k,J}$ (see Cor.~\ref{cor-stopmodel3}).

\subsection{Relations with $F_\k(\cdot\,;J)$ and $Q_{\k, J}$} 
We will now express the operator $T_o$ and calculate $V_k, V^*$ for the states in $(B \times C) \cap \Gamma$. 
We focus on the set of functions 
\begin{equation}
 \big\{ V \in A(S^o) \mid V(\infty) = 0, V(x,u) = J(x), \, (x, u) \in (S \setminus B) \times C \big\}, \label{eq-setV}
\end{equation} 
which includes $V^*, V_k$ (cf.\ Eq.~(\ref{eq-stop5})). 
For any $V$ in this set, we can express $T_o(V)$ as
\begin{equation} \label{eq-ToV}
  T_o(V)(x, u)  = \min \left\{ J(x) \, , \, G_V(x,u) \right\},   \qquad   (x, u) \in (B \times C) \cap \Gamma,
\end{equation} 
where the first term $J(x)$ is the stopping cost, and the second term $G_V(x,u)$ is associated with the continuation action, given by
\begin{align}
    G_V(x,u) & = g(x, u) + \alpha \int_{S \times C} V(z) \, q^o \big(dz \!\mid (x, u), 1 \big) \notag \\
    & = g(x, u) + \alpha \int_{S \setminus B} J(x') \, q(dx' \!\mid x, u)  +  \alpha \int _B\int_C V(x',u') \,  \mu(du' \!\mid x') \, q(dx' \!\mid x, u). \label{eq-ToVb}
\end{align}    
In the above, to obtain Eq.~(\ref{eq-ToVb}) from the first equality, we used the expression of $q^o(dz\!\mid (x,u),1)$ given in~(\ref{eq-stopmodel-q}), which implies that for any $(x, u) \in S \times C$ and function $V$ in the set (\ref{eq-setV}),
\begin{align}
  \int_{S \times C} V(z) \, q^o \big(dz \!\mid (x, u), 1 \big)  
 &  =       \int_{S \setminus B} J(x') \, q(dx' \!\mid x, u) + \int _B\int_C V (x', u') \,  \mu(du' \!\mid x') \, q(dx' \!\mid x, u), \label{eq-stop5a} 
\end{align}   
Equations~(\ref{eq-ToV})-(\ref{eq-ToVb}) express the optimality equation $V = T_o(V)$ for $V$ in the set (\ref{eq-setV}).

Using the fact $V^* = T_o(V^*)$, we then obtain
\begin{equation} \label{eq-stop-opteq}
   V^* (x,u) = T_o \big(V^*\big)(x, u) = \min \left\{ J(x) \, , \, f^*(x,u) \right\},  \qquad  \forall \, (x, u) \in (B \times C) \cap \Gamma, 
 \end{equation}
where $f^*(x,u) = G_{V^*}(x,u)$, and it is the optimal expected future cost for continuation and can be expressed in several equivalent ways:
\begin{align} 
f^*(x, u) & = g(x, u) + \alpha \int_{S \times C} V^*(z) \, q^o \big(dz \!\mid (x, u), 1 \big)  \label{eq-f*} \\
& = g(x, u) + \alpha \int_{S \setminus B} J(x') \, q(dx' \!\mid x, u)  +  \alpha \int _B\int_C V^*(x',u') \,  \mu(du' \!\mid x') \, q(dx' \!\mid x, u)  \notag \\ 
   & = g(x, u) + \alpha \int_{S \setminus B} J(x') \, q(dx' \!\mid x, u) + \alpha \int _B\int_C \min \big\{ J(x') \, , \, f^*\big(x',u'\big) \big\} \,  \mu(du' \!\mid x') \, q(dx' \!\mid x, u). \label{eq-stop6b}
\end{align}
Here in deriving Eq.~(\ref{eq-stop6b}), 
we used the fact that for any $x' \in B$, 
\begin{equation}
   \int_C V^*(x',u') \,  \mu(du' \!\mid x')  = \int_C \min \big\{ J(x') \, , \, f^*\big(x',u'\big) \big\} \,  \mu(du' \!\mid x'), \label{eq-stop6c}
\end{equation}
which holds since $\mu$ satisfies the control constraint of the original problem: $\mu\big(U(x') \!\mid x' \big) = 1$, and for $x' \in B$ and $u' \in U(x')$, $V^*(x', u')$ can be expressed as in~(\ref{eq-stop-opteq}).

Similar to the preceding derivation, we can calculate the optimal $k$-stage cost functions $V_k$, $k \geq 1$, and
define functions $f_k$ on $(B \times C) \cap \Gamma$ associated with the continuation action, for $k \geq 0$, by
\begin{equation}
f_k(x, u)  = g(x, u) + \alpha \int_{S \times C} V_k(z) \, q^o \big(dz \!\mid (x, u), 1 \big), \qquad (x, u) \in (B \times C) \cap \Gamma, \ k \geq 0. \label{eq-fk}
\end{equation}
From the recursive relations, 
$$ V_{k+1} (x,u)  = T_o \big(V_{k}\big)(x, u) = \min \left\{ J(x) \, , \, f_k(x,u) \right\}, \qquad (x, u) \in (B \times C) \cap \Gamma, \ k \geq 0,$$
we obtain that the functions $f_k$, $k \geq 1$, satisfy the recursion (\ref{eq-stop6b}) with $f_k$ replacing $f^*$ on the left-hand side and with $f_{k-1}$ replacing $f^*$ in the right-hand side.

We recognize the expression on the right-hand side of Eq.~(\ref{eq-stop6b}) as the same expression that defines $F_\k(f^*; J)(x,u)$ (cf.\ Eq.~(\ref{def-F})).
To be more precise, since $F_\k(\cdot\,; J)$ is a mapping on $A(\Gamma)$ and $f^*$ is defined on $(B \times C) \cap \Gamma$,
we will adopt the following convention: for any function $f$ defined on $(B \times C) \cap \Gamma$, by $F_\k( f\,; J)$ we mean any $F_\k( f_e \,; J)$ where $f_e$ is an (arbitrary) extension of $f$ to $\Gamma$.
This is valid because by definition $F_\k ( Q \, ; J)$ is completely determined by the function $Q$ restricted to $(B \times C) \cap \Gamma$. 
In other words, denoting $\Gamma_B = (B \times C) \cap \Gamma$, we have
 \begin{equation}
     Q \,|_{\Gamma_B} = Q' |_{\Gamma_B} \quad \Longrightarrow \quad F_\k ( Q \, ; J) = F_\k ( Q'  ; J). \label{eq-equiv-FQ}
 \end{equation}    

Based on the equivalence between Eq.~(\ref{eq-stop6b}) and $F_\k(f^*; J)(x,u)$,  we can relate the optimal cost functions $V^*$, $V_k$ of the optimal stopping problem to the mapping $F_\k(\cdot\,; J)$ and the function $Q_{\k, J} = \lim_{k \to \infty} F_\k^k(\0\,; J)$ as follows.

\mysmallskip
\begin{lem}  \label{lma-stopmodel1}
{\rm (D)(N)(P)} \ Let $\Gamma_B = (B \times C) \cap \Gamma$, and let $f^*, f_k : \Gamma_B \to [-\infty, \infty]$, $k \geq 0$, be the minimal future cost functions associated with continuation, given by Eqs.~(\ref{eq-f*}) and (\ref{eq-fk}) respectively; in particular,  
$f_0=g |_{\Gamma_B}$. Then  
$$ f^* = F_\k \big( f^*; J \big) \big|_{\Gamma_B}, \qquad f_k = F_\k \big(f_{k-1}; J \big) \, \big|_{\Gamma_B}, \ \ \ k \geq 1,$$
and $f_k \to f^*$. Moreover,
\begin{equation} 
Q_{\k, J} \big|_{\Gamma_B} = f^*, \qquad \ Q_{\k, J} = F_\k \big(f^*; J). \label{eq-rel-Qf*}
\end{equation}
\end{lem}

\begin{proof}
The recursive relations for $f^*, f_k$ were derived earlier. The fact $f_k \to f^*$ follows from Eqs.~(\ref{eq-f*}) and (\ref{eq-fk}) by applying the bounded convergence theorem in case (D), and the monotone convergence theorem in cases (N)(P), using the convergence $V_k \to V^*$ in each of these cases.

We now prove the relation~(\ref{eq-rel-Qf*}) between $Q_{\k,J}$ and $f^*$. 
Since $f_k \to f^*$, using the relation $f_k = F_\k \big(f_{k-1}; J \big) \, \big|_{\Gamma_B}$ and Eq.~(\ref{eq-equiv-FQ}), we have $f_k=F_\k^k (g\,; J) \big|_{\Gamma_B} \to f^*$.
Suppose we have proved $F_\k^k (g; J)\to Q_{\k,J}$. Then it will follow that 
$Q_{\k, J} \big|_{\Gamma_B}   =  f^*.$
In turn, this will imply $F_\k(Q_{\k, J} ; J) = F_\k (f^*; J )$ by Eq.~(\ref{eq-equiv-FQ}), and hence $Q_{\k,J} = F_\k (f^*; J )$ since $Q_{\k,J} = F_\k(Q_{\k, J} ; J)$ by Prop.~\ref{prp-optstop-cost}.
Thus it is sufficient to prove $F_\k^k (g; J)\to Q_{\k,J}$.
For (D), this follows from the proof of Prop.~\ref{prp-optstop-cost}(D), which established the contraction property of $F_\k (\cdot; J)$. 
For (N), we have $g \leq 0$ and $J \leq 0$. By a direct calculation, $F_\k(\0\,; J) \leq g \leq 0$,
so we have, by the monotonicity of $F_\k(\cdot\,; J)$, 
$$  F_\k^k(\0\,; J)  \leq F_\k^{k-1}(g ; J) \leq F_\k^{k-1}(\0\,; J), \quad k \geq 1.$$
Since $F_\k^k(\0\,; J) \downarrow Q_{\k,J}$ by Prop.~\ref{prp-optstop-cost}, we have $F_\k^k (g; J)\downarrow Q_{\k,J}$.
The convergence $F_\k^k (g; J)\to Q_{\k,J}$ in case (P) follows from a symmetrical argument.
\end{proof}
%\mysmallskip

Based on Lemma~\ref{lma-stopmodel1}, we may view $F_\k(\cdot\,; J)$ as an optimal cost operator for the minimal future costs associated with the continuation action in the optimal stopping problem.
For states $(x,u) \in \Gamma_B$, based on the relation $Q_{\k, J} \big|_{\Gamma_B} = f^*$ in~(\ref{eq-rel-Qf*}), we can also interpret $Q_{\k, J}(x,u)$ as the minimal costs at $(x,u)$ with continuation at the first stage.

We now give several expressions of $Q_{\k, J}$ in terms of $V^*$ in the following corollary.
For each $(x,u) \in \Gamma$, we will consider the optimal stopping problem starting with an initial state distribution $p^o$ given by $q^o( \cdot \mid (x, u), 1)$, the transition distribution for $(x,u)$ under the continuation action.

\mysmallskip
\begin{cor} \label{cor-stopmodel3}
{\rm (D)(N)(P)} \ For all $(x, u) \in \Gamma$,
\begin{align*}
  Q_{\k, J}(x, u) & = g(x, u) + \alpha \int_{S\times C} V^*(z) \, q^o(dz \mid (x, u), 1). 
\end{align*}
In particular, if in the optimal stopping problem associated with $(\k, J)$, an optimal policy $\pi^{o*}$ exists (as is true under (D)(P)), then for all $(x, u) \in \Gamma$,
\begin{align*}
   Q_{\k, J}(x, u)  = g(x, u) + \alpha \, \E^{\pi^{o*}\!, p^o} \left\{ \sum_{k=0}^{\tau -1}  \alpha^{k} g(x_k, u_k) + \alpha^\tau J(x_{\tau}) \right\},
\end{align*}   
where $\tau = \min \{ k \geq 0 \mid u^o_k = 0 \}$ with $\tau = \infty$ if this set is empty, and the expectation is with respect to the probability measure induced by $\pi^{o*}$ and the initial distribution $p^o$ of $(x_0, u_0)$, given by $p^o(\cdot) = q^o(\cdot \mid (x, u), 1)$.
\end{cor}  

\begin{proof}
Since $Q_{\k, J} = F_\k (f^*; J)$ (Lemma~\ref{lma-stopmodel1}), using the definition of $F_\k(\cdot\,; J)$ and Eq.~(\ref{eq-stop6c}), we have that for all $(x,u) \in \Gamma$, 
\begin{equation} \label{eq-prf-expQ}
  Q_{\k, J}(x,u) = g(x, u) + \alpha \int_{S \setminus B} J(x') \, q(dx' \!\mid x, u)  +  \alpha \int _B\int_C V^*(x',u') \,  \mu(du' \!\mid x') \, q(dx' \!\mid x, u),
\end{equation}
which together with (\ref{eq-stop5a}) implies the first expression for  $Q_{\k, J}(x,u)$ in the corollary. 
From this expression and Eq.~(\ref{eq-expV}) for policy $\pi^{o*}$, we obtain the second expression for $Q_{\k, J}$ in the corollary. 
\end{proof}

\subsection{A Useful Linear Program for Case (P)} \label{appsec-optstop-lp}

As Cor.\ \ref{cor-stopmodel3} shows, we can obtain $Q_{\k, J}$ from the optimal cost function $V^*$ of the optimal stopping problem. 
For case (D) (resp.\ case (N)), the function $V^*$ is the maximal solution to $V \leq T_o (V)$ among the set of bounded lower semi-analytic functions (resp.\ the set of nonpositive lower semi-analytic functions) \cite[Props.\ 9.10, 9.15]{bs}. The inequality $V \leq T_o (V)$ can be expressed as a system of linear inequalities, so under suitable conditions, $V^*$ can be obtained by solving a linear program. 
(See \cite[Chap.\ 6]{HL96} for standard linear programming formulations for DP problems with infinite state space.)

In case (P), however, $V^*$ is the minimal nonnegative lower semi-analytic solution to $V \geq T_o (V)$ \cite[Prop.\ 9.10(P)]{bs}, and this in general does not admit a linear programming formulation. 
We consider below a linear program with linear constraints based on the inequality $V \leq T_o (V)$ instead. 
While it does not yield $V^*$ in general, under an assumption to be given shortly, we can use it to obtain an upper bound on $V^*$ (in an almost-everywhere sense) and then an upper bound on $Q_{\k, J}$ (see Lemma~\ref{lma-lp}).
This bound on $Q_{\k, J}$ can be used in a  mixed value and policy iteration algorithm given in Section~\ref{secP5.2}, which is convergent under certain initial conditions for case (P), as shown by Theorem~\ref{thm-P-mixed2}.

Let $\Gamma_B = (B \times C) \cap \Gamma$ as earlier. Let $\U$ denote the universal $\sigma$-algebra on $S \times C$.

\mysmallskip
\begin{assumption} \label{assm-Plp}
{\rm (P)} \ 
There exists a $\sigma$-finite measure $\rho$ on $(S \times C, \,\U)$ such that\moveup\moveup
\begin{itemize}
\item[(i)] $\int_{\Gamma_B} J(x) \rho\big(d(x,u)\big) < \infty$; and\moveup
\item[(ii)] 
for each $(x,u) \in \Gamma_B$, the measure $\rho_{x,u}$ on $(S \times C, \, \U)$ given by 
$$\rho_{x,u}(D) = \int_B \int_C \b1_D(x', u') \, \mu(d u' \!\mid x') \, q(d x' \!\mid x, u), \qquad D \in \U,$$ 
is absolutely continuous with respect to $\rho$ (i.e., $\rho(D) = 0 \Rightarrow \rho_{x,u}(D) = 0$).\moveup
\end{itemize}
\end{assumption}
%\mysmallskip

Suppose Assumption~\ref{assm-Plp} holds (which is the case if $S$, $C$ are countable and $J$ is finite on $B$; see Remark~\ref{rmk-lp}).
Let $A_+(\Gamma_B)$ denote the set of nonnegative, lower semi-analytic functions on $\Gamma_B$.
Let $\Gamma_{B,\rho} \subset \Gamma_B$ be such that $\rho(\Gamma_B \setminus \Gamma_{B,\rho}) = 0$.
We consider a linear program on the space $A_+(\Gamma_B)$:
\begin{align}
   \mathop{\text{Maximize}}_{V \in A_+(\Gamma_B)} \ \ &   \int_{\Gamma_B} V(x, u) \, \rho\big(d(x,u)\big) \label{eq-lp} \\ 
    \text{Subject to:} \ \ &     V(x, u) \leq J(x),  \quad \forall \, (x, u) \in \Gamma_{B,\rho}, \notag \\
   &   V(x, u) \leq g(x,u) + \int_{S \setminus B} J(x') \, q(d x' \!\mid x, u) \notag \\
    & \qquad \qquad \ \ + \int_B \int_C V(x',u') \, \mu(d u' \!\mid x') \, q(d x' \!\mid x, u), \quad  \forall \, (x, u) \in \Gamma_{B,\rho}. \notag
\end{align}
As can be seen from the expression (\ref{eq-ToV})-(\ref{eq-ToVb}) for the operator $T_o$, this linear program corresponds to the following maximization problem:
\begin{align*}
   \mathop{\text{Maximize}}_{V \in A_+(\Gamma_B)} \ \ &   \int_{\Gamma_B} V(x, u) \, \rho\big(d(x,u)\big) \\ 
    \text{Subject to:} \ \ & V(x,u) \leq T_o(V^e)(x,u),  \quad \rho\text{-almost every} \ (x,u) \in \Gamma_B,
\end{align*}    
where $V^e$ is any extension of $V$ on $S^o$ satisfying Eq.~(\ref{eq-setV}) (i.e., $V^e(\infty) = 0$, $V^e(x,u) = J(x)$ for $(x,u) \in (S \setminus B) \times C$). 

Corresponding to any optimal solution $\bar V$ of (\ref{eq-lp}),
we define $\bar Q \in A_+(\Gamma)$ by the expression on the right-hand side of the second constraint in (\ref{eq-lp}), with $\bar V$ in place of $V$ and for all $(x,u)$ in $\Gamma$:
\begin{align}
\bar Q(x,u) & = g(x,u) + \int_{S \setminus B} J(x') \, q(d x' \!\mid x, u)  
    + \int_B \int_C \bar V(x',u') \, \mu(d u' \!\mid x') \, q(d x' \!\mid x, u). \label{eq-lpsol-Q}
\end{align}
The next lemma shows that $\bar Q$ satisfies a property needed for the convergence analysis of the mixed value and policy iteration algorithm (\ref{alg-1av})-(\ref{alg-1bv}) discussed in Section~\ref{secP5.2}.

\mysmallskip
\begin{lem} \label{lma-lp}
{\rm (P)} \ Let Assumption~\ref{assm-Plp} hold. 
Then an optimal solution $\bar V$ of the linear program (\ref{eq-lp}) exists, and 
the function $\bar Q \in A_+(\Gamma)$ given by Eq.~(\ref{eq-lpsol-Q}) satisfies
$$\bar Q \leq F_{\k}\big(\bar Q\,; J\big), \qquad \bar Q \geq Q_{\k, J}.$$
\end{lem}

\begin{proof}
Since $V^*=T_o(V^*)$, the optimal cost function $V^*$ restricted to $\Gamma_B$ is a feasible solution of (\ref{eq-lp}), so the feasible set of (\ref{eq-lp}) is nonempty. By Assumption~\ref{assm-Plp}(i), the optimal objective value $v^*$ of (\ref{eq-lp}) is finite. 
Let $\bar V_n, n \geq 1$, be a sequence of feasible solutions with their objective values approaching $v^*$. Then the function resulting from taking pointwise supremum, $\sup_{n} \bar V_n$, lies in $A_+(\Gamma_B)$ \cite[Lemma 7.30(2)]{bs}, satisfies the constraints of (\ref{eq-lp}), 
and achieves the optimal value $v^*$. It is hence an optimal solution of (\ref{eq-lp}). This shows that an optimal solution $\bar V$ of (\ref{eq-lp}) exists.

The function $\max \{ V^*, \bar V \}$ on $\Gamma_B$ is then an optimal solution of (\ref{eq-lp}) as well.
This implies that
\begin{equation}
  V^*(x, u) \leq \bar V(x, u) \qquad \text{for $\rho$-almost every} \ (x, u) \in \Gamma_B, \label{eq-ae-bdV}
\end{equation} 
for otherwise, by Assumption~\ref{assm-Plp}(i) we would have 
$$ \infty >  \int_{\Gamma_B} \max \big\{ V^*(x,u), \bar V(x,u) \big\} \, \rho\big(d(x,u)\big) > \int_{\Gamma_B} \bar V(x,u) \, \rho\big(d(x,u)\big),$$
a contradiction to the optimality of $\bar V$. 
We now show $\bar Q \geq Q_{\k, J}$.
By Eq.~(\ref{eq-prf-expQ}), for all $(x, u) \in \Gamma$, 
$Q_{\k,J}(x,u)$ equals the right-hand side of Eq.~(\ref{eq-lpsol-Q}) with $V^*$ in place of $\bar V$.
This, together with Assumption~\ref{assm-Plp}(ii) and the relation (\ref{eq-ae-bdV}), implies $Q_{\k, J} \leq \bar Q$.
To show $\bar Q \leq F_{\k}\big(\bar Q\,; J\big)$, notice that by the feasibility of $\bar V$ for (\ref{eq-lp}) and the definition of $\bar Q$, 
$$\bar V(x, u) \leq \min \left\{ J(x) \, , \, \bar Q(x,u)\right\}, \qquad \forall \, (x,u) \in \Gamma_{B,\rho}.$$ 
We use this relation to upper-bound $\bar V$ $\rho$-almost everywhere on $\Gamma_B$, in the integral on the right-hand side of (\ref{eq-lpsol-Q}), which defines $\bar Q$.
Using also Assumption~\ref{assm-Plp}(ii), we then obtain that for all $(x, u) \in \Gamma$,
\begin{align*}
   \bar Q(x,u) & \leq g(x,u) + \int_{S \setminus B} J(x') \, q(d x' \!\mid x, u) + 
    \int_B \int_C \min \left\{ J(x) \,, \, \bar Q(x',u') \right\} \, \mu(d u' \!\mid x') \, q(d x' \!\mid x, u), 
\end{align*}
which is the inequality $\bar Q \leq F_{\k}\big(\bar Q\,; J\big)$. This completes the proof.
\end{proof}
%\smallskip

\begin{rem} \label{rmk-lp}
Assumption~\ref{assm-Plp} holds in particular when the state and control spaces $S$ and $C$ are countable sets and the function $J$ is finite on $B$. Without loss of generality, suppose $S=C=\{1,2,\ldots\}$. 
Denote by $\rho(x,u)$ the mass assigned to a point $(x,u) \in S \times C$ by the measure $\rho$ in Assumption~\ref{assm-Plp}.
Then Assumption~\ref{assm-Plp} is satisfied by letting $\rho(x, u) = 2^{-(x+u)}/(J(x)+ 1)$ if $(x, u) \in \Gamma_B$, and $\rho(x, u) = 0$ otherwise, for instance.

If $\mu$ is a nonrandomized policy, we may let $\rho(x,\mu(x)) = 2^{-x}/(J(x)+ 1)$ if $x \in B$ and let $\rho(x,u) = 0$ for all the other $(x,u)$. Then, with $\Gamma_{B,\rho} =\{ (x, \mu(x) ) \!\mid x \in B \}$, the linear program (\ref{eq-lp}) involves only the variables $V(x, \mu(x)), x \in B$, and with the change of variable $W(x) = V(x, \mu(x))$, it becomes:
\begin{align*}
   \mathop{\text{Maximize}}_{W \geq 0} \ \ &   \sum_{x \in B} W(x) \, \rho \big(x, \mu(x) \big)  \\ 
    \text{Subject to:} \ \ &     W(x) \leq J(x),  \quad \forall \, x \in B, \notag \\
   &   W(x) \leq g \big(x,\mu(x) \big) + \sum_{x' \in S \setminus B} J(x') \, q \big( x' \!\mid x, \mu(x) \big) 
    + \sum_{x' \in B}  W(x') \, q \big( x' \!\mid x, \mu(x) \big), \quad \forall \, x \in B. 
\end{align*}
Although $S$ is countable, if $B$ is a finite set, this is a finite-dimensional linear program.\qed
\end{rem}
%\mysmallskip

\appsection{Proof of $Q_{\k,J^*} = Q^*$ for Nonnegative Case (P)}    \label{appsec2}
\markboth{\rm Appendix \thesection. Proof of $Q_{\k,J^*} = Q^*$ for Case (P)}{\rm Appendix \thesection. Proof of $Q_{\k,J^*} = Q^*$ for Case (P)}
    
In this appendix we prove for the nonnegative case (P) that for any $\k \in \K$, the function $Q_{\k,J^*} = \lim_{k \to \infty} F^k_\k(\0\,; J^*)$ 
is $Q^*$ given in Eq.~(\ref{eq-Q}). 
This establishes Prop.~\ref{prp-propertyF}(d) for (P), which is also used in the lower bound part of Lemma~\ref{lma-bound-NP} for (P).

\mysmallskip
\begin{prop} \label{lma-P-optstop}
{\rm (P)} \ Let $\k = (\mu, B) \in \K$. We have $Q_{\k,J^*} = Q^*.$ 
\end{prop}
\mysmallskip

Since $Q^* \geq 0$ and  $F_\k(Q^*; J^*) = Q^*$ (Prop.~\ref{prp-propertyF}(a)), we have by the monotonicity of $F_\k$ (cf.~Eq.~(\ref{eq-Fmono})),
\begin{equation} \label{eq-prfP-Qineq0}
Q_{\k,J^*} = \lim_{n\to \infty} F^n_\k (\0\,; J^*) \leq Q^*.
\end{equation}
Thus to prove Prop.~\ref{lma-P-optstop}, we need to show $Q_{\k,J^*} \geq Q^*$. We will prove this by showing that for each $(x, u) \in \Gamma$ and any $\epsilon > 0$,
\begin{equation} \label{eq-prfP-Qineq}
 Q_{\k,J^*}(x,u) \geq Q^*(x,u) - \epsilon. 
\end{equation} 
In the proof we will use the correspondence between the optimal stopping problem associated with $\k =(\mu, B)$ and $J^*$, as defined in Appendix~\ref{secA.1}, and a controller for the original problem. 

We need some notations and an expression of $Q_{\k, J^*}$ to be used in the proof.
Fix $(\bar x, \bar u) \in \Gamma$. 
For the optimal stopping problem associated with $\k =(\mu, B)$ and $J^*$, 
by \cite[Cor.~9.17.1]{bs}, there exists an optimal stationary nonrandomized (universally measurable) policy $\mu^o : S^o=(S \times C) \cup \{\infty\} \to \{0,1\}$. 
Let the optimal stopping problem start from time $1$, and consider the stochastic process $(z_1, u^o_1), (z_2, u^o_2), \ldots$, where $z_k \in S^o$ and $u^o_k \in \{0,1\}$, induced by $\mu^o$ and the initial distribution of $z_1$ given by $q^o(\cdot \!\mid (\bar x, \bar u), 1)$ (cf.\ Eq.~(\ref{eq-stopmodel-q})).
For each $k \geq 1$, define $(x_k,v_k) = z_k$ if $z_k \in S \times C$, and define $(x_k,v_k)$ to be some fixed point in $S \times C$ if $z_k = \infty$ (the absorbing state). Here for clarity, we are using $v_k$ instead of $u_k$ to denote the component of $z_k$ in $C$, since we will use $u_k$ later for the controls applied in the original problem.
By Cor.~\ref{cor-stopmodel3} we have
\begin{equation} \label{eq-prfP-Q}
   Q_{\k, J^*}(\bar x, \bar u) = g(\bar x, \bar u) + \E^{\mu^o} \left\{ \sum_{k=1}^{\tau-1}  g(x_k, v_k) + J^*(x_{\tau}) \right\}, 
\end{equation}   
where $\tau$ is the time the process is stopped: $\tau = \min \big\{ k \geq 1 \mid \mu^o(x_k, v_k) = 0 \big\}$ ($\infty$ if the set is empty), and $\E^{\mu^o}$ denotes expectation with respect to the probability measure induced by $\mu^o$ and the initial distribution of $(x_1, v_1)$, given by $q^o(\cdot \!\mid (\bar x, \bar u), 1)$. 

To simplify notation, let 
$$ D = \{ (x,v) \in S \times C \mid \mu^o(x, v) = 0 \}, 
\qquad D_x = \{ v \in C \mid (x, v) \in D \}, \quad x \in S. $$
($D$ is the subset of $S \times C$ 
on which $\mu^o$ stops the process.)
Since $\mu^o$ is universally measurable, $D$ and hence $D_x$, $x \in S$, are universally measurable sets \cite[Lemma 7.29]{bs}.
Note that expressed in terms of these sets,
$\tau = \min \big\{ k \geq 1 \mid v_k \in D_{x_k} \big\}$
and
\begin{align}  \label{eq-appB-0}
  \tau & = m \qquad  \Longleftrightarrow  \qquad  v_1 \not\in D_{x_1}, \  \cdots, \ v_{m-1} \not\in D_{x_{m-1}}, \ v_m \in D_{x_m}, \\
  \tau & > m \qquad  \Longleftrightarrow  \qquad  v_1 \not\in D_{x_1}, \  \cdots, \ v_{m-1} \not\in D_{x_{m-1}}, \ v_m \not\in D_{x_m}. \label{eq-tau-2}
\end{align}  

We consider also the probability measure on the space of $(x_1, v_1, x_2, v_2, \ldots)$ induced by the policy $\mu$ and the initial distribution $q(dx_1 \!\mid \bar x, \bar u)$ of $x_1$, and the expectation
$\E^\mu \left\{ \sum_{k=1}^{\tau-1}  g(x_k, v_k) +  J^*(x_{\tau}) \right\} $ with respect to this probability measure, 
where $\tau$ is the same as defined earlier.

\mysmallskip
\begin{lem} \label{lma-appB-1}
We have 
$$\E^{\mu^o} \left\{ \sum_{k=1}^{\tau-1}  g(x_k, v_k) +  J^*(x_{\tau}) \right\} = \E^{\mu} \left\{ \sum_{k=1}^{\tau-1}  g(x_k, v_k) + J^*(x_{\tau}) \right\}.$$
\end{lem} 
\smallskip

We first prove Prop.~\ref{lma-P-optstop}, using Lemma~\ref{lma-appB-1}. We will then give the proof of this lemma.

\begin{proof}[Proof of Prop.~\ref{lma-P-optstop}]
Fix $(\bar x, \bar u) \in \Gamma$ and let $\epsilon > 0$. 
Let $\pi^\epsilon = (\pi^\epsilon_0, \pi^\epsilon_1, \ldots) \in \Pi$ be an $\epsilon$-optimal Markov policy for the original control problem; such a policy exists by \cite[Prop.~9.19]{bs}.
We use the policies $\mu^o$, $\pi^\epsilon$, and $\mu$ (the stationary policy defining $F_\k$ and the associated optimal stopping problem) to define a controller $\hat \pi$ for the original problem, such that it applies control $\bar u$ at state $\bar x$ at the first stage, and its expected total cost for state $\bar x$ is no greater than $Q_{\k, J^*}(\bar x, \bar u)  + \epsilon$. From this the desired inequality (\ref{eq-prfP-Qineq}) for establishing the proposition will be shown to follow. 

Roughly speaking, the controller $\hat \pi$ follows the policy $\mu$ before it switches to following policy $\pi^\epsilon$. 
To decide when to make the switch, it generates at each time $k \geq 1$, an auxiliary variable $v_k \in C$ to ``simulate'' a control that $\mu$ might apply at the current state and ``query'' the optimal-stopping policy $\mu^o$ about whether that control suggested by $\mu$ should be followed. 
The history at time $k \geq 1$ for the controller is
$$ (x_0, u_0, x_1, v_1, u_1, \ldots, x_k, v_k ) \in (S \times C) \times (S \times C^2)^{k-1} \times (S \times C),$$
including the auxiliary variables $v_j,  1\leq j\leq k$, as well as the past states $x_j, j \leq k$, and past controls $u_j, j < k$.
The controller is denoted $\hat \pi = (\hat \mu_0, \hat \mu_1, \ldots )$, where each $\hat \mu_k$ is a universally measurable stochastic kernel on $C$ given the respective space of histories. We now define 
$\hat \mu_k, k \geq 0$.

For $k=0$, let $\hat \mu_0$ be a universally measurable stochastic kernel on $C$ given $S$, such that $\hat \mu_0$ satisfies the control constraint $U$ and
$\hat \mu_0(du_0 \!\mid \bar x) = \delta_{\bar u}.$
For each $k \geq 1$, the auxiliary variable $v_k$ is generated according to the stochastic kernel $\mu$ 
given the state $x_k$:
\begin{equation} \label{eq-Pprf-stop1}
         \mu\big(dv_k \!\mid x_k\big).
\end{equation}
The stochastic kernels $\hat \mu_k$, $k \geq 1$, are defined as follows. For each $k \geq 1$, define a universally measurable function $\tau_k : (S \times C)^{k} \to \{0, 1, 2, \ldots \}$ by 
\begin{equation} \label{eq-Pprf-stop2}
  \tau_k(x_1, v_1, x_2, v_2, \ldots, x_k, v_k) = \begin{cases}
\min \big\{ m  \mid v_m \in D_{x_m},  \, 1 \leq m \leq k \big\} & \text{if such $m$ exists}, \\
 0  & \text{otherwise}.
 \end{cases} 
\end{equation} 
Let $\hat \mu_k$ be a universally measurable stochastic kernel on $C$ given $(S \times C) \times (S \times C^2)^{k-1} \times (S \times C)$, given by
\begin{equation} \label{eq-Pprf-stop3}
\hat \mu_k (du_k \!\mid x_0, u_0, x_1, v_1, u_1, \ldots, x_k, v_k) = \begin{cases} 
    \delta_{v_k} & \text{if} \  \, \tau_k(x_1, v_1, \ldots, x_k, v_k) = 0, \\
     \pi^\epsilon_{k-m}( du_k \! \mid x_{k}) & \text{if} \ \, \tau_k(x_1, v_1, \ldots, x_k, v_k) = m \geq 1.
           \end{cases}
\end{equation}
(Thus, $\hat \pi$ ``copies'' the control $v_k$ if it has not yet switched to applying policy $\pi^\epsilon$, and the switch happens the first time $v_m \in D_{x_m}$.)

The controller $\hat \pi = (\hat \mu_0, \hat \mu_1, \ldots )$ induces a probability measure on the space $(S \times C) \times (S \times C^2)^\infty$ of 
$(x_0, u_0, x_1, v_1, u_1, x_2, v_2, u_2, \ldots )$ (with the universal $\sigma$-algebra). 
With respect to this probability, the expected total cost of $\hat \pi$ for state $\bar x$ is 
$$\hat J_{\hat \pi}(\bar x)  = g(\bar x, \bar u) + \E^{\hat \pi}  \left\{ \sum_{k=1}^{ \infty}  g(x_k, u_k) \right\}.$$
Let $\tau = \min \big\{ k \geq 1  \mid v_k \in D_{x_k} \big\}$ with $\tau = \infty$ if the set in the definition is empty.
Using the definition of conditional expectation and a formula for conditional expectation given the sub-$\sigma$-algebra associated with the stopping time $\tau$ \cite[Prop.~II-1-3]{Nev75}, it follows that
\begin{align} \label{eq-prfP-stop4}
 \hat J_{\hat \pi}(\bar x) & =  g(\bar x, \bar u) + \E^{\hat \pi} \left\{ \sum_{k=1}^{ \tau -1}  g(x_k, u_k) +  J_{\pi^\epsilon}(x_{ \tau}) \right\}  
    =  g(\bar x, \bar u) + \E^{\hat \pi} \left\{ \sum_{k=1}^{\tau -1}  g(x_k, v_k) +  J_{\pi^\epsilon}(x_\tau) \right\},
\end{align} 
where the second equality follows from the fact $u_k = v_k$ for $k < \tau$ (cf.\ Eqs.~(\ref{eq-Pprf-stop2}), (\ref{eq-Pprf-stop3})).
We have
$$ J_{\pi^\epsilon}(x) \leq J^*(x) + \epsilon, \qquad \forall \, x \in S,$$
since $\pi^\epsilon$ is an $\epsilon$-optimal policy of the original problem and $J^* \geq 0$.
Then by Eq.~(\ref{eq-prfP-stop4}),
$$ \hat J_{\hat \pi}(\bar x) \leq g(\bar x, \bar u) + \E^{\hat \pi} \left\{ \sum_{k=1}^{\tau -1}  g(x_k, v_k) +  J^*(x_\tau) \right\} + \epsilon. $$

On the other hand, based on the definition of $\hat \pi$, $\{v_k\}$ are generated according to $\mu$ (cf.~Eq.~(\ref{eq-Pprf-stop1})) and they are the controls applied before time $\tau$ (cf.\ Eqs.~(\ref{eq-Pprf-stop2}), (\ref{eq-Pprf-stop3})), 
so we have
\begin{align*}
 \E^{\hat \pi} \left\{ \sum_{k=1}^{\tau -1} g(x_k, v_k) +  J^*(x_\tau) \right\} & = \E^{\mu} \left\{ \sum_{k=1}^{\tau -1}  g(x_k, v_k) +  J^*(x_\tau) \right\} 
 = \E^{\mu^o} \left\{ \sum_{k=1}^{\tau -1}  g(x_k, v_k) +  J^*(x_\tau) \right\},
\end{align*} 
where the second equality follows from Lemma~\ref{lma-appB-1}.
Combining the preceding two relations with Eq.~(\ref{eq-prfP-Q}), we obtain
\begin{equation} \label{eq-prfP-stop5}
 \hat J_{\hat \pi}(\bar x) \leq g(\bar x, \bar u) + \E^{\mu^o} \left\{ \sum_{k=1}^{\tau -1}  g(x_k, v_k) +  J^*(x_\tau)  \right\} + \epsilon = Q_{\k,J^*}(\bar x, \bar u) + \epsilon.
 \end{equation}
 
Although the controller $\hat \pi$ uses the additional auxiliary variables $\{v_k\}$ for control, it does not have advantages over the set of policies in $\Pi$, in the sense that we can construct a semi-Markov randomized policy in $\Pi$ such that it applies control $\bar u$ at the first stage if $\bar x$ is the initial state, and it has the same expected total cost $\hat J_{\hat \pi}(\bar x)$ for the state $\bar x$. (Such a construction is similar to that used to prove Props.\ 8.1, 9.1 in \cite{bs}.)
This means that
$\hat J_{\hat \pi}(\bar x) \geq Q^*(\bar x, \bar u)$ (cf.\ Eq.~(\ref{eq-Qalt})),  so by Eq.~(\ref{eq-prfP-stop5}),
$$ Q^*(\bar x, \bar u) \leq Q_{\k,J^*}(\bar x, \bar u) + \epsilon.$$
Since $\epsilon$ is arbitrary, we have $Q_{\k, J^*}(\bar x, \bar u) \geq Q^*(\bar x, \bar u)$. Then by Eq.~(\ref{eq-prfP-Qineq0}), we have $Q_{\k, J^*} = Q^*$ .
\end{proof}   
%\mysmallskip

We now establish Lemma~\ref{lma-appB-1}.

%\mysmallskip
\begin{proof}[Proof of Lemma~\ref{lma-appB-1}]
To prove
$\E^{\mu^o} \left\{ \sum_{k=1}^{\tau-1}  g(x_k, v_k) +  J^*(x_{\tau}) \right\} = \E^{\mu} \left\{ \sum_{k=1}^{\tau-1}  g(x_k, v_k) + J^*(x_{\tau}) \right\},$
we first introduce some functions to rewrite the two expectations.
In view of (\ref{eq-appB-0}) (i.e., $\tau = m \Leftrightarrow  v_1 \not\in D_{x_1},   \cdots,  v_{m-1} \not\in D_{x_{m-1}},  v_m \in D_{x_m}$),
we define for each $m \geq 1$, $\phi_m: (S \times C)^m \to [0,\infty]$ by
\begin{align*}
  \phi_m(x_1, v_1, \ldots, x_m, v_m) : = & \left( \prod_{i =1}^{m-1} \b1_{C \setminus D_{x_i}}(v_i) \right) 
 \cdot \b1_{D_{x_m}}(v_m) \cdot \left( \sum_{k=1}^{m-1}  g(x_k, v_k) +  J^*(x_{m}) \right) \\
 = & \  \b1_{\{\tau = m\}}(x_1, v_1, \ldots) \cdot \left(   \sum_{k=1}^{\tau-1}  g(x_k, v_k) +  J^*(x_{\tau}) \right).
\end{align*}
We define for $m = \infty$, $\phi_\infty: (S \times C)^\infty \to [0,\infty]$ by
\begin{align*}
\phi_\infty(x_1, v_1, x_2, v_2, \ldots ) : = &  \left( \prod_{i =1}^\infty \b1_{C \setminus D_{x_i}}(v_i) \right) \cdot  \sum_{k=1}^{\infty} g(x_k, v_k) 
=  \b1_{\{\tau = \infty \}}(x_1, v_1, \ldots) \cdot   \sum_{k=1}^{\infty}  g(x_k, v_k).
\end{align*}
Since $g \geq 0, J^* \geq 0$, we may write
\begin{align}
 \E^{\mu^o} \left\{ \sum_{k=1}^{\tau-1}  g(x_k, v_k) +  J^*(x_{\tau}) \right\} 
 & = \sum_{m=1}^\infty \E^{\mu^o} \left \{ \b1_{\{\tau = m\}}(x_1, v_1, \ldots) \cdot \left( \sum_{k=1}^{m-1}  g(x_k, v_k) +  J^*(x_{m}) \right) \right\}  \notag \\
 & \ \ \ \ + \E^{\mu^o} \left \{ \b1_{\{\tau = \infty \}}(x_1, v_1, \ldots) \cdot  \sum_{k=1}^{\infty}  g(x_k, v_k)  \right\} \notag \\
& = \sum_{m \, \in \{ 1, 2, \ldots \} \cup \{ \infty \}} \E^{\mu^o} \big\{ \phi_m(x_1, v_1, \ldots, x_m, v_m)  \big\},  \label{eq-appB-1a} 
\end{align}
and similarly,
\begin{align}
\E^{\mu} \left\{ \sum_{k=1}^{\tau-1}  g(x_k, v_k) +  J^*(x_{\tau}) \right\} 
& = \sum_{m \, \in \{ 1, 2, \ldots \} \cup \{ \infty \}} \E^{\mu} \big\{ \phi_m(x_1, v_1, \ldots, x_m, v_m)  \big\}. \label{eq-appB-1b} 
\end{align}

To prove that (\ref{eq-appB-1a}) and (\ref{eq-appB-1b}) are equal, we will proceed in four steps.

\vspace*{0.07in}
\noindent (i) First, we show that for each integer $m \geq 1$, 
\begin{equation} \label{eq-appB-1}
 \E^{\mu^o} \big\{ \phi_m(x_1, v_1, \ldots, x_m, v_m)  \big\} = \E^{\mu} \big\{ \phi_m(x_1, v_1, \ldots, x_m, v_m)  \big\}.
\end{equation}
Note that $\phi_m(x_1, v_1, \ldots, x_m, v_m)=0$ on $\{\tau \not=m\}$, 
and  $\tau \not=m$ if $x_i \not\in B$ for some $i < m$ (since in the optimal stopping problem, the only control that $\mu^o$ can take for states in $(S \setminus B) \times C$ is to stop, $x_i \not\in B$ implies $\tau \leq i$).
Using these facts together with the definition of $\E^{\mu}$, we have that
$\E^{\mu} \!\big\{ \phi_m(x_1, v_1, \ldots, x_m, v_m) \big\}$ is equal to
\begin{align}
 &  \int_B \int_C \ldots  \int_B \int_C \left[ \int_S \int_C  \phi_m(x_1, v_1, \ldots, x_m, v_m) \, \mu(dv_m \!\mid x_m) \, q(d x_m \!\mid x_{m-1}, v_{m-1}) \right]  \cdot \notag \\
   & \qquad  \qquad \qquad \quad  \mu(dv_{m-1} \!\mid x_{m-1}) \, q(d x_{m-1} \!\mid x_{m-2}, v_{m-2}) \, \cdots  \mu(dv_1 \!\mid x_1) \, q(d x_1 \!\mid \bar x, \bar u). \label{eq-appBprf-int0}
\end{align}
Using the same facts just mentioned, and using also the definition of the optimal stopping problem (Appendix~\ref{secA.1}), 
we have that 
$\E^{\mu^o} \!\big\{ \phi_m(x_1, v_1, \ldots, x_m, v_m) \big\}$ is equal to
\begin{align*}
 &  \int_B \int_C \ldots  \int_B \int_C \left[ \int_S \int_C  \phi_m(x_1, v_1, \ldots, x_m, v_m) \, \tilde \mu(dv_m \!\mid x_m) \, q(d x_m \!\mid x_{m-1}, v_{m-1}) \right]  \cdot \\
   & \qquad  \qquad \qquad \quad  \tilde \mu(dv_{m-1} \!\mid x_{m-1}) \, q(d x_{m-1} \!\mid x_{m-2}, v_{m-2}) \, \cdots  \tilde \mu(dv_1 \!\mid x_1) \, q(d x_1 \!\mid \bar x, \bar u).
\end{align*}
Since $\tilde \mu(\cdot \!\mid x)= \mu(\cdot \!\mid x)$ for $x \in B$ by the definition of the optimal stopping problem, the above integral in turn equals
\begin{align}
 &  \int_B \int_C \ldots  \int_B \int_C \left[ \int_S \int_C  \phi_m(x_1, v_1, \ldots, x_m, v_m) \, \tilde \mu(dv_m \!\mid x_m) \, q(d x_m \!\mid x_{m-1}, v_{m-1}) \right]  \cdot \notag \\
   & \qquad  \qquad \qquad \quad \mu(dv_{m-1} \!\mid x_{m-1}) \, q(d x_{m-1} \!\mid x_{m-2}, v_{m-2}) \, \cdots \mu(dv_1 \!\mid x_1) \, q(d x_1 \!\mid \bar x, \bar u). \label{eq-appBprf-int}
\end{align}
Consider now the inner-most integral in (\ref{eq-appBprf-int}). 
If $x_m \in S \setminus B$, then in view of the control constraint $U^o$ of the optimal stopping problem (cf.\ Appendix~\ref{secA.1}), 
we have $D_{x_m} = C$. Hence $\b1_{D_{x_m}}(v_m)=1$ for all $v_m \in C$, so
$$  \phi_m(x_1, v_1, \ldots, x_m, v_m) 
= \left( \prod_{i =1}^{m-1} \b1_{C \setminus D_{x_i}}(v_i) \right) 
 \cdot \left( \sum_{k=1}^{m-1}  g(x_k, v_k) +  J^*(x_{m}) \right)$$
does not depend on $v_m$. Consequently,
\begin{align*}
 \int_C  \phi_m(x_1, v_1, \ldots, x_m, v_m) \, \tilde \mu(dv_m \!\mid x_m) 
=  \int_C  \phi_m(x_1, v_1, \ldots, x_m, v_m) \,  \mu(dv_m \!\mid x_m), \quad x_m \in S \setminus B. 
\end{align*}
If $x_m \in B$, then since $\tilde \mu(\cdot \!\mid x)= \mu(\cdot \!\mid x)$ for $x \in B$, we have 
\begin{align*}
 \int_C  \phi_m(x_1, v_1, \ldots, x_m, v_m) \, \tilde \mu(dv_m \!\mid x_m) 
=  \int_C  \phi_m(x_1, v_1, \ldots, x_m, v_m) \, \mu(dv_m \!\mid x_m), \quad x_m \in B.
\end{align*}
The preceding two equalities together imply that the value of the integral (\ref{eq-appBprf-int}) remains unchanged if we replace $\tilde \mu(dv_m \! \mid x_m)$ in the inner-most integral in (\ref{eq-appBprf-int}) by $\mu(d v_m \!\mid x_m)$. Hence the integral (\ref{eq-appBprf-int}) is equal to the integral (\ref{eq-appBprf-int0}), and this proves the desired equality (\ref{eq-appB-1}) for $m \geq 1$.

\vspace*{0.07in}
\noindent (ii) By arguments similar to the ones in the preceding proof, we have that
for all $m \geq 1$ and $n \geq m$,
\begin{equation} \label{eq-appB-4}
 \E^{\mu^o} \left\{ \left( \prod_{i =1}^n \b1_{C \setminus D_{x_i}}(v_i) \right) \cdot \sum_{k=1}^{m}  g(x_k, v_k) \right\} = \E^{\mu} \left\{ \left( \prod_{i =1}^n \b1_{C \setminus D_{x_i}}(v_i) \right) \cdot \sum_{k=1}^{m}  g(x_k, v_k) \right\}.
\end{equation} 
In particular, observing that $\prod_{i =1}^n \b1_{C \setminus D_{x_i}}(v_i) = 0$ if $x_i \not\in B$ for some $i \leq n$,
an analysis similar to the first half of the proof in (i) then shows that
both sides of (\ref{eq-appB-4}) are equal to
\begin{align*}
    \int_B \int_C \ldots  \int_B \int_C &  \left( \prod_{i =1}^n \b1_{C \setminus D_{x_i}}(v_i) \right) \cdot \left( \sum_{k=1}^{m}  g(x_k, v_k)  \right) \, \mu(dv_n \!\mid x_n) \, q(d x_n \!\mid x_{n-1}, v_{n-1})  \cdot \\
    & \cdots \mu(dv_1 \!\mid x_1) \, q(d x_1 \!\mid \bar x, \bar u).
\end{align*}     
We will need (\ref{eq-appB-4}) shortly in the proof.

\vspace*{0.07in}
\noindent (iii) 
Let us now consider the two terms corresponding to $m = \infty$ in Eqs.~(\ref{eq-appB-1a}) and (\ref{eq-appB-1b}). 
We examine when they are equal, i.e., when
\begin{equation}
   \E^{\mu^o} \!\big\{ \phi_\infty(x_1, v_1,x_2, v_2, \ldots ) \big\} = \E^{\mu}  \!\big\{ \phi_\infty(x_1, v_1,x_2, v_2, \ldots ) \big\}. \label{eq-appB-2}
\end{equation}   
From the definition of $\phi_\infty$, we have, by the monotone convergence theorem, that as $m \to \infty$,
$$\E^{\mu^o} \left\{  \left( \prod_{i =1}^\infty \b1_{C \setminus D_{x_i}}(v_i) \right) \cdot \sum_{k=1}^{m}  g(x_k, v_k)  \right\} \  \big\uparrow \ \, \E^{\mu^o} \!\big\{ \phi_\infty(x_1, v_1,x_2, v_2, \ldots ) \big\},$$
$$ \E^{\mu} \left\{  \left( \prod_{i =1}^\infty \b1_{C \setminus D_{x_i}}(v_i) \right) \cdot \sum_{k=1}^{m}  g(x_k, v_k)  \right\} \  \big\uparrow \ \, \E^{\mu}  \!\big\{ \phi_\infty(x_1, v_1,x_2, v_2, \ldots ) \big\}.$$
Thus Eq.~(\ref{eq-appB-2}) holds if for each integer $m \geq 1$,
\begin{equation} \label{eq-appB-3}
   \E^{\mu^o} \left\{ \left( \prod_{i =1}^\infty \b1_{C \setminus D_{x_i}}(v_i) \right) \cdot \sum_{k=1}^{m}  g(x_k, v_k)   \right\} = \E^{\mu} \left\{  \left( \prod_{i =1}^\infty \b1_{C \setminus D_{x_i}}(v_i) \right) \cdot \sum_{k=1}^{m}  g(x_k, v_k)  \right\}.
\end{equation}    
Now, for each positive integer $m$, we have the following pointwise convergence of functions as $n \to \infty$:
$$   \left( \prod_{i =1}^n \b1_{C \setminus D_{x_i}}(v_i) \right) \cdot \sum_{k=1}^{m}  g(x_k, v_k)  \, \  \big\downarrow \   \left( \prod_{i =1}^\infty \b1_{C \setminus D_{x_i}}(v_i) \right) \cdot \sum_{k=1}^{m}  g(x_k, v_k).$$
We also have the equality (\ref{eq-appB-4}) for all $n \geq m$.
Hence, if for each $m$ there exists some $n \geq m$ for which the quantity in (\ref{eq-appB-4}) is less than $\infty$, then by the dominated convergence theorem \cite[p.\ 132]{Dud02}, (\ref{eq-appB-3}) holds for each $m$, and hence the desired equality (\ref{eq-appB-2}) holds, which together with (\ref{eq-appB-1}) implies that (\ref{eq-appB-1a}) and (\ref{eq-appB-1b}) are equal.

\vspace*{0.05in}
\noindent (iv) The only case left now is that for some integer $m \geq 1$ and all $n \geq m$, the quantity in (\ref{eq-appB-4}) is $\infty$. But in view of Eq.~(\ref{eq-tau-2}) (i.e., $\tau > m   \Leftrightarrow  v_1 \not\in D_{x_1},   \cdots,  v_{m-1} \not\in D_{x_{m-1}},  v_m \not\in D_{x_m}$), this would imply
$$ \E^{\mu^o} \left\{  \b1_{\{\tau > m \}}(x_1, v_1, \ldots) \sum_{ k = 1}^{\tau -1}  g(x_k, v_k) \right\} = \infty,  \qquad \E^{\mu} \left\{  \b1_{\{\tau > m\}}(x_1, v_1, \ldots) \sum_{ k = 1}^{\tau -1}  g(x_k, v_k) \right\} = \infty,$$
and hence both (\ref{eq-appB-1a}) and (\ref{eq-appB-1b}) equal $\infty$. This completes the proof.
\end{proof}

\appsection{About the Existence of Borel Measurable Policies}  \label{appsec-borelpolicy}
\markboth{\rm Appendix \thesection. Existence of Borel Measurable Policies}{\rm Appendix \thesection. Existence of Borel Measurable Policies}

In this appendix, we mention some useful facts about the existence of Borel measurable policies, to demonstrate the application range of Algorithm III given in Section~\ref{sec3.2}, which uses Borel measurable policies.

Recall that the graph of the control constraint $U$, 
$\Gamma = \big\{ (x, u) \mid x \in S, u \in U(x) \big\},$
is an analytic subset of the product of two Polish spaces (of which $S$ and $C$ are Borel subsets).
The question whether a Borel measurable nonrandomized stationary policy exists in our control problem 
is equivalently whether the set $\Gamma$ admits a Borel measurable function $f: S \to C$ whose graph lies in $\Gamma$ (i.e., $f(x) \in U(x)$ for all $x$).

Suppose $\Gamma$ is a Borel subset of $S \times C$.
It can still happen that such $f$ does not exist~\cite{Blk-borel}. 
Then there exists 
no Borel measurable stationary policy, randomized or nonrandomized.
(Because if a randomized Borel measurable stationary policy were to exist, a nonrandomized one must also exist by the selection theorem of Blackwell and Ryll-Nardzewski \cite{BlkR63}.)
Nevertheless, when $\Gamma$ is Borel, a number of selection theorems for Borel sets in the product of two Polish spaces can be applied for a Borel measurable selection in $\Gamma$, under fairly general conditions on the control constraint $U$ (see e.g.,~\cite{Sriv-borel}). 
Below are some examples. 

Let $Y$ be the Polish space of which $C$ is a Borel subset. Assume $\Gamma$ is Borel. 
Then in each of the following cases, a Borel measurable nonrandomized stationary policy exists:%\moveup\moveup
\begin{enumerate}
\item[(a)] For every $x$, $U(x)$ is a countable set (by a theorem of Lusin, \cite[Theorem 5.8.11]{Sriv-borel}).\moveup%\moveup
\item[(b)] For every $x$, $U(x)$ contains a nonempty open set in $Y$ (by theorems of Kechris and Sarbadhikari, \cite[Theorem 5.8.5]{Sriv-borel}).\moveup%\moveup
\item[(c)] For every $x$, $U(x)$ is a $\sigma$-compact set in $Y$ (by a theorem of Arsenin and Kunugui, \cite[Theorem 5.12.1]{Sriv-borel}), which is true in particular when $Y$ is $\sigma$-compact and each $U(x)$ is a countable union of open or closed sets in $Y$.\moveup%\moveup 
\item[(d)] $U$ is a Borel measurable multifunction (i.e., set-valued function) and for every $x$, $U(x)$ is a closed set in $Y$ (by the Kuratowski and Ryll-Nardzewski selection theorem, \cite[Theorem 5.2.1]{Sriv-borel}).\moveup%\moveup
\end{enumerate}
Thus for many general control constraints, 
the policy iteration-like algorithm III can be applied.

\appsection{An Illustrative Example for Value Iteration in Case (P)} \label{appsec-Pexample}
\markboth{\rm Appendix \thesection. An Illustrative Example for Value Iteration in Case (P)}{\rm Appendix \thesection. An Illustrative Examples for Value Iteration in Case (P)}

In this appendix we use an example to illustrate Theorem~\ref{thm-valite-P}(b) for the convergence of value iteration in case (P).
This example is from Strauch \cite[Example 6.2, p.\ 881]{Str-negative} and also described in Maitra and Sudderth \cite[p.\ 930]{MS92}. 
Our description below closely follows \cite{MS92}. As mentioned in Remark~\ref{rem-P2}, in this example value iteration converges in one iteration if we know the set $\{ x\in S \mid J^*(x) = 0\}$ and hence can choose an initial function $J$ to satisfy $J^* \leq J \leq c J^*$ for some constant $c$, fulfilling the convergence condition in Theorem~\ref{thm-valite-P}(b). In what follows we give a detailed description and calculation to supplement this remark.

Let $\mathcal{R}_{(0,1)}$ denote the set of rationals in $(0,1)$ with its usual ordering, and index its elements by $r_1, r_2, \ldots$ (in an arbitrary way).
Let $\big\{W_r \mid r \in \mathcal{R}_{(0,1)} \big\}$ be a collection of Borel subsets of $[0,1]$ (called a Borel sieve). 
Correspondingly, define for each $z \in [0,1]$,
$$ M_z = \big\{ r \in  \mathcal{R}_{(0,1)} \mid z \in W_r \big\},$$
and let
$$D = \big\{ z \in [0,1] \, \big| \, M_z \text{ is not well-ordered} \big\}.$$
Fix the sets $\{W_r\}$ such that the set $D$ is not Borel measurable.
Define the control problem as follows.

Let $S = \big\{ (z, r) \mid 0 \leq z \leq 1, \  0 \leq r \leq 1, \ r \text{ rational} \big\} \cup \{ t\}$. 
Let $C = \{ 1, 2, \ldots \}$ and $U(x) = C$ for every state $x \in S$.
State transitions are deterministic. The successor state $f(x,u)$ when applying control $u$ at state $x$ is given by
$$ f(t, u) = t, \qquad f\big( (z, r), u \big) = \begin{cases}  (z, r_u) & \text{if} \ r_u < r \ \text{and} \ z \in W_{r_u}, \\  t & \text{otherwise}. \end{cases}$$
The cost $\hat g(x, u, x')$ of transition from state $x$ to state $x'$ is given by
$$ \hat{g}(t, u, t) = 0, \qquad \hat g\big( (z, r), u, x' \big) = \begin{cases}  0 & \text{if} \ x' \not= t, \\ 1 & \text{otherwise}. \end{cases} $$ 
(Equivalently, the one-stage costs are given by $g(t, u) = 0$, $g\big( (z, r), u \big) = 0$ if $r_u < r$ and $z \in W_{r_u}$, and $g\big( (z, r), u \big) = 1$ otherwise.)
The optimal cost function $J^*$ takes only values $0$ or $1$, and it is not Borel measurable \cite{Str-negative}. 
As shown by \cite{Str-negative}, $J^*(t) = 0$ and  for states $(z, 1)$ where $z \in [0,1]$, 
\begin{equation} \label{eq-exC1}
    J^*(z,1) = \begin{cases}  0 & \text{if} \ z \in D, \\
1 & \text{if} \ z \in [0,1] \setminus D. \end{cases} 
\end{equation}
Value iteration starting from the constant function zero requires uncountably many iterations to converge to $J^*$, 
as shown by Maitra and Sudderth \cite[p.\ 930]{MS92}.

Denote $G = \big\{ x \in S \mid J^*(x) = 0 \big\}.$ Then $T^k(J) \to J^*$ if we let $J$ be $J(t) = 0$ and for states $x = (z, r)$,
$$ J(z, r) = \begin{cases}  0 & \text{if} \ (z, r) \in G, \\
v & \text{otherwise}, \end{cases} \qquad \quad \text{for some constant $a \geq 1$}.$$
Since $J^* \leq J \leq a J^*$, this function $J$ satisfies the condition of Theorem~\ref{thm-valite-P}(b) for the convergence of value iteration.

Indeed $T(J) = J^*$, as can be verified directly. 
Consider any $(z, r) \in S$ where $z \in [0,1]$. 
If $(z, r) \in G$, then by the definition of the control problem given above and by the relation $J^* = T(J^*)$, 
there must exists some $u \in C$ such that with $x' = f\big((z, r), u\big)$, 
$$J^*(z, r) = \hat g\big( (z, r), u, x'\big) + J^*(x') = 0,$$
which implies that $\hat g\big( (z, r), u, x'\big) = 0$ and $x' \in G$. Consequently, we have $T(J)(z, r) = 0 = J^*(z,r)$.

Suppose $(z, r) \not\in G$, i.e., $J^*(z, r) = 1$. Then, by the relation $J^* = T(J^*)$ and the binary nature of the costs,
we must have that for each $u \in C$, either (i) or (ii) can happen:\moveup\moveup
\begin{itemize}
\item[(i)] $f\big((z, r), u\big) = t$ and  $\hat g\big( (z, r), u, t\big)  = 1$, in which case
$\hat g\big( (z, r), u, t\big) + J(t) = 1$.\moveup 
\item[(ii)] $x'=f\big((z, r), u\big) \not= t$, $\hat g\big( (z, r), u, x'\big)  = 0$, and $J^*(x') = 1$ (i.e., $x' \not\in G$), in which case
$\hat g\big( (z, r), u, x'\big) + J(x') = a \geq 1.$\moveup\moveup 
\end{itemize}
Therefore, if there exists $u$ satisfying (i), then $T(J)(z, r) = 1 = J^*(z, r)$. 

Now, if $r = 0$, then only case (i) can happen, since $r_u > 0$ for all $u$.
If $r \in (0,1)$, then there exists $u$ with $r_u = r$, and this $u$ satisfies (i).
Finally, suppose $r = 1$. Then the assumption $J^*(z, r) = 1$ implies that $z \not\in D$ (cf.\ Eq.~(\ref{eq-exC1})). By the definition of $D$, this means that $M_z$ is well-ordered and therefore has a smallest element $\bar r \in \mathcal{R}_{(0,1)}$ if $M_z \not= \emptyset$.
Then, there exists a positive rational number $r_u$ such that $r_u < \bar r$ if $M_z \not= \emptyset$, and by the definition of $M_z$, $z \not\in W_{r_u}$. 
The corresponding index $u$ satisfies (i). Thus, we have shown $T(J) = J^*$.

\appsection{An Alternative Proof of Theorem~\ref{thm-valite-P}} \label{appsec-altprfP}
\markboth{\rm Appendix \thesection. Alternative Proof of Theorem~\ref{thm-valite-P}}{\rm Appendix \thesection. Alternative Proof of Theorem~\ref{thm-valite-P}}

In this Appendix, we give an alternative proof of Theorem~\ref{thm-valite-P} about the convergence of value iteration in the nonnegative cost case (P). 
This proof combines Lemma~\ref{lma-P2} with arguments that are motivated by similar arguments used in the book by Meyn~\cite[Chap.\ 9.4]{Mey08} to analyze average cost and total cost problems. Recently, this proof approach has also been used by the first author to extend Theorem~\ref{thm-valite-P} to a class of total cost problems without sign constraints on the one-stage costs~\cite{Yu-vi14}.

It is also worth mentioning that there are two alternative proofs for Whittle's bridging condition, Whittle's original proof \cite{Whit79}, and Hartley's proof \cite{Har80} based on a concavity argument. For Theorem~\ref{thm-valite-P}, in terms of the arguments involved, the proof given in Section~\ref{sec5} can be compared with Hartley's proof, and the alternative proof given below can be compared with Whittle's original proof of the bridging condition.

\mysmallskip
\begin{lem} \label{lma-P-altprf}
{\rm (P)} For a given $x \in S$, let $\pi$ be a policy with $J_\pi(x) < \infty$. Then with $x_0 = x$, we have
$$\lim_{n \to \infty} \, \E^{\pi} \big\{ J^*(x_n) \big\} = 0.$$
\end{lem}
\begin{proof}
For any $n \geq 0$, we have
$$ J_\pi(x) = \E^\pi \left\{ \sum_{k=0}^{n-1} g(x_k, u_k) \right\} + \E^\pi \left\{ \sum_{k=n}^{\infty} g(x_k, u_k) \right\} \geq \E^\pi \left\{ \sum_{k=0}^{n-1} g(x_k, u_k) \right\} + \E^\pi \big\{ J^*(x_n) \big\}.$$
Since $\E^\pi \left\{ \sum_{k=0}^{n-1} g(x_k, u_k) \right\} \uparrow J_\pi(x) < \infty$ as $n \to \infty$, it follows that
$\limsup_{n \to \infty} \, \E^\pi \big\{ J^*(x_n) \big\} = 0$ and hence $\lim_{n \to \infty} \, \E^\pi \big\{ J^*(x_n) \big\} = 0$.
\end{proof}

\mysmallskip
\begin{proof}[Proof of Theorem~\ref{thm-valite-P}]
As discussed immediately after the theorem, it is sufficient to prove part (a). Denote $J^\infty =  \lim_{k \to \infty} T^k(c J^*)$. By Lemma~\ref{lma-P2}, $T (J^\infty) = J^\infty$ and $J^* \leq J^\infty \leq cJ^*.$ We need to show $J^\infty = J^*$. The equality $J^\infty(x) = J^*(x)$ holds certainly for all states $x$ with $J^*(x) = \infty$. So we consider an arbitrary $x \in S$ with $J^*(x) < \infty$, and we show that for an arbitrary $\epsilon > 0$, $J^\infty(x) \leq J^*(x) + \epsilon$.  

Let $\pi=(\mu_0, \mu_1, \ldots)$ be an $\epsilon$-optimal Markov policy for the state $x$. Since $J^\infty =T (J^\infty)$, we have for any $n \geq 1$, with $x_0=x$,
\begin{align}
 J^\infty(x) = T^n (J^\infty)(x) & \leq \big(T_{\mu_0} \circ T_{\mu_1} \circ \cdots T_{\mu_{n-1}} \big) (J^\infty)(x) \notag \\
       & = \E^{\pi} \left\{ \sum_{k=0}^{n-1} g(x_k, u_k) \right\} + \E^\pi \big\{ J^\infty(x_n) \big\} \notag \\
       & \leq J^*(x) + \epsilon + c \, \E^\pi \big\{ J^*(x_n) \big\}, \label{eq-prf-altprfP}
\end{align}
where we used the $\epsilon$-optimality of $\pi$ for $x$ and the fact $J^\infty \leq cJ^*$ in the last inequality. Since $J_\pi(x) \leq J^*(x) + \epsilon < \infty$, we have $\lim_{n \to \infty} \E^\pi \big\{ J^*(x_n) \big\} = 0$ by Lemma~\ref{lma-P-altprf}, so it follows from the inequality~(\ref{eq-prf-altprfP}) that $J^\infty(x) \leq J^*(x) + \epsilon$. Since $\epsilon$ is arbitrary, we obtain $J^\infty(x) \leq J^*(x)$ for all $x$ with $J^*(x) < \infty$. In view of the fact $J^\infty \geq J^*$, equality must hold and this completes the proof. 
\end{proof}

\end{appendices}

\end{document}